\documentclass[12pt]{article}
\usepackage{amssymb, amsmath, amsthm, amscd}
\usepackage[pdftex]{graphics}
\usepackage[utf8]{inputenc}
\usepackage[all,cmtip]{xy}
\usepackage{enumitem}
\usepackage{setspace}
\usepackage{mathdots}

\usepackage[mathscr]{euscript}
\usepackage[colorlinks=true]{hyperref}
\hypersetup{colorlinks   = true}
\hypersetup{linkcolor=blue}
\usepackage{tikz}

\begin{document}

\newtheorem{The}{Theorem}[section]
\newtheorem{Lem}[The]{Lemma}
\newtheorem{Prop}[The]{Proposition}
\newtheorem{Cor}[The]{Corollary}
\newtheorem{Rem}[The]{Remark}
\newtheorem{Obs}[The]{Observation}
\newtheorem{SConj}[The]{Standard Conjecture}
\newtheorem{Titre}[The]{\!\!\!\! }
\newtheorem{Conj}[The]{Conjecture}
\newtheorem{Question}[The]{Question}
\newtheorem{Prob}[The]{Problem}
\newtheorem{Def}[The]{Definition}
\newtheorem{Not}[The]{Notation}
\newtheorem{Claim}[The]{Claim}
\newtheorem{Conc}[The]{Conclusion}
\newtheorem{Ex}[The]{Example}
\newtheorem{Fact}[The]{Fact}
\newtheorem{Formula}[The]{Formula}
\newtheorem{Formulae}[The]{Formulae}
\newtheorem{The-Def}[The]{Theorem and Definition}
\newtheorem{Prop-Def}[The]{Proposition and Definition}
\newtheorem{Lem-Def}[The]{Lemma and Definition}
\newtheorem{Cor-Def}[The]{Corollary and Definition}
\newtheorem{Conc-Def}[The]{Conclusion and Definition}
\newtheorem{Terminology}[The]{Note on terminology}
\newcommand{\C}{\mathbb{C}}
\newcommand{\R}{\mathbb{R}}
\newcommand{\N}{\mathbb{N}}
\newcommand{\Z}{\mathbb{Z}}
\newcommand{\Q}{\mathbb{Q}}
\newcommand{\Proj}{\mathbb{P}}
\newcommand{\Rc}{\mathcal{R}}
\newcommand{\Oc}{\mathcal{O}}
\newcommand{\Vc}{\mathcal{V}}
\newcommand{\Id}{\operatorname{Id}}
\newcommand{\pr}{\operatorname{pr}}
\newcommand{\rk}{\operatorname{rk}}
\newcommand{\del}{\partial}
\newcommand{\delbar}{\bar{\partial}}
\newcommand{\Cdot}{{\raisebox{-0.7ex}[0pt][0pt]{\scalebox{2.0}{$\cdot$}}}}
\newcommand\nilm{\Gamma\backslash G}
\newcommand\frg{{\mathfrak g}}
\newcommand{\fg}{\mathfrak g}
\newcommand{\Oh}{\mathcal{O}}
\newcommand{\Kur}{\operatorname{Kur}}
\newcommand\gc{\frg_\mathbb{C}}

\begin{center}

{\Large\bf Twisted Adiabatic Limit for Complex Structures}

\end{center}

\begin{center}

{\large Dan Popovici}

\end{center}

\vspace{1ex}

\noindent{\small{\bf Abstract.} Given a complex manifold $X$ and a smooth positive function $\eta$ thereon, we perturb the standard differential operator $d=\partial + \bar\partial$ acting on differential forms to a first-order differential operator $D_\eta$ whose principal part is $\eta\partial + \bar\partial$. The role of the zero-th order part is to force the integrability property $D_\eta^2=0$ that leads to a cohomology isomorphic to the de Rham cohomology of $X$, while the components of types $(0,\,1)$ and $(1,\,0)$ of $D_\eta$ induce cohomologies isomorphic to the Dolbeault and conjugate-Dolbeault cohomologies. We compute Bochner-Kodaira-Nakano-type formulae for the Laplacians induced by these operators and a given Hermitian metric on $X$. The computations throw up curvature-like operators of order one that can be made (semi-)positive under appropriate assumptions on the function $\eta$. As applications, we obtain vanishing results for certain harmonic spaces on complete, non-compact, manifolds and for the Dolbeault cohomology of compact complex manifolds that carry certain types of functions $\eta$. This study continues and generalises the one of the operators $d_h=h\partial + \bar\partial$ that we introduced and investigated recently for a positive constant $h$ that was then let to converge to $0$ and, more generally, for constants $h\in\C$. The operators $d_h$ had, in turn, been adapted to complex structures from the well-known adiabatic limit construction for Riemannian foliations. Allowing now for possibly non-constant functions $\eta$ creates positivity in the curvature-like operator that stands one in good stead for various kinds of applications.}

\vspace{1ex}

\section{Introduction}\label{section:introd} Let $\pi:{\cal X}\longrightarrow B$ be a holomorphic family of compact complex manifolds $X_t:=\pi^{-1}(t)\subset{\cal X}$ with $t$ varying in a small open ball $B$ about the origin in some $\C^N$. This means that ${\cal X}$ is a complex manifold and $\pi$ is a proper holomorphic submersion.

It is standard that the degeneration at the first page of the Fr\"olicher spectral sequence (FSS) is a deformation-open property, namely that $E_1(X_0) = E_\infty(X_0)$ implies the analogous property $E_1(X_t) = E_\infty(X_t)$ on the nearby fibres $X_t$ when $t$ is sufficiently close to $0$. This follows at once from the Kodaira-Spencer theory, especially from the upper semi-continuity of the functions $B\ni t\mapsto h^{p,\,q}_{\bar\partial}(t)$ whose values are the Hodge numbers (i.e. the $\C$-vector space dimensions of the Dolbeault cohomology groups $H^{p,\,q}_{\bar\partial}(X_t,\,\C)$ ) of the fibres $X_t$.

However, this is no longer true for the higher pages of the FSS where the analogous numbers $e_r^{p,\,q}(t):=\mbox{dim}_\C E_r^{p,\,q}(X_t)$ need not vary semi-continuously with $t\in B$ when $r\geq 2$. For an example of such a pathological behaviour, see [COUV16, Example 4.8. and Corollary 4.9.] where Ceballos, Otal, Ugarte and Villacampa consider a nilmanifold $M$ with underlying Lie algebra $\frak{h}_{15}$, endowed with a family of invariant complex structures $J_t$, and observe, as a consequence of their classification of invariant complex structures on $6$-nilmanifolds, that the numbers $e_2^{1,\,1}(t)$ and $e_3^{1,\,1}(t)$ are not upper semi-continuous functions of $t$, while the numbers $e_2^{0,\,2}(t)$ and $e_3^{0,\,2}(t)$ are not lower semi-continuous functions of $t$.

Nevertheless, there are quite a few examples of classes of compact complex manifolds whose Fr\"olicher spectral sequence degenerates at $E_2$ (though not at $E_1$) and for which this property persists in their small deformations. Our intuition is that this ought to be due to some geometric property of the central fibre that forces the FSS to behave well under small deformations, so we propose the following issue to ponder.

\begin{Prob}\label{Prob:E_r_def} Let $r\geq 2$ be an integer. Find a geometric property $(P)$ that certain compact complex manifolds $X$ satisfy such that, whenever $X_0$ has property $(P)$ and has its Fr\"olicher spectral sequence degenerate at $E_r$, the Fr\"olicher spectral sequence of every fibre $X_t$ with $t\in B$ close enough to $0$ degenerates again at $E_r$.

\end{Prob}

We refrain from speculating on the nature (metric? cohomological?) of any such property $(P)$, but we stress the need for it to be readily verifiable on concrete examples of manifolds. An analytic such property, in terms of the decay rate to $0$ of the small eigenvalues of certain Laplacians, was given in [Pop19], but that result seems hard to apply in explicit families of manifolds. One of the goals of the present paper is to generalise the main construction of [Pop19] (whose key points we recall in outline in $\S$\ref{subsection:introd_adiabatic_h_reminder} for the reader's convenience) in order to make it more flexible.

A solution to Problem \ref{Prob:E_r_def} is expected to play a central role in various contexts, including in the further development of a non-K\"ahler mirror symmetry theory that started in [Pop18] and has continued with [PSU22], [PSU21a] and [PSU21b].

\subsection{Reminder of the main facts from [Pop19]}\label{subsection:introd_adiabatic_h_reminder} Let $X$ be a complex manifold with $\mbox{dim}_\C X=n$. 

By adapting to the setting of complex structures the {\it adiabatic limit} construction for Riemannian foliations initiated and subsequently studied in e.g. [Wit85] and [MM90], we defined in [Pop19], for constants $h>0$, the first-order differential operators $$d_h:=h\partial + \bar\partial:C^\infty_k(X,\,\C)\longrightarrow C^\infty_{k+1}(X,\,\C),  \hspace{5ex} k\in\{0,\dots , 2n\},$$ and the zero-th order pointwise operators $$\theta_h:\Lambda^{p,\,q}T^\star X\longrightarrow\Lambda^{p,\,q}T^\star X, \hspace{5ex} u\longmapsto\theta_hu:=h^p\,u,$$ that are then extended by linearity to operators $\theta_h:\Lambda^k T^\star X\longrightarrow\Lambda^k T^\star X$ for every $k$.

The equality $d_h=\theta_h\,d\,\theta_h^{-1}$ implies the {\it integrability} property $d_h^2=0$ and the vector-space {\it isomorphisms} $$H^k_{dR}(X,\,\C)\stackrel{\simeq}{\longrightarrow}H^k_{d_h}(X,\,\C), \hspace{5ex} \{u\}_{dR}\longmapsto\{\theta_hu\}_{d_h},$$ for $k\in\{0,\dots , 2n\}$, between the de Rham cohomology groups of $X$ and those of $d_h$-cohomology defined as $\ker d_h/\mbox{Im}\,d_h$.

Fix now a Hermitian metric $\omega$ on $X$. The pointwise inner product $\langle\,\,,\,\,\rangle_\omega$ induced by $\omega$ on the differential forms on $X$ can be rescaled in the following way on $(p,\,q)$-forms: $$\langle u,\,v\rangle_{\omega_h}:= h^{2p}\,\langle u,\,v\rangle_\omega, \hspace{5ex} h>0;\hspace{1ex} u,v\in\Lambda^{p,\,q}T^\star X,$$ for every bidegree $(p,\,q)$. This rescaling defines a new Hermitian metric $$\omega_h = \frac{1}{h^2}\,\omega, \hspace{5ex} h>0,$$ on the holomorphic tangent bundle $T^{1,\,0}X$ of $X$, or equivalently, a rescaled $C^\infty$ positive definite $(1,\,1)$-form $\omega_h=h^{-2}\,\omega$ on $X$. This, in turn, induces a $C^\infty$ positive definite volume form $$dV_{\omega_h}:=\frac{\omega_h^n}{n!} = \frac{1}{h^{2n}}\,\frac{\omega^n}{n!} =  \frac{1}{h^{2n}}\,dV_\omega$$ on $X$ which, together with the pointwise inner product $\langle\,\,,\,\,\rangle_{\omega_h}$, defines an $L^2$-inner product $$\langle\langle u,\,v\rangle\rangle_{\omega_h}:=\int\limits_X\langle u,\,v\rangle_{\omega_h}\,dV_{\omega_h} = \frac{1}{h^{2n}}\,\langle\langle\theta_h u,\,\theta_h v\rangle\rangle_\omega$$ for all forms $u,v\in C^\infty_{p,\,q}(X,\,\C)$ and all bidegrees $(p,\,q)$.

The two rescalings (of the operator $d$ and of the metric $\omega$) lead to two different Laplace-type operators $\Delta_h,\,\Delta_{\omega_h}:C^\infty_k(X,\,\C)\longrightarrow C^\infty_k(X,\,\C)$ defined respectively by $$\Delta_h:=d_hd_h^\star + d_h^\star d_h  \hspace{5ex}\mbox{and}\hspace{5ex} \Delta_{\omega_h}:=dd_{\omega_h}^\star + d_{\omega_h}^\star d$$ on the smooth $k$-forms for every $k\in\{0,\dots , 2n\}$. By $d_h^\star$ we mean the formal adjoint of $d_h$ w.r.t. the $L^2$-inner product induced by the orginal metric $\omega$, while $d_{\omega_h}^\star$ stands for the formal adjoint of the original operator $d$ w.r.t. the $L^2$-inner product induced by the rescaled metric $\omega_h$.

The second-order differential operators $\Delta_h$ and $\Delta_{\omega_h}$ are elliptic, self-adjoint and non-negative, they have the same principal part and are related by the formula \begin{eqnarray}\label{eqn:link_Delta_h-Delta_omega-h}\Delta_h = \theta_h\Delta_{\omega_h}\theta_h^{-1},  \hspace{5ex} h>0.\end{eqnarray}

This formula implies that $\Delta_h$ and $\Delta_{\omega_h}$ have the same spectrum and their respective eigenspaces $E_{\Delta_{\omega_h}}(\lambda)$ and $E_{\Delta_h}(\lambda)$ are obtained from each other via the rescaling isometry $\theta_h$: $$\theta_h\bigg(E_{\Delta_{\omega_h}}(\lambda)\bigg) = E_{\Delta_h}(\lambda) \hspace{5ex}\mbox{for every}\hspace{2ex}\lambda\in\mbox{Spec}(\Delta_h) = \mbox{Spec}(\Delta_{\omega_h}).$$

The main result of [Pop19] expressed, for every positive integer $r$ and every degree $k$, the dimension of the $\C$-vector space $E_r^k(X):=\oplus_{p+q=k}E_r^{p,\,q}(X)$ (representing the direct sum of the spaces of total degree $k$ featuring on the $r^{th}$ page of the Fr\"olicher spectral sequence of $X$) in terms of the number of eigenvalues, counted with multiplicities, of any of the operators $\Delta_h$ and $\Delta_{\omega_h}$ that decay sufficiently fast to $0$ (at a rate depending on $r$) as $h\downarrow 0$.

\begin{The}([Pop19, Theorem 1.3.])\label{The:Pop_19_main} Let $(X,\,\omega)$ be a compact Hermitian manifold with $\mbox{dim}_\C X=n$. For every $r\in\N^\star$ and every $k=0,\dots , 2n$, the following identity holds:

\begin{equation}\label{eqn:main_identity}\mbox{dim}_\C E_r^k(X) = \sharp\bigg\{i\,\mid\,\lambda_i^k(h)\in O(h^{2r}) \hspace{2ex} \mbox{as} \hspace{1ex} h\downarrow 0\bigg\},\end{equation}

\noindent where $0\leq \lambda_1^k(h)\leq \lambda_2^k(h)\leq\dots\leq\lambda_i^k(h)\leq\dots$ are the eigenvalues, counted with multiplicities, of the rescaled Laplacian $\Delta_h : C^\infty_k(X,\,\C)\longrightarrow C^\infty_k(X,\,\C)$ ($=$ those of $\Delta_{\omega_h}: C^\infty_k(X,\,\C)\longrightarrow C^\infty_k(X,\,\C)$) acting on $k$-forms. As usual, $\sharp$ stands for the cardinal of a set.

\end{The}

Now, $0$ is always an eigenvalue of $\Delta_h : C^\infty_k(X,\,\C)\longrightarrow C^\infty_k(X,\,\C)$ of multiplicity exactly equal to the {\it $k$-th Betti number} $b_k = b_k(X)$ of $X$ (since $\Delta_h$ and the standard Laplacian $\Delta=dd^\star + d^\star d$ have {\it isomorphic} kernels, each of these two kernels being isomorphic to the corresponding cohomology space, $H_{d_h}^k(X,\,\C)$, resp. $H_{dR}^k(X,\,\C)$, and these cohomology spaces being mutually {\it isomorphic}, as recalled above). Thus, the smallest {\it positive} eigenvalue of $\Delta_h : C^\infty_k(X,\,\C)\longrightarrow C^\infty_k(X,\,\C)$ is $\lambda_{b_k + 1}^k(h)$. We will denote it by $\delta_h^{(k)}:=\lambda_{b_k + 1}^k(h)$.

Meanwhile, it is standard that we always have $$b_k(X)\leq\mbox{dim}_\C E_r^k(X)  \hspace{5ex}\mbox{for every}\hspace{1ex} k\in\{0,\dots , 2n\}$$ and that the Fr\"olicher spectral sequence of $X$ {\it degenerates} at a given page $E_r$ if and only if all these inequalities are {\it equalities}, namely if and only if $$b_k(X) = \mbox{dim}_\C E_r^k(X)  \hspace{5ex}\mbox{for every}\hspace{1ex} k\in\{0,\dots , 2n\}.$$

Therefore, the above Theorem \ref{The:Pop_19_main} equates the degeneration at $E_r$ of the FSS of $X$ to the fact that $0$ is the only eigenvalue of $\Delta_h$ in every degree $k$ that decays to $0$ at least as fast as $Const\cdot h^{2r}$ as $h\downarrow 0$. In other words, we have the following

\begin{Prop}([Pop19, Proposition 5.3.])\label{Prop:Pop_19_main-conseq} Let $(X,\,\omega)$ be a compact complex Hermitian manifold with $\mbox{dim}_\C X=n$. For every constant $h>0$, let $\delta_h^{(k)}>0$ be the smallest positive eigenvalue of $\Delta_h:C^\infty_k(X,\,\C)\longrightarrow C^\infty_k(X,\,\C)$.

Then, for every $r\in\N^\star$, the Fr\"olicher spectral sequence of $X$ degenerates at $E_r$ if and only if $$\limsup\limits_{h\rightarrow 0}\frac{\delta_h^{(k)}}{h^{2r}} = +\infty, \hspace{3ex} \mbox{for all}\hspace{2ex} k\in\{1,\dots , n\}.$$ 

\end{Prop}

In the context of a holomorphic family $(X_t)_{t\in B}$ of compact complex manifolds on which a $C^\infty$ family $(\omega_t)_{t\in B}$ of Hermitian metrics on the fibres has been fixed, Proposition \ref{Prop:Pop_19_main-conseq} shows that, for some degree $k\in\{1,\dots , n\}$, the decay rate to $0$ of $\delta_{h}^{(k)}>0$ can accelerate when one moves from $X_0$ to the nearby fibres $X_t$ if $E_r(X_0) = E_\infty(X_0)$ but $E_r(X_t)\neq E_\infty(X_t)$ for $t\neq 0$.

\subsection{Constructions and results obtained in this paper}\label{subsection:introd_constructions-results_this} If one aims at solving Problem \ref{Prob:E_r_def}, the takeaway from Proposition \ref{Prop:Pop_19_main-conseq} is that one needs to control the bottom of the positive part of the spectrum of $\Delta_h$, possibly in terms of some (curvature-like) positive quantity whose dependence on the fibre $X_t$ should be at least continuous.

On the other hand, Bochner-Kodaira-Nakano-type (in-)equalities involving the operator $\Delta_h$ cannot produce the needed positivity when $h$ remains constant. (See e.g. [BP18].) It is this quest for positivity that would translate into information on the smallest positive eigenvalues of Laplacians like $\Delta_h$ that motivates the main thrust of this paper (and hopefully of its future sequels): transforming the constant $h$ into a $C^\infty$ function whose derivatives would produce positive curvature-like quantities.

\vspace{2ex}

Let $X$ be a connected complex manifold with $\mbox{dim}_\C X=n$. Fix a $C^\infty$ real-valued function $\eta$ on $X$ such that $\eta>0$ or $\eta < 0$ at every point of $X$.

The immediate analogue of $d_h$ when the constant $h$ has been replaced by $\eta$ is $d_\eta:=\eta\partial + \bar\partial$. However, if $\eta$ is not constant, $d_\eta$ need not be integrable (in the sense that it need not square to $0$), so it need not define a cohomology. The way round this takes us to first defining the pointwise analogue $\theta_\eta$ (cf. (\ref{eqn:theta_eta_def})) of $\theta_h$ and then defining $D_\eta: = \theta_\eta d\theta_\eta^{-1}$ in every degree $k$ (cf. Definition \ref{Def:D_eta}) by the analogue of the formula $d_h = \theta_h d\theta_h^{-1}$ satisfied by $d_h$.

If we fix a Hermitian metric $\omega$ on $X$, we consider the formal adjoint $D_\eta^\star$ of $D_\eta$ w.r.t. the $L^2$-inner product induced by $\omega$ and then the associated twisted Laplacian $$\Delta_\eta = D_\eta D_\eta^\star + D_\eta^\star D_\eta:C^\infty_k(X,\,\C)\longrightarrow C^\infty_k(X,\,\C)$$ in every degree $k$. (See Definition \ref{Def:Delta_eta}.)

The twisted metric $\omega_\eta$ is defined analogously to the $\omega_h$ of [Pop19], but the definition of the twisted Laplacian $\Delta_{\omega_\eta}$ is not as straightforward as in the case where $\eta$ is constant. It is necessary to first compute an operator $T^\star_\eta : C^\infty_k(X,\,\C)\longrightarrow C^\infty_{k-1}(X,\,\C)$ such that $D_\eta^\star = \theta_\eta T^\star_\eta\theta_\eta^{-1}$ in every degree $k$ (cf. Proposition \ref{Prop:Laplacians-link}) that we then use to define (cf. Definition \ref{Def:Delta_omega-eta}) the twisted Laplacian $$\Delta_{\omega_\eta} = d T_\eta^\star + T_\eta^\star d$$ which is then shown to relate to the first twisted Laplacian $\Delta_\eta$ via the identity (cf. Corollary \ref{Cor:Laplacians-link}) $$\Delta_\eta = \theta_\eta\Delta_{\omega_\eta}\theta_\eta^{-1}$$ analogous to (\ref{eqn:link_Delta_h-Delta_omega-h}) of the constant $\eta$ case. This then implies that $\Delta_\eta$ and $\Delta_{\omega_\eta}$ have the same spectrum and that their respective eigenspaces are isomorphic via $\theta_\eta$ as in the case where $\eta$ is constant.

\vspace{2ex}

In $\S$\ref{section:twisted-commutation_D-eta}, we compute the formal adjoint of $D_\eta$ and observe that it depends on $\overline{D_{-\eta}}$, a fact that leads to {\bf $\eta$-twisted commutation relations} (Proposition \ref{Prop:eta-twisted_commutation-relations}) and eventually to two {\bf Bochner-Kodaira-Nakano-type identities} relating the Laplacians $\Delta_\eta$ and $\overline{\Delta}_{-\eta}$ (cf. Proposition \ref{Prop:eta-BKN_rough} giving the {\it rough} version of the identity), respectively the Laplacians $\Delta_\eta$ and $[\overline{D}_{-\eta} + \overline\tau_{-\eta},\,\overline{D}_{-\eta}^\star + \overline\tau_{-\eta}^\star]$ (cf. Theorem \ref{The:eta-BKN}), where $\tau_\eta$ is the zero-th order operator defined as $[\Lambda,\,D_\eta\omega\wedge\cdot\,]$. The latter identity is the {\it refined} version of the former. It absorbs some of the torsion terms into the Laplacian to which $\Delta_\eta$ is compared and that can be neglected in applications due to its non-negativity. Both versions of the $\eta$-BKN identity expressing $\Delta_\eta$ in terms of another Laplacian throw up the first-order {\it curvature operator} $i\,[[D_\eta,\,\overline{D}_{-\eta}],\,\Lambda]$ that is then computed in Proposition \ref{Prop:curvature_computation}.


\vspace{2ex}

In $\S$\ref{section:twisted-commutation_1-0_0-1}, we run analogous computations starting from the $(0,\,1)$-part $D_\eta^{0,\,1}$ of $D_\eta = D_\eta^{1,\,0} + D_\eta^{0,\,1}$ and its conjugate. This setting has the advantage that $D_\eta^{0,\,1}$ and $D_\eta^{1,\,0}$ anti-commute (cf. (\ref{eqn:D_eta_1-0_0-1_prop})), unlike $D_\eta$ and $\overline{D_{-\eta}}$ of the previous $\S$\ref{section:twisted-commutation_D-eta}. Meanwhile, $D_\eta^{0,\,1}$ and $D_\eta^{1,\,0}$ are integrable and define cohomologies isomorphic to the Dolbeault, respectively conjugate-Dolbeault, cohomologies of $X$ (see (\ref{eqn:cohomologies_1-0_0-1_isom})).

We consider the Laplacians \begin{eqnarray*}\Delta''_\eta:=[D_\eta^{0,\,1},\,(D_\eta^{0,\,1})^\star]  \hspace{5ex}\mbox{and}\hspace{5ex}  \Delta'_\eta:=[\overline{D_\eta^{0,\,1}},\,\overline{D_\eta^{0,\,1}}^{\,\star}]\end{eqnarray*} that we then relate to each other in what we call the {\bf rough $\eta$-Bochner-Kodaira-Nakano ($\eta$-BKN) identity} (cf. Proposition \ref{Prop:eta-BKN_rough_1-0_0-1}). Further commutation relations involving torsion terms (cf. Lemma \ref{Lem:preliminary_BKN-refined_1-0_0-1}) lead to the {\bf refined $\eta$-Bochner-Kodaira-Nakano ($\eta$-BKN) identity} for this setting (cf. Theorem \ref{The:eta-BKN_1-0_0-1}). In bidegree $(p,\,q)$, the zero-th order part of the first-order {\it curvature operator} $i\,[[D^{0,\,1}_\eta,\,\overline{D^{0,\,1}_\eta}],\,\Lambda]$ turns out to be (cf. Proposition \ref{Prop:curvature_computation_1-0_0-1}) the operator \begin{eqnarray*}(p-q)\,\bigg[\bigg(\frac{2}{\eta^2}\,i\partial\eta\wedge\bar\partial\eta - \frac{1}{\eta}\,i\partial\bar\partial\eta\bigg)\wedge\cdot\,,\,\Lambda\bigg].\end{eqnarray*}

We give two applications of our {\bf $\eta$-BKN identities} based on the shape of the curvature operator.

\vspace{2ex}

In $\S$\ref{section:vanishing_harmonic}, we deal with {\it non-compact complete} manifolds $X$ on which we assume the existence of a $C^\infty$ function $\eta:X\longrightarrow(0,\,\infty)$ that induces a {\bf positive definite} {\it curvature form} \begin{eqnarray*}\gamma_\eta:=\frac{2}{\eta^2}\,i\partial\eta\wedge\bar\partial\eta - \frac{1}{\eta}\,i\partial\bar\partial\eta >0.\end{eqnarray*}

This key positivity assumption, which amounts to requiring the smooth real $(1,\,1)$-form $\gamma_\eta$ to define a Hermitian metric on $X$, is impossible on compact manifolds due to the maximum principle, but is quite natural in the non-compact setting given the fact that it constitutes the zero-th order term of the operator that plays the role of the curvature in our $\eta$-BKN identities. The result we obtain can be loosely formulated as follows (see Theorem \ref{The:vanishing_harmonic} for the precise statement).

\begin{The}\label{The:introd_vanishing_harmonic} Under the extra assumptions that the Hermitian metric $\gamma_\eta$ on $X$ is {\bf complete} and the pointwise $\gamma_\eta$-norm $|\partial\eta| = |\partial\eta|_{\gamma_\eta}$ of the $(1,\,0)$-form $\partial\eta$ is {\bf small} relative to $\eta$, there exist no non-zero $\Delta''_\eta$-harmonic $L^2_{\gamma_\eta}$-forms of bidegree $(p,\,q)$ on $X$ whenever either $\bigg(p>q \hspace{1ex}\mbox{and}\hspace{1ex} p+q\geq n+1\bigg)$ or $\bigg(p<q \hspace{1ex}\mbox{and}\hspace{1ex} p+q\leq n-1\bigg)$.

\end{The}

This resut (Theorem \ref{The:vanishing_harmonic}) can be compared to Gromov's Main Theorem $2.5.$ in [Gro91]. They are both vanishing theorems for certain spaces of $L^2$ harmonic forms on a complete manifold carrying a K\"ahler metric. The particular shape of our complete metric $\gamma_\eta$ can be viewed as a (non-standard) exactness property of $\gamma_\eta$, while Gromov imposed the exactness condition requiring the metric to be $d(bounded)$. The vanishing conclusion in Gromov's theorem applies to all the bidegrees $(p,\,q)$ with $p+q\neq n$ (i.e. outside the middle degree), while our conclusion, though applicable to fewer bidegrees, also avoids the middle degree. Moreover, much as in Gromov's case, our Theorem \ref{The:vanishing_harmonic} is mainly intended for applications where the manifold $X$ is the universal cover of a compact complex manifold $Y$ and the complete metric $\gamma_\eta$ on $X$ is the pullback of some Hermitian metric on $Y$.

We hope that Theorem \ref{The:vanishing_harmonic} will have a role to play in the further development of a (possibly non-K\"ahler) hyperbolicity theory linking the existence of certain types of special Hermitian metrics having a certain exactness property on the universal cover of the compact complex manifold on which they are defined to the non-existence of entire holomorphic maps with a relatively small growth from some $\C^p$ (with $p$ possibly $>1$) into the given manifold. The K\"ahler case with $p=1$ was treated in [Gro91], while the more general Hermitian case with $p$ allowed to be $>1$ was initiated in [MP21], [MP22] and [KP23].

\vspace{2ex}

In $\S$\ref{section:vanishing_cohomology_compact}, we deal with {\it compact} manifolds $X$ on which we assume the existence of a $C^\infty$ function $\eta:X\longrightarrow(0,\,\infty)$ whose induced first-order {\it curvature operator} $F_\eta:=i\,[[D^{0,\,1}_\eta,\,\overline{D^{0,\,1}_\eta}],\,\Lambda]$ is assumed to be sufficiently positive. Exploiting the fact that $\ker\Delta''_\eta$ is isomorphic to the Dolbeault cohomology group $H^{p,\,q}_{\bar\partial}(X,\,\C)$ in every bidegree $(p,\,q)$, our {\bf $\eta$-BKN identities} can be applied to give the following result (see Theorem \ref{The:vanishing_cohomology_compact} for the precise statement).

\begin{The}\label{The:introd_vanishing_cohomology_compact} Let $(X,\,\omega)$ be a compact complex Hermitian manifold with $\mbox{dim}_\C X=n$. If there exists a $C^\infty$ function $\eta:X\longrightarrow(0,\,\infty)$ such that $\sup_X|\bar\partial\omega|_\omega$ and $\sup_X\frac{|\partial\eta|_\omega}{\eta}$ are sufficiently small, then for every bidegree $(p,\,q)$ such that $$F_\eta:=i\,\bigg[[D^{0,\,1}_\eta,\,\overline{D^{0,\,1}_\eta}],\,\Lambda\bigg]>C(\eta)\,\mbox{Id}$$ in bidegree $(p,\,q)$ for an appropriate constant $C(\eta)>0$ depending on $\eta$, we have $H^{p,\,q}_{\bar\partial}(X,\,\C) = \{0\}$.

\end{The}

We hope further applications of these results will be obtained in future work, either to investigate Problem \ref{Prob:E_r_def} or to study cohomological and metric properties of complex manifolds, for example by choosing particular types of functions $\eta$ supported by certain such manifolds.

\section{The twisted operators and metrics}\label{section:twisted_def}

Let $X$ be an $n$-dimensional complex manifold, where $n\geq 2$. 

\subsection{The twisted operator $D_\eta$}\label{subsection:twisted-d} With every $C^\infty$ function $\eta$ on $X$ such that $\eta>0$ or $\eta<0$, we associate the bijective linear operators: \begin{eqnarray}\label{eqn:theta_eta_def}\nonumber\theta_{\eta(x)}:\Lambda^{p,\,q}T^\star_xX & \longrightarrow & \Lambda^{p,\,q}T^\star_xX\\
  u(x) & \longmapsto & \eta^p(x)\,u(x)\end{eqnarray} defined pointwise on the $(p,\,q)$-forms of $X$ at every point $x\in X$ and in every bidegree $(p,\,q)$. Thus, we get an automorphism $\theta_\eta:\Lambda^{p,\,q}T^\star X \longrightarrow \Lambda^{p,\,q}T^\star X$ of the vector bundle of $(p,\,q)$-forms on $X$ defined by $\theta_\eta(u) = \eta^p u$. 

We then extend $\theta_\eta$ by linearity to an automorphism \begin{eqnarray*}\theta_\eta:\Lambda^kT^\star X & \longrightarrow & \Lambda^kT^\star X \\
  u=\sum\limits_{p+q=k}u^{p,\,q} & \longmapsto & \sum\limits_{p+q=k}\theta_\eta(u^{p,\,q})=\sum\limits_{p+q=k}\eta^p\,u^{p,\,q}\end{eqnarray*} of the vector bundle of $k$-forms on $X$ for every $k\in\{0,\dots , 2n\}$.

The conjugate of $\theta_\eta$ is defined by the condition $\overline\theta_\eta(u) = \overline{\theta_\eta(\bar{u})}$ for every form $u$. 

\begin{Obs}\label{Obs:theta-eta_inverse} The inverse of $\theta_\eta$ is given by the formula: \begin{eqnarray}\label{eqn:theta-eta_inverse}\theta_\eta^{-1} = \theta_{\eta^{-1}},\end{eqnarray} where $\eta^{-1} = 1/\eta$, while the conjugate of $\theta_\eta$ is given  by the formula: \begin{eqnarray}\label{eqn:theta-eta_conjugate}\overline\theta_\eta(u^{p,\,q}) = \eta^q\,u^{p,\,q}\end{eqnarray} for every $(p,\,q)$-form $u^{p,\,q}$.

  In particular, $\theta_\eta\neq\overline\theta_\eta$ in any bidegree $(p,\,q)$ with $p\neq q$.

\end{Obs}

\noindent {\it Proof.} Immediate verification.  \hfill $\Box$

\begin{Def}\label{Def:D_eta} For any $C^\infty$ function $\eta$ on $X$ such that $\eta>0$ or $\eta<0$ and any $k\in\{0,\dots , 2n\}$, let $D_\eta:C^\infty_k(X,\,\C)\longrightarrow C^\infty_{k+1}(X,\,\C)$ be the differential operator defined on the $C^\infty$ $k$-forms on $X$ by \begin{eqnarray}\label{eqn:D_eta}D_\eta = \theta_\eta d\theta_\eta^{-1}.\end{eqnarray}

\end{Def}

Since $D_\eta^2 = \theta_\eta d^2\theta_\eta^{-1} =0$, the operator $D_\eta$ induces a cohomology space in every degree $k\in\{0,\dots , 2n\}$ on the $n$-dimensional manifold $X$ by \begin{eqnarray}\label{eqn:D_h-cohomology_def}H^k_{D_\eta}(X,\,\C):=\frac{\ker\bigg(D_\eta:C^\infty_k(X,\,\C)\longrightarrow C^\infty_{k+1}(X,\,\C)\bigg)}{\mbox{Im}\,\bigg(D_\eta:C^\infty_{k-1}(X,\,\C)\longrightarrow C^\infty_k(X,\,\C)\bigg)}.\end{eqnarray} 

An immediate observation is that the $D_\eta$-cohomology is canonically (i.e. in a way depending only on the complex structure of $X$ and the function $\eta$, but independent of any metric on $X$) isomorphic to the De Rham cohomology of $X$.

\begin{Prop}\label{Prop:cohomologies_isom} Let $X$ be an $n$-dimensional complex manifold. For any $C^\infty$ function $\eta$ on $X$ such that $\eta>0$ or $\eta<0$ and any $k\in\{0,\dots , 2n\}$, the linear map \begin{eqnarray}\label{eqn:cohomologies_isom}\theta_\eta:H^k_{DR}(X,\,\C)\longrightarrow H^k_{D_\eta}(X,\,\C), \hspace{5ex} \{u\}_{DR}\longmapsto\{\theta_\eta u\}_{D_\eta},\end{eqnarray} is {\bf well defined} and an {\bf isomorphism}.

\end{Prop}

\noindent {\it Proof.} To prove well-definedness, we need to prove the inclusions: \begin{eqnarray*}\theta_\eta(\ker d)\subset\ker D_\eta  \hspace{5ex}\mbox{and}\hspace{5ex} \theta_\eta(\mbox{Im}\, d)\subset\mbox{Im}\, D_\eta.\end{eqnarray*} These inclusions follow, respectively, from the equivalences: \begin{eqnarray*}D_\eta u = 0 \iff d(\theta_\eta^{-1}u) = 0 \hspace{5ex}\mbox{and}\hspace{5ex} u=D_\eta v \iff \theta_\eta^{-1}u = d(\theta_\eta^{-1}v).\end{eqnarray*} These equivalences actually amount to the above inclusions being equalities: \begin{eqnarray*}\theta_\eta(\ker d) = \ker D_\eta  \hspace{5ex}\mbox{and}\hspace{5ex} \theta_\eta(\mbox{Im}\, d) = \mbox{Im}\, D_\eta.\end{eqnarray*} 

Thanks to these equalities of vector spaces and to $\theta_\eta:C^\infty_k(X,\,\C)\longrightarrow C^\infty_k(X,\,\C)$ being an isomorphism at the level of differential forms, the well-defined linear map induced by $\theta_\eta$ in cohomology is an isomorphism.  \hfill $\Box$

\begin{Prop-Def}\label{Prop-Def:D_eta_explicit} For every $C^\infty$ function $\eta$ on $X$ such that $\eta>0$ or $\eta<0$ and every bidegree $(p,\,q)$, the operator $D_\eta:C^\infty_{p,\,q}(X,\,\C)\longrightarrow C^\infty_{p+q+1}(X,\,\C)$ of Definition \ref{Def:D_eta} arises explicitly as \begin{eqnarray}\label{eqn:D_eta_explicit}D_\eta = \eta\partial + \bar\partial - \frac{p}{\eta}\bigg(\eta\partial + \bar\partial\bigg)\eta\wedge\cdot = d_\eta - \frac{p}{\eta}(d_\eta\eta)\wedge\cdot,\end{eqnarray} where the operator $d_\eta:C^\infty_k(X,\,\C)\longrightarrow C^\infty_{k+1}(X,\,\C)$ is defined in any degree $k$ by $d_\eta=\eta\partial + \bar\partial$.

We let $D_\eta^{1,\,0}:C^\infty_{p,\,q}(X,\,\C)\longrightarrow C^\infty_{p+1,\,q}(X,\,\C)$ and $D_\eta^{0,\,1}:C^\infty_{p,\,q}(X,\,\C)\longrightarrow C^\infty_{p,\,q+1}(X,\,\C)$ be the differential operators: \begin{eqnarray}\label{eqn:D_eta_1-0-1_def}D_\eta^{1,\,0} = \eta\partial - p\,\partial\eta\wedge\cdot  \hspace{3ex} \mbox{and} \hspace{3ex} D_\eta^{0,\,1} = \bar\partial - \frac{p}{\eta}\,\bar\partial\eta\wedge\cdot\end{eqnarray} that are the components of bidegrees $(1,\,0)$ and $(0,\,1)$ of $D_\eta$ in the decomposition $D_\eta = D_\eta^{1,\,0} + D_\eta^{0,\,1}$.

\end{Prop-Def}

\noindent {\it Proof.} For every form $u\in C^\infty_{p,\,q}(X,\,\C)$, we have: \begin{eqnarray*}D_\eta u & = & (\theta_\eta d\theta_\eta^{-1})(u) = \theta_\eta d\bigg(\frac{1}{\eta^p}\,u\bigg) = \theta_\eta\bigg(\frac{1}{\eta^p}\,du - \frac{p}{\eta^{p+1}}\,d\eta\wedge u\bigg) \\
& = & \frac{1}{\eta^p}\,\theta_\eta(\partial u) + \frac{1}{\eta^p}\,\theta_\eta(\bar\partial u) - \frac{p}{\eta^{p+1}}\,\theta_\eta(\partial\eta\wedge u) - \frac{p}{\eta^{p+1}}\,\theta_\eta(\bar\partial\eta\wedge u) \\
  & = & \eta\,\partial u + \bar\partial u - p\,\partial\eta\wedge u - \frac{p}{\eta}\,\bar\partial\eta\wedge u,\end{eqnarray*} where the last equality follows from the definition of $\theta_\eta$ and the fact that the forms $\partial u$ and $\partial\eta\wedge u$ are of bidegree $(p+1,\,q)$, while the forms $\bar\partial u$ and $\bar\partial\eta\wedge u$ are of bidegree $(p,\,q+1)$. This proves (\ref{eqn:D_eta_explicit}). \hfill $\Box$

\vspace{2ex}

As with any operator, we define the {\it conjugate} $\overline{D}_\eta$ of $D_\eta$ by requiring the equality \begin{eqnarray*}\overline{D_\eta u} = \overline{D}_\eta(\bar{u})\end{eqnarray*} to hold for every form $u$ on $X$. We now compute $\overline{D}_{-\eta}$ as it will be needed later on. If $u^{p,\,q}$ is a $(p,\,q)$-form, by conjugating the expression of $D_\eta u^{p,\,q}$ we get: $\overline{D_\eta u^{p,\,q}} = \eta\,\bar\partial\overline{u^{p,\,q}} + \partial\overline{u^{p,\,q}} - p\,\bar\partial\eta\wedge\overline{u^{p,\,q}} - \frac{p}{\eta}\,\partial\eta\wedge\overline{u^{p,\,q}}.$ Replacing $\eta$ with $-\eta$, requiring $\overline{D_{-\eta}u^{p,\,q}} = \overline{D}_{-\eta}(\overline{u^{p,\,q}})$, using the fact that $\overline{u^{p,\,q}}$ is of type $(q,\,p)$ and then permuting $p$ and $q$, we get \begin{eqnarray}\label{eqn:D_bar_minus-eta}\overline{D}_{-\eta} = \partial - \eta\bar\partial - \frac{q}{\eta}\,\partial\eta\wedge\cdot + q\,\bar\partial\eta\wedge\cdot \hspace{6ex} \mbox{in bidegree}\hspace{2ex} (p,\,q).\end{eqnarray}

Thus, the coefficients of the zero-th order terms of $\overline{D}_\eta$ and $\overline{D}_{-\eta}$ depend on the anti-holomorphic degree ($q$), while the coefficients of the zero-th order terms of $D_\eta$ and $D_{-\eta}$ depend on the holomorphic degree ($p$), of the bidegree in which they act. Meanwhile, \begin{eqnarray}\label{eqn:d_bar_minus-eta}\overline{d}_{-\eta} = \partial - \eta\bar\partial\end{eqnarray} in every (bi-)degree.

\subsection{The twisted metric $\omega_\eta$}\label{subsection:twisted-omega} Let $\omega$ be a Hermitian metric on $X$. For every $C^\infty$ function $\eta:X\longrightarrow(0,\,+\infty)$, we define the following twisting of the induced pointwise inner product $\langle\,\,,\,\,\rangle_\omega$ on the $(p,\,q)$-forms on $X$: \begin{eqnarray}\label{eqn:eta-inner-product_def}\langle u,\,v\rangle_{\omega_\eta}:=\eta^{2p}\,\langle u,\,v\rangle_\omega = \langle\theta_\eta(u),\,\theta_\eta(v)\rangle_\omega,  \hspace{5ex} u,v\in\Lambda^{p,\,q}T^\star X,\end{eqnarray} in every bidegree $(p,\,q)$. In this way, we get a Hermitian metric $\omega_\eta$ on $X$ induced by the smooth positive definite $(1,\,1)$-form \begin{eqnarray}\label{eqn:eta-twisted-metric_def}\omega_\eta = \frac{1}{\eta^2}\,\omega.\end{eqnarray} 

In particular, the volume forms on $X$ induced by the metrics $\omega_\eta$ and $\omega$ are related by the formula \begin{eqnarray}\label{eqn:volume-forms_eta-twisted}dV_{\omega_\eta} = \frac{1}{\eta^{2n}}\,dV_\omega.\end{eqnarray}

This leads to the $L^2$-inner products induced by $\omega_\eta$ and $\omega$ being related as follows: \begin{eqnarray}\label{eqn:eta-L2-inner-product_def}\langle\langle u,\,v\rangle\rangle_{\omega_\eta} = \int\limits_X\langle u,\,v\rangle_{\omega_\eta}\, dV_{\omega_\eta} = \int\limits_X\frac{1}{\eta^{2(n-p)}}\,\langle u,\,v\rangle_\omega\, dV_\omega,  \hspace{5ex} u,v\in C^\infty_{p,\,q}(X,\,\C).\end{eqnarray}

\subsection{The twisted Laplacians $\Delta_\eta$ and $\Delta_{\omega_\eta}$}\label{subsection:twisted-Laplacians} Suppose that $(X,\,\omega)$ is a complex Hermitian manifold with $\mbox{dim}_\C X=n$. 

\begin{Def}\label{Def:Delta_eta} For every $C^\infty$ function $\eta$ on $X$ such that $\eta>0$ or $\eta<0$ and any $k\in\{0,\dots , 2n\}$, the $D_\eta$-Laplacian $\Delta_\eta:C^\infty_k(X,\,\C)\longrightarrow C^\infty_k(X,\,\C)$ is the differential operator defined by \begin{eqnarray}\label{eqn:Delta_eta}\Delta_\eta = D_\eta D_\eta^\star + D_\eta^\star D_\eta,\end{eqnarray} where $D_\eta^\star$ is the formal adjoint of $D_\eta$ with respect to the $L^2$-inner product $\langle\langle\,\cdot\,,\,\cdot\,\rangle\rangle_\omega$ induced by $\omega$.

\end{Def}

Note that the principal part of $\Delta_\eta$ is the second-order differential operator $\eta^2\,\Delta' + \Delta''$, proving that $\Delta_\eta$ is {\bf elliptic} (since $\Delta'$ and $\Delta''$ are known to be so). If, moreover, $X$ is {\bf compact}, we get, in every degree $k$, the {\bf Hodge isomorphism} \begin{eqnarray}\label{eqn:Hodge_isom_Delta_eta}H^k_{D_\eta}(X,\,\C)\simeq{\cal H}_{\Delta_\eta}^k(X,\,\C):=\ker\bigg(\Delta_\eta:C^\infty_k(X,\,\C)\longrightarrow C^\infty_k(X,\,\C)\bigg)\end{eqnarray} mapping every $D_\eta$-cohomology class to its unique $\Delta_\eta$-harmonic representative.

\vspace{2ex}

We will now compute $D_\eta^\star$ in terms of the formal adjoint $d^\star_{\omega_\eta}$ of $d$ with respect to the twisted metric $\omega_\eta$. Specifically, we will prove the following

\begin{Prop}\label{Prop:Laplacians-link} For every $C^\infty$ function $\eta:X\longrightarrow(0,\,+\infty)$ on a compact complex Hermitian manifold $(X,\,\omega)$ with $\mbox{dim}_\C X=n$ and for every $k\in\{0,\dots , 2n\}$, the operator $D_\eta^\star : C^\infty_k(X,\,\C)\longrightarrow C^\infty_{k-1}(X,\,\C)$ is given by the formula: \begin{eqnarray}\label{eqn:Laplacians-link}D_\eta^\star = \theta_\eta T^\star_\eta\theta_\eta^{-1},\end{eqnarray} where $T^\star_\eta : C^\infty_k(X,\,\C)\longrightarrow C^\infty_{k-1}(X,\,\C)$ is the linear operator \begin{eqnarray*}T^\star_\eta\bigg(\sum\limits_{p+q=k}u^{p,\,q}\bigg) = \sum\limits_{p+q=k} T^\star_{p,\,\eta}(u^{p,\,q}),  \hspace{5ex} \mbox{with}\hspace{1ex} u^{p,\,q}\in C^\infty_{p,\,q}(X,\,\C),\end{eqnarray*} and, for every bidegree $(p,\,q)$, $T^\star_{p,\,\eta}:C^\infty_{p,\,q}(X,\,\C)\longrightarrow C^\infty_{p+q-1}(X,\,\C)$ is the linear operator defined by \begin{eqnarray*}T^\star_{p,\,\eta} = \frac{1}{\eta^{2p}}\,\bigg(d^\star_{\omega_\eta} - 2(n-p)\,\bigg[\Lambda,\,i\bigg(\frac{1}{\eta}\,\partial\eta - \eta\,\bar\partial\eta\bigg)\wedge\cdot\bigg]\bigg)\circ\theta_\eta^2,\end{eqnarray*} where $d^\star_{\omega_\eta}$ is the formal adjoint of $d$ with respect to the $L^2$-inner product induced by the twisted metric $\omega_\eta$ as in (\ref{eqn:eta-L2-inner-product_def}).   

\end{Prop}

Before proving this result, we prove a few formulae that will be needed.

\begin{Lem}\label{Lem:del-star_both-metrics} Let $(X,\,\omega)$ be a complex Hermitian manifold with $\mbox{dim}_\C X=n$. Then:

  \vspace{1ex}

  (i)\, for any $C^1$ function $\rho:X\longrightarrow\R$, any bidegree $(p,\,q)$ and any smooth $(p,\,q)$-form $\gamma$ on $X$, we have: \begin{eqnarray}\label{eqn:del-star_del-bar-star_function}\partial^\star(\rho\gamma) = \rho\,\partial^\star\gamma + [\Lambda,\,i\bar\partial\rho\wedge\cdot](\gamma) \hspace{3ex}\mbox{and}\hspace{3ex} \bar\partial^\star(\rho\gamma) = \rho\,\bar\partial^\star\gamma - [\Lambda,\,i\partial\rho\wedge\cdot](\gamma).\end{eqnarray}

  \vspace{1ex}

  (ii)\, for any $C^\infty$ function $\eta:X\longrightarrow(0,\,+\infty)$, any bidegree $(p,\,q)$ and any smooth $(p,\,q)$-form $u$ on $X$, we have: \begin{eqnarray}\label{eqn:del-star_both-metrics}\partial^\star_{\omega_\eta} u = \eta^2\partial^\star_\omega u - 2(n-p)\eta\,[\Lambda,\,i\bar\partial\eta\wedge\cdot]\,u \hspace{3ex}\mbox{and}\hspace{3ex} \bar\partial^\star_{\omega_\eta} u = \bar\partial^\star_\omega u + \frac{2(n-p)}{\eta}\,[\Lambda,\,i\partial\eta\wedge\cdot]\,u.\end{eqnarray}

\end{Lem}  

\noindent{\it Proof.} (i)\, It suffices to prove the latter equality in (\ref{eqn:del-star_del-bar-star_function}) since the former follows from it by conjugation. Using the Hermitian commutation identity (ii) of (\ref{eqn:standard-comm-rel}), we get the first equality below: \begin{eqnarray*}\bar\partial^\star(\rho\gamma) & = & -i\,[\Lambda,\,\partial](\rho\gamma) - \bar\tau^{\star}(\rho\gamma) = -i\,\Lambda(\rho\,\partial\gamma + \partial\rho\wedge\gamma) + i\partial(\rho\,\Lambda\gamma) - \rho\,\bar\tau^{\star}(\gamma)\\
  & = & -i\rho\,\bigg(\Lambda\partial\gamma - \partial\Lambda\gamma\bigg) -i\,\bigg(\Lambda(\partial\rho\wedge\gamma) - \partial\rho\wedge\Lambda\gamma\bigg) - \rho\,\bar\tau^{\star}(\gamma) \\
  & = & \rho\,\bigg(-i\,[\Lambda,\,\partial] - \bar\tau^{\star}\bigg)(\gamma) - [\Lambda,\,i\partial\rho\wedge\cdot](\gamma) = \rho\,\bar\partial^\star\gamma - [\Lambda,\,i\partial\rho\wedge\cdot](\gamma),\end{eqnarray*} where the last equality follows again from the Hermitian commutation identity (ii) of (\ref{eqn:standard-comm-rel}).

\vspace{1ex}

(ii)\, $\bullet$ To prove the former equality in (\ref{eqn:del-star_both-metrics}), we will first prove the formula \begin{eqnarray}\label{eqn:del-star_omega_twisted-omega}\partial^\star_\omega = \frac{1}{\eta^{2(n-p+1)}}\, \partial^\star_{\omega_\eta}\bigg(\eta^{2(n-p)}\cdot\bigg)   \hspace{5ex} \mbox{on} \hspace{1ex} (p,\,q)\mbox{-forms}.\end{eqnarray}

To this end, let $\alpha\in C^\infty_{p-1,\,q}(X,\,\C)$ and $\beta\in C^\infty_{p,\,q}(X,\,\C)$. We have: \begin{eqnarray*}\langle\langle\alpha,\,\partial_\omega^\star\beta\rangle\rangle_\omega & = & \langle\langle\partial\alpha,\,\beta\rangle\rangle_\omega  = \int\limits_X\langle\partial\alpha,\,\beta\rangle_\omega\,dV_\omega = \int\limits_X\frac{1}{\eta^{2p}}\,\langle\partial\alpha,\,\beta\rangle_{\omega_\eta}\,\eta^{2n}\,dV_{\omega_\eta} = \langle\langle\partial\alpha,\,\eta^{2(n-p)}\beta\rangle\rangle_{\omega_\eta} \\
  & = & \langle\langle\alpha,\,\partial^\star_{\omega_\eta}(\eta^{2(n-p)}\beta)\rangle\rangle_{\omega_\eta} = \int\limits_X \langle\alpha,\,\partial^\star_{\omega_\eta}(\eta^{2(n-p)}\beta)\rangle_{\omega_\eta}\,dV_{\omega_\eta} \\
  & = & \int\limits_X \eta^{2(p-1)}\,\langle\alpha,\,\partial^\star_{\omega_\eta}(\eta^{2(n-p)}\beta)\rangle_\omega\,\frac{1}{\eta^{2n}}\,dV_\omega = \langle\langle\alpha,\,\frac{1}{\eta^{2(n-p+1)}}\,\partial^\star_{\omega_\eta}(\eta^{2(n-p)}\beta)\rangle\rangle_\omega.\end{eqnarray*} This proves (\ref{eqn:del-star_omega_twisted-omega}).

Now, let $u\in C^\infty_{p,\,q}(X,\,\C)$ be arbitrary. Identity (\ref{eqn:del-star_omega_twisted-omega}) gives the first equality below: \begin{eqnarray*}\partial^\star_{\omega_\eta}u = \eta^{2(n-p+1)}\,\partial_\omega^\star\bigg(\frac{1}{\eta^{2(n-p)}}\,u\bigg) = \eta^{2(n-p+1)}\,\bigg(\frac{1}{\eta^{2(n-p)}}\,\partial_\omega^\star u + \bigg[\Lambda,\,i\bar\partial\bigg(\frac{1}{\eta^{2(n-p)}}\bigg)\wedge\cdot\bigg]\,u\bigg),\end{eqnarray*} while the second equality follows from the former formula in (\ref{eqn:del-star_del-bar-star_function}). Since \begin{eqnarray*}\bar\partial\bigg(\frac{1}{\eta^{2(n-p)}}\bigg) = -\frac{2(n-p)}{\eta^{2(n-p)+1} }\,\bar\partial\eta,\end{eqnarray*} we get the former formula in (\ref{eqn:del-star_both-metrics}).

\vspace{1ex}

$\bullet$ To prove the latter equality in (\ref{eqn:del-star_both-metrics}), we will first prove the formula \begin{eqnarray}\label{eqn:del-bar-star_omega_twisted-omega}\bar\partial^\star_\omega = \frac{1}{\eta^{2(n-p)}}\, \bar\partial^\star_{\omega_\eta}\bigg(\eta^{2(n-p)}\cdot\bigg)   \hspace{5ex} \mbox{on} \hspace{1ex} (p,\,q)\mbox{-forms}.\end{eqnarray}

To this end, let $\alpha\in C^\infty_{p,\,q-1}(X,\,\C)$ and $\beta\in C^\infty_{p,\,q}(X,\,\C)$. We have: \begin{eqnarray*}\langle\langle\alpha,\,\bar\partial_\omega^\star\beta\rangle\rangle_\omega & = & \langle\langle\bar\partial\alpha,\,\beta\rangle\rangle_\omega  = \int\limits_X\langle\bar\partial\alpha,\,\beta\rangle_\omega\,dV_\omega = \int\limits_X\frac{1}{\eta^{2p}}\,\langle\bar\partial\alpha,\,\beta\rangle_{\omega_\eta}\,\eta^{2n}\,dV_{\omega_\eta} = \langle\langle\bar\partial\alpha,\,\eta^{2(n-p)}\beta\rangle\rangle_{\omega_\eta} \\
  & = & \langle\langle\alpha,\,\bar\partial^\star_{\omega_\eta}(\eta^{2(n-p)}\beta)\rangle\rangle_{\omega_\eta} = \int\limits_X \langle\alpha,\,\bar\partial^\star_{\omega_\eta}(\eta^{2(n-p)}\beta)\rangle_{\omega_\eta}\,dV_{\omega_\eta} \\
  & = & \int\limits_X \eta^{2p}\,\langle\alpha,\,\bar\partial^\star_{\omega_\eta}(\eta^{2(n-p)}\beta)\rangle_\omega\,\frac{1}{\eta^{2n}}\,dV_\omega = \langle\langle\alpha,\,\frac{1}{\eta^{2(n-p)}}\,\bar\partial^\star_{\omega_\eta}(\eta^{2(n-p)}\beta)\rangle\rangle_\omega.\end{eqnarray*} This proves (\ref{eqn:del-bar-star_omega_twisted-omega}).

Now, let $u\in C^\infty_{p,\,q}(X,\,\C)$ be arbitrary. Identity (\ref{eqn:del-bar-star_omega_twisted-omega}) gives the first equality below: \begin{eqnarray*}\bar\partial^\star_{\omega_\eta}u = \eta^{2(n-p)}\,\bar\partial_\omega^\star\bigg(\frac{1}{\eta^{2(n-p)}}\,u\bigg) = \eta^{2(n-p)}\,\bigg(\frac{1}{\eta^{2(n-p)}}\,\bar\partial_\omega^\star u - \bigg[\Lambda,\,i\partial\bigg(\frac{1}{\eta^{2(n-p)}}\bigg)\wedge\cdot\bigg]\,u\bigg),\end{eqnarray*} while the second equality follows from the latter formula in (\ref{eqn:del-star_del-bar-star_function}). Since \begin{eqnarray*}\partial\bigg(\frac{1}{\eta^{2(n-p)}}\bigg) = -\frac{2(n-p)}{\eta^{2(n-p)+1} }\,\partial\eta,\end{eqnarray*} we get the latter formula in (\ref{eqn:del-star_both-metrics}).  \hfill $\Box$

\vspace{3ex}

\noindent{\it Proof of Proposition \ref{Prop:Laplacians-link}.} On the one hand, taking conjugates in the definition (\ref{eqn:D_eta}) of $D_\eta$, we get: \begin{eqnarray}\label{eqn:Laplacians-link_proof_1}\theta_\eta^{-1}D_\eta^\star\theta_\eta = \theta_\eta^{-2}d_\omega^\star\theta_\eta^2\end{eqnarray} since the formal adjoint of $\theta_\eta$ with respect to any metric (in particular, in this case, with respect to $\omega$) is $\theta_\eta$ itself.

On the other hand, we will express the right-hand side of (\ref{eqn:Laplacians-link_proof_1}) in a different way. For any bidegree $(p,\,q)$ and any form $u\in C^\infty_{p,\,q}(X,\,\C)$, formulae (\ref{eqn:del-star_both-metrics}) give the latter equality below: \begin{eqnarray*}d_\omega^\star u = \partial_\omega^\star u + \bar\partial_\omega^\star u = \frac{1}{\eta^2}\,\partial_{\omega_\eta}^\star u + \bar\partial_{\omega_\eta}^\star u - 2\,\frac{n-p}{\eta}\,\bigg[\Lambda,\,i(\partial\eta - \bar\partial\eta)\wedge\cdot\bigg]\,u.\end{eqnarray*}

Replacing $u$ with $\theta_\eta^2 u = \eta^{2p}u$ and taking $\theta_\eta^{-2}$ on both sides, the last equality transforms to: \begin{eqnarray}\label{eqn:Laplacians-link_proof_2}(\theta_\eta^{-2}d_\omega^\star\theta_\eta^2)(u) = \frac{1}{\eta^{2p}}\,\partial_{\omega_\eta}^\star(\theta_\eta^2 u) + \frac{1}{\eta^{2p}}\,\bar\partial_{\omega_\eta}^\star(\theta_\eta^2 u) - 2\,\frac{n-p}{\eta^{2p}}\,\bigg[\Lambda,\,i\bigg(\frac{1}{\eta}\,\partial\eta - \eta\,\bar\partial\eta\bigg)\wedge\cdot\bigg]\,(\theta_\eta^2 u).\end{eqnarray}

Putting (\ref{eqn:Laplacians-link_proof_1}) and (\ref{eqn:Laplacians-link_proof_2}) together, we infer that, for every $u\in C^\infty_{p,\,q}(X,\,\C)$, we have: \begin{eqnarray*}(\theta_\eta^{-1}D_\eta^\star\theta_\eta)u = \frac{1}{\eta^{2p}}\,\bigg(d_{\omega_\eta}^\star - 2(n-p)\,\bigg[\Lambda,\,i\bigg(\frac{1}{\eta}\,\partial\eta - \eta\,\bar\partial\eta\bigg)\wedge\cdot\bigg]\bigg)\,(\theta_\eta^2 u) = T^\star_{p,\,\eta}u.\end{eqnarray*} 

This proves formula (\ref{eqn:Laplacians-link}) in bidegree $(p,\,q)$. By linearity, we get (\ref{eqn:Laplacians-link}) in any degree $k$.  \hfill $\Box$

\vspace{3ex}

Prompted by Proposition \ref{Prop:Laplacians-link}, we introduce the following

\begin{Def}\label{Def:Delta_omega-eta} Let $(X,\,\omega)$ be a {\bf compact} complex Hermitian manifold with $\mbox{dim}_\C X=n$.  

For any $C^\infty$ function $\eta:X\longrightarrow(0,\,\infty)$ and any $k\in\{0,\dots , 2n\}$, the twisted $d$-Laplacian $\Delta_{\omega_\eta}:C^\infty_k(X,\,\C)\longrightarrow C^\infty_k(X,\,\C)$ with respect to the twisted metric $\omega_\eta$ is the differential operator defined by \begin{eqnarray}\label{eqn:Delta_omega-eta}\Delta_{\omega_\eta} = d T_\eta^\star + T_\eta^\star d,\end{eqnarray} where $T_\eta^\star$ is the first-order differential operator introduced in Proposition \ref{Prop:Laplacians-link}.

\end{Def}  

In the special case where $\eta$ is {\it constant}, we have $T_\eta^\star = d^\star_{\omega_\eta}$, so $\Delta_{\omega_\eta}$ coincides with the usual $d$-Laplacian $dd^\star_{\omega_\eta} + d^\star_{\omega_\eta}d$ with respect to $\omega_\eta$. In general, for an arbitrary smooth function $\eta>0$, $d^\star_{\omega_\eta}$ is the principal part of $T_\eta^\star$, so the usual $d$-Laplacian with respect to $\omega_\eta$ is the principal part of $\Delta_{\omega_\eta}$.

\begin{Cor}\label{Cor:Laplacians-link} Let $(X,\,\omega)$ be a {\bf compact} complex Hermitian manifold with $\mbox{dim}_\C X=n$.

For any $C^\infty$ function $\eta:X\longrightarrow(0,\,\infty)$ and any $k\in\{0,\dots , 2n\}$, the twisted Laplacians $\Delta_\eta,\,\Delta_{\omega_\eta}:C^\infty_k(X,\,\C)\longrightarrow C^\infty_k(X,\,\C)$ are related by the formula: \begin{eqnarray}\label{eqn:Laplacians-link}\Delta_\eta = \theta_\eta\Delta_{\omega_\eta}\theta_\eta^{-1}.\end{eqnarray}

In particular, they have the same spectrum: \begin{eqnarray}\label{eqn:Laplacians_same-spectrum}\mbox{Spec}(\Delta_\eta) = \mbox{Spec}(\Delta_{\omega_\eta})\end{eqnarray} and for every eigenvalue $\lambda$, the linear map \begin{eqnarray}\label{eqn:Laplacians_eigenspaces_isom}E^k_{\Delta_\eta}(\lambda)\ni u\longmapsto\theta_\eta^{-1}u\in E^k_{\Delta_{\omega_\eta}}(\lambda)\end{eqnarray} is an isomorphism between the corresponding eigenspaces.

\end{Cor}

\noindent {\it Proof.} Formulae (\ref{eqn:D_eta}) and (\ref{eqn:Delta_omega-eta}) yield the second equality below: \begin{eqnarray*}\Delta_\eta = D_\eta D_\eta^\star + D_\eta^\star D_\eta = (\theta_\eta d\theta_\eta^{-1})(\theta_\eta T_\eta^\star\theta_\eta^{-1}) + (\theta_\eta T_\eta^\star\theta_\eta^{-1})(\theta_\eta d\theta_\eta^{-1}) = \theta_\eta(dT_\eta^\star + T_\eta^\star d)\theta_\eta^{-1} = \theta_\eta\Delta_{\omega_\eta}\theta_\eta^{-1}.\end{eqnarray*}

This proves (\ref{eqn:Laplacians-link}), from which the other two claims follow at once.  \hfill $\Box$

\section{Twisted commutation relations for $D_\eta$ and $\overline{D}_{-\eta}$}\label{section:twisted-commutation_D-eta} Let $X$ be a (possibly non-compact) complex manifold with $\mbox{dim}_\C X=n\geq 2$. We fix a Hermitian metric $\omega$ on $X$ and denote by $\langle\,\,,\,\,\rangle$, resp. $\langle\langle\,\,,\,\,\rangle\rangle$, the pointwise, resp. $L^2$, inner product induced by $\omega$ on $\C$-valued differential forms on $X$.

In this section, we give commutation relations for $D_\eta$, $\overline{D}_{-\eta}$ and their formal adjoints, as well as the identities of the Bochner-Kodaira-Nakano-type they induce. The identities we obtain are the analogues in our twisted and possibly non-K\"ahler context of the classical K\"ahler commutation relations that were subsequently given Hermitian versions in [Gri69], [Ohs82] and [Dem84]  for the standard operators $\partial$ and $\bar\partial$ and then for the operators $d_h$ and $\overline{d}_{-h}$ twisted by a constant $h\in\C$ in [BP18]. In our present case, the twisting is by a possibly non-constant function $\eta$. 

We start with a preliminary computation.

\begin{Lem}\label{Lem:D_eta_star_prelim} For any $C^\infty$ function $\eta$ on $X$ such that $\eta>0$ or $\eta<0$ and any $k\in\{0,\dots , 2n\}$, the formal adjoint $D_\eta^\star:C^\infty_{k+1}(X,\,\C)\longrightarrow C^\infty_k(X,\,\C)$ w.r.t. $\langle\langle\,\,,\,\,\rangle\rangle$ of the operator $D_\eta:C^\infty_k(X,\,\C)\longrightarrow C^\infty_{k+1}(X,\,\C)$ introduced in Definition \ref{Def:D_eta} is given by the formula: \begin{eqnarray}\label{eqn:D_eta_star_prelim}D_\eta^\star v = \partial^\star(\eta v) + \bar\partial^\star v - \sum\limits_{p+q=k}p\,(\partial\eta\wedge\cdot)^\star v^{p+1,\,q} - \sum\limits_{p+q=k}\frac{p}{\eta}\,(\bar\partial\eta\wedge\cdot)^\star v^{p,\,q+1}\end{eqnarray} for every form $v=\sum\limits_{r+s=k+1}v^{r,\,s}\in C^\infty_{k+1}(X,\,\C)$.

\end{Lem}  

\noindent {\it Proof.} For any forms $u = \sum\limits_{p+q=k}u^{p,\,q}\in C^\infty_k(X,\,\C)$ and $v=\sum\limits_{r+s=k+1}v^{r,\,s}\in C^\infty_{k+1}(X,\,\C)$, formula (\ref{eqn:D_eta_explicit}) for $D_\eta u$ and the fact that the inner product of any two pure-type forms of different types vanishes lead to the following equivalences: \begin{eqnarray*}\langle\langle D_\eta u,\,v\rangle\rangle = \langle\langle u,\,D_\eta^\star v\rangle\rangle\end{eqnarray*} \begin{eqnarray*} & \iff & \sum\limits_{p+q=k}\langle\langle\eta\partial u^{p,\,q} - p\,\partial\eta\wedge u^{p,\,q},\,v^{p+1,\,q}\rangle\rangle + \sum\limits_{p+q=k}\langle\langle\bar\partial u^{p,\,q} - \frac{p}{\eta}\,\bar\partial\eta\wedge u^{p,\,q},\,v^{p,\,q+1}\rangle\rangle = \langle\langle u,\,D_\eta^\star v\rangle\rangle \\
    & \iff & \sum\limits_{p+q=k}\langle\langle u^{p,\,q},\,\partial^\star(\eta\,v^{p+1,\,q}) - p\,(\partial\eta\wedge\cdot)^\star v^{p+1,\,q} + \bar\partial^\star v^{p,\,q+1} - \frac{p}{\eta}\,(\bar\partial\eta\wedge\cdot)^\star v^{p,\,q+1} \rangle\rangle   = \langle\langle u,\,D_\eta^\star v\rangle\rangle.\end{eqnarray*}

  Now, the last term translates to \begin{eqnarray*}\langle\langle u,\,D_\eta^\star v\rangle\rangle = \sum\limits_{p+q=k}\langle\langle u^{p,\,q},\,(D_\eta^\star v)^{p,\,q}\rangle\rangle,\end{eqnarray*} so we get: \begin{eqnarray*}(D_\eta^\star v)^{p,\,q} = \partial^\star(\eta\,v^{p+1,\,q}) + \bar\partial^\star v^{p,\,q+1} - p\,(\partial\eta\wedge\cdot)^\star v^{p+1,\,q} - \frac{p}{\eta}\,(\bar\partial\eta\wedge\cdot)^\star v^{p,\,q+1}\end{eqnarray*} for every bidegree $(p,\,q)$. Summing up over $p+q=k$, we get (\ref{eqn:D_eta_star_prelim}). \hfill $\Box$

\vspace{3ex}

Now, we fix a degree $k\in\{0,\dots , 2n\}$ and a form $v=\sum\limits_{r+s=k+1}v^{r,\,s}\in C^\infty_{k+1}(X,\,\C)$. Expressing $\partial^\star$, $\bar\partial^\star$, $(\partial\eta\wedge\cdot)^\star$ and $(\bar\partial\eta\wedge\cdot)^\star$ by means of the commutation relations of Lemmas \ref{Lem:com_1} and \ref{Lem:com_2}, formula (\ref{eqn:D_eta_star}) reads: \begin{eqnarray*}D_\eta^\star v & = & i[\Lambda,\,\bar\partial](\eta v) - \eta\,\tau^\star(v) - i[\Lambda,\,\partial] (v) - \bar\tau^\star(v) \\
  & + & i\,\sum\limits_{p+q=k}p\,[\Lambda,\,\bar\partial\eta\wedge\cdot](v^{p+1,\,q}) - i\,\sum\limits_{p+q=k}\frac{p}{\eta}\,[\Lambda,\,\partial\eta\wedge\cdot](v^{p,\,q+1}).\end{eqnarray*} Since \begin{eqnarray*}[\Lambda,\,\bar\partial](\eta v) = \Lambda(\eta\,\bar\partial v + \bar\partial\eta\wedge v) - \eta\,\bar\partial\Lambda v - \bar\partial\eta\wedge\Lambda v = \eta\,[\Lambda,\,\bar\partial](v)  + [\Lambda,\,\bar\partial\eta\wedge\cdot](v),\end{eqnarray*} the above equality translates to \begin{eqnarray}\label{eqn:D_eta_star_sequel_1}\nonumber D_\eta^\star v + (\eta\,\tau + \bar\tau)^\star(v) & = & -i\,[\Lambda,\,\partial-\eta\bar\partial](v) + [\Lambda,\,i\bar\partial\eta\wedge\cdot](v)\\
 & - & \sum\limits_{p+q=k}\frac{p}{\eta}\,[\Lambda,\,i\partial\eta\wedge\cdot](v^{p,\,q+1}) + \sum\limits_{p+q=k}p\,[\Lambda,\,i\bar\partial\eta\wedge\cdot](v^{p+1,\,q}).\end{eqnarray}
  
On the other hand, using formula (\ref{eqn:D_bar_minus-eta}) for $\overline{D}_{-\eta}$, we get:  \begin{eqnarray*}[\Lambda,\,\overline{D}_{-\eta}](v) & = & \sum\limits_{r+s=k+1} [\Lambda,\,\overline{D}_{-\eta}](v^{r,\,s}) \\
  & = & -\sum\limits_{r+s=k+1}\Lambda\bigg((\eta\bar\partial - \partial)v^{r,\,s} - s\,\bar\partial\eta\wedge v^{r,\,s} + \frac{s}{\eta}\,\partial\eta\wedge v^{r,\,s}\bigg) \\
  & + & \sum\limits_{r+s=k+1}(\eta\bar\partial - \partial)(\Lambda v^{r,\,s}) - \sum\limits_{r+s=k+1}(s-1)\,\bar\partial\eta\wedge\Lambda v^{r,\,s} + \sum\limits_{r+s=k+1}\frac{s-1}{\eta}\,\partial\eta\wedge\Lambda v^{r,\,s},\end{eqnarray*} which translates to \begin{eqnarray}\label{eqn:D_eta_star_sequel_2}\nonumber -i[\Lambda,\,\overline{D}_{-\eta}](v) & = & -i[\Lambda,\,\partial - \eta\bar\partial](v) - \sum\limits_{r+s=k+1}s\,[\Lambda,\,i\bar\partial\eta\wedge\cdot]\,v^{r,\,s} -i\bar\partial\eta\wedge\Lambda(v) \\
  & +  & \sum\limits_{r+s=k+1}\frac{s}{\eta}\,[\Lambda,\,i\partial\eta\wedge\cdot]\,v^{r,\,s} + \sum\limits_{r+s=k+1}\frac{i}{\eta}\,\partial\eta\wedge\Lambda(v^{r,\,s}).\end{eqnarray}

The conclusion of this computation is a preliminary $\eta$-twisted commutation relation.

\begin{Lem}\label{Lem:D_eta_star} For any $C^\infty$ function $\eta$ on $X$ such that $\eta>0$ or $\eta<0 $ and any $k\in\{0,\dots , 2n\}$, the formal adjoint $D_\eta^\star:C^\infty_{k+1}(X,\,\C)\longrightarrow C^\infty_k(X,\,\C)$ of $D_\eta:C^\infty_k(X,\,\C)\longrightarrow C^\infty_{k+1}(X,\,\C)$ is given by \begin{eqnarray}\label{eqn:D_eta_star}D_\eta^\star + (\eta\,\tau + \bar\tau)^\star = -i\,[\Lambda,\,\overline{D}_{-\eta}] - \frac{k+1}{\eta}[\Lambda,\,i\,\bar{d}_{-\eta}\eta\wedge\cdot] -\frac{i}{\eta}\,\bar{d}_{-\eta}\eta\wedge\Lambda.\end{eqnarray}

\end{Lem}  

\noindent {\it Proof.} Let $v=\sum\limits_{r+s=k+1}v^{r,\,s}\in C^\infty_{k+1}(X,\,\C)$.

Plugging into (\ref{eqn:D_eta_star_sequel_1}) the expression we obtain for $-i[\Lambda,\,\partial - \eta\bar\partial](v)$ from (\ref{eqn:D_eta_star_sequel_2}), we get:  \begin{eqnarray}\label{eqn:D_eta_star_sequel_3}\nonumber D_\eta^\star v + (\eta\,\tau + \bar\tau)^\star(v) & = & -i[\Lambda,\,\overline{D}_{-\eta}](v) + \sum\limits_{r+s=k+1}s\,[\Lambda,\,i\bar\partial\eta\wedge\cdot\,]\,(v^{r,\,s}) + i\bar\partial\eta\wedge\Lambda(v) \\
  & - & \sum\limits_{r+s=k+1}\frac{s}{\eta}\,[\Lambda,\,i\partial\eta\wedge\cdot\,]\,(v^{r,\,s}) - \sum\limits_{r+s=k+1}\frac{i}{\eta}\,\partial\eta\wedge\Lambda(v^{r,\,s}) \\
\nonumber & +  & [\Lambda,\,i\bar\partial\eta\wedge\cdot\,]\,(v)  -  \sum\limits_{p+q=k}\frac{p}{\eta}\,[\Lambda,\,i\partial\eta\wedge\cdot\,](v^{p,\,q+1}) + \sum\limits_{p+q=k}p\,[\Lambda,\,i\bar\partial\eta\wedge\cdot\,](v^{p+1,\,q}).\end{eqnarray}

We now form pairs using certain sums featuring on the r.h.s. of (\ref{eqn:D_eta_star_sequel_3}). We get: \begin{eqnarray*}\sum\limits_{r+s=k+1}s\,[\Lambda,\,i\bar\partial\eta\wedge\cdot\,]\,(v^{r,\,s}) + \sum\limits_{p+q=k}p\,[\Lambda,\,i\bar\partial\eta\wedge\cdot\,]\,(v^{p+1,\,q}) = k\,[\Lambda,\,i\bar\partial\eta\wedge\cdot\,]\,(v)\end{eqnarray*} after renaming, in the latter sum, $p$ as $r-1$ and $q$ as $s$.

Similarly, we get: \begin{eqnarray*} - \sum\limits_{r+s=k+1}\frac{s}{\eta}\,[\Lambda,\,i\partial\eta\wedge\cdot\,]\,(v^{r,\,s}) - \sum\limits_{p+q=k}\frac{p}{\eta}\,[\Lambda,\,i\partial\eta\wedge\cdot\,](v^{p,\,q+1}) = - \frac{k+1}{\eta}\,[\Lambda,\,i\partial\eta\wedge\cdot]\,(v)\end{eqnarray*} after renaming, in the latter sum, $p$ as $r$ and $q$ as $s-1$.

Finally, after using formula (\ref{eqn:d_bar_minus-eta}) giving $\overline{d}_{-\eta}$ for the latter equality below, we get: \begin{eqnarray*}i\bar\partial\eta\wedge\Lambda(v) - \sum\limits_{r+s=k+1}\frac{i}{\eta}\,\partial\eta\wedge\Lambda(v^{r,\,s}) = i\bar\partial\eta\wedge\Lambda(v) - \frac{i}{\eta}\,\partial\eta\wedge\Lambda(v) = - \frac{i}{\eta}\,\bar{d}_{-\eta}\eta\wedge\Lambda(v).\end{eqnarray*}

Putting all these expressions together and using (\ref{eqn:d_bar_minus-eta}), we see that (\ref{eqn:D_eta_star_sequel_3}) becomes the desired (\ref{eqn:D_eta_star}). \hfill $\Box$

\vspace{3ex}

Note that the second and third terms on the r.h.s. of (\ref{eqn:D_eta_star}) are of order zero. Thus, the principal part of $D_\eta^\star$ is contained in $-i\,[\Lambda,\,\overline{D}_{-\eta}]$. 

\vspace{2ex}

We now introduce the following

\begin{Def}\label{Def:eta-torsion}


  For every degree $k\in\{0,\dots , 2n\}$, the {\bf $\eta$-torsion operator} $\tau_\eta:\Lambda^kT^\star X\longrightarrow\Lambda^{k+1}T^\star X$ is defined pointwise on the $\C$-valued $k$-forms on $X$ by \begin{eqnarray}\label{eqn:eta-torsion}\tau_\eta = [\Lambda,\,D_\eta\omega\wedge\cdot\,].\end{eqnarray}

\end{Def}

A straightforward computation yields the following

\begin{Lem}\label{Lem:eta-torsion_computed} For every $k\in\{0,\dots , 2n\}$, the $\eta$-torsion operator $\tau_\eta:\Lambda^kT^\star X\longrightarrow\Lambda^{k+1}T^\star X$ is explicitly given by \begin{eqnarray}\label{eqn:eta-torsion_computed}\tau_\eta = \eta\tau + \bar\tau + i\,\omega\wedge(\bar\partial\eta\wedge\cdot\,)^\star - \frac{i}{\eta}\,(\partial\eta\wedge\cdot\,)^\star\,(\omega\wedge\cdot\,) - (n-k-1)\,\partial\eta\wedge\cdot - (n-k)\,\frac{1}{\eta}\,\bar\partial\eta\wedge\cdot.\end{eqnarray}

\end{Lem}

\noindent {\it Proof.} From (\ref{eqn:D_eta_explicit}), we get: \begin{eqnarray*}D_\eta\omega = \eta\partial\omega + \bar\partial\omega - \partial\eta\wedge\omega - \frac{1}{\eta}\,\bar\partial\eta\wedge\omega.\end{eqnarray*} Hence, for every $k$-form $u$, we get: \begin{eqnarray*}\tau_\eta u & = & \Lambda\bigg(\eta\,\partial\omega\wedge u +  \bar\partial\omega\wedge u - \partial\eta\wedge\omega\wedge u - \frac{1}{\eta}\,\bar\partial\eta\wedge\omega\wedge u\bigg)  \\
  & - & \bigg(\eta\,\partial\omega +  \bar\partial\omega -\partial\eta\wedge\omega - \frac{1}{\eta}\,\bar\partial\eta\wedge\omega\bigg)\wedge\Lambda u \\
  & = & \eta\,[\Lambda,\,\partial\omega\wedge\cdot\,]\,u + [\Lambda,\,\bar\partial\omega\wedge\cdot\,]\,u - [\Lambda,\,\omega\wedge\cdot\,]\,(\partial\eta\wedge u) - \omega\wedge[\Lambda,\,\partial\eta\wedge\cdot\,]\,u \\
  & - & \frac{1}{\eta}\,[\Lambda,\,\bar\partial\eta\wedge\cdot\,]\,(\omega\wedge u) - \frac{1}{\eta}\,\bar\partial\eta\wedge[\Lambda,\,\omega\wedge\cdot\,]\,u \\
  & = & \eta\,\tau(u) + \bar\tau(u) - (n-k-1)\,\partial\eta\wedge u + i\,\omega\wedge(\bar\partial\eta\wedge\cdot\,)^\star (u)\\
  & - &  \frac{1}{\eta}\,i (\partial\eta\wedge\cdot\,)^\star\,(\omega\wedge u) - \frac{1}{\eta}\,\bar\partial\eta\wedge (n-k)\,u,\end{eqnarray*} which is (\ref{eqn:eta-torsion_computed}).

We have used the definition of $\tau$, formulae (a) of Lemma \ref{Lem:com_2} and the standard identity $[\Lambda,\,\omega\wedge\cdot\,] = (n-k)\,\mbox{Id}$ on $k$-forms.  \hfill $\Box$  

\vspace{2ex}

As a consequence of these computations, we get the following analogue of Lemma 4.1. of [BP18] (which dealt with the case where the function $\eta$ was a constant $h\in\R\setminus\{0\}$). 

\begin{Prop}\label{Prop:eta-twisted_commutation-relations} Let $(X,\,\omega)$ be a complex Hermitian manifold with $\mbox{dim}_\C X=n\geq 2$. For any $C^\infty$ function $\eta$ on $X$ such that $\eta>0$ or $\eta<0$, the following {\bf $\eta$-twisted commutation relations} hold on differential forms of any degree on $X$: \begin{eqnarray}\label{eqn:eta-twisted_commutation-relations}\nonumber & (a)& D_\eta^\star +\tau_\eta^\star = -i\,[\Lambda,\, \overline{D}_{-\eta}] + \frac{n}{\eta}\,[i\,\overline{d}_{-\eta}\eta\wedge\cdot\,,\Lambda];  \hspace{2ex}   (b)\hspace{2ex}\overline{D}_{-\eta}^\star +\overline{\tau}_{-\eta}^\star = i\,[\Lambda,\, D_\eta] - \frac{n}{\eta}\,[i\,d_\eta\eta\wedge\cdot\,,\Lambda];\\
    \nonumber & (c) & D_\eta + \tau_\eta = i\,[\overline{D}_{-\eta}^\star,\,\omega\wedge\cdot\,] - \frac{n}{\eta}\,[\omega\wedge\cdot\,,i\,(\overline{d}_{-\eta}\eta\wedge\cdot\,)^\star];  \\
    & (d) & \overline{D}_{-\eta} + \overline{\tau}_{-\eta} = -i\,[D_\eta^\star,\,\omega\wedge\cdot\,] + \frac{n}{\eta}\,[\omega\wedge\cdot\,,i\,(d_\eta\eta\wedge\cdot\,)^\star].\end{eqnarray}

\end{Prop}

\noindent {\it Proof.} It suffices to prove (a), since (b) will then follow by taking conjugates and replacing $\eta$ with $-\eta$, (c) will follow by taking adjoints in (a), while (d) will follow by taking adjoints in (b).

Taking adjoints in (\ref{eqn:eta-torsion_computed}), we get: \begin{eqnarray*}\tau_\eta^\star = \eta\tau^\star + \bar\tau^\star - i\,\bar\partial\eta\wedge\Lambda + \frac{i}{\eta}\,\Lambda(\partial\eta\wedge\cdot\,) - (n-k-1)\,(\partial\eta\wedge\cdot\,)^\star - (n-k)\,\frac{1}{\eta}\,(\bar\partial\eta\wedge\cdot\,)^\star\end{eqnarray*} in degree $k+1$. (Note that the switch from the degree $k$ of (\ref{eqn:eta-torsion_computed}) to the current degree $k+1$ is due to $\tau_\eta u$ being a $(k+1)$-form whenever of $u$ is a $k$-form.)

  Plugging the value for $\eta\tau^\star + \bar\tau^\star$ given by the above equality into (\ref{eqn:D_eta_star}), we get:  \begin{eqnarray}\label{eqn:eta-twisted_commutation-relations_proof_1}\nonumber D_\eta^\star + \tau_\eta^\star = -i\,[\Lambda,\,\overline{D}_{-\eta}] & - & \frac{k+1}{\eta}\,\Lambda(i\,\partial\eta\wedge\cdot\,) + \frac{k+1}{\eta}\,\Lambda(i\eta\,\bar\partial\eta\wedge\cdot\,) + \frac{k+1}{\eta}\,i\partial\eta\wedge\Lambda - \frac{k+1}{\eta}\,\eta\,i\bar\partial\eta\wedge\Lambda \\
    \nonumber   & - & \frac{i}{\eta}\,\partial\eta\wedge\Lambda + \frac{i}{\eta}\,\eta\,\bar\partial\eta\wedge\Lambda - i\bar\partial\eta\wedge\Lambda + \frac{i}{\eta}\,\Lambda(\partial\eta\wedge\cdot\,) \\
         & - & (n-k-1)\,(\partial\eta\wedge\cdot\,)^\star - (n-k)\,\frac{1}{\eta}\,(\bar\partial\eta\wedge\cdot\,)^\star\end{eqnarray} in degree $k+1$.

  Now, using the commutation relations (a) of Lemma \ref{Lem:com_2}, we see that the last line translates to \begin{eqnarray*}-(n-k-1)\,i\bar\partial\eta\wedge\Lambda + (n-k-1)\,\Lambda(i\bar\partial\eta\wedge\cdot\,) + \frac{n-k}{\eta}\,i\partial\eta\wedge\Lambda - \frac{n-k}{\eta}\,\Lambda(i\partial\eta\wedge\cdot\,).\end{eqnarray*}

  Therefore, (\ref{eqn:eta-twisted_commutation-relations_proof_1}) becomes: \begin{eqnarray*}D_\eta^\star + \tau_\eta^\star = -i\,[\Lambda,\,\overline{D}_{-\eta}] + \frac{n}{\eta}\,[i\partial\eta\wedge\cdot\,,\Lambda] - n\,[i\bar\partial\eta\wedge\cdot\,,\Lambda],\end{eqnarray*} which is nothing but (a). \hfill $\Box$

\vspace{2ex}

We are now in a position to give the first main result of this section. It generalises Corollary 4.2. of [BP18] (which dealt with the case where the function $\eta$ was a constant $h\in\R\setminus\{0\}$).

\begin{Prop}\label{Prop:eta-BKN_rough} Let $(X,\,\omega)$ be a complex Hermitian manifold with $\mbox{dim}_\C X=n\geq 2$.  For any $C^\infty$ function $\eta$ on $X$ such that $\eta>0$ or $\eta<0$ and any $k\in\{0,\dots , 2n\}$, the following {\bf rough $\eta$-Bochner-Kodaira-Nakano ($\eta$-BKN) identity} holds on the $C^\infty$ forms on $X$: \begin{eqnarray}\label{eqn:eta-BKN_rough}\nonumber\Delta_\eta = \overline{\Delta}_{-\eta} + i\,[[D_\eta,\,\overline{D}_{-\eta}],\,\Lambda] & + & [\overline{D}_{-\eta},\,\overline\tau_{-\eta}^\star] - [D_\eta,\,\tau_\eta^\star] \\
    & + & n\,\bigg[D_\eta,\,\frac{1}{\eta}\,[i\overline{d}_{-\eta}\eta\wedge\cdot\,,\,\Lambda]\bigg] + n\,\bigg[\overline{D}_{-\eta},\,\frac{1}{\eta}\,[id_\eta\eta\wedge\cdot\,,\,\Lambda]\bigg].\end{eqnarray}

\end{Prop}  

\noindent {\it Proof.} Using the expression for $D_\eta^\star$ given in (a) of (\ref{eqn:eta-twisted_commutation-relations}), we get the second equality below: \begin{eqnarray}\label{eqn:eta-BKN_proof_1}\Delta_\eta & = & [D_\eta,\,D_\eta^\star] =  -i\,[D_\eta,\,[\Lambda,\,\overline{D}_{-\eta}]] - [D_\eta,\,\tau_\eta^\star] + n\,\bigg[D_\eta,\,\frac{1}{\eta}\,[i\overline{d}_{-\eta}\eta\wedge\cdot\,,\,\Lambda]\bigg].\end{eqnarray}

Now, the Jacobi identity yields the former equality below (expressing the first term on the r.h.s. of (\ref{eqn:eta-BKN_proof_1})), while the $\eta$-twisted commutation relation (b) of (\ref{eqn:eta-twisted_commutation-relations}) yields the latter equality: \begin{eqnarray}\label{eqn:eta-BKN_proof_2}\nonumber -i\,[D_\eta,\,[\Lambda,\,\overline{D}_{-\eta}]] & = & [\overline{D}_{-\eta},\,i\,[\Lambda,\,D_\eta]] + i\,[[\overline{D}_{-\eta},\,D_\eta],\,\Lambda] \\
  & = & i\,[[D_\eta,\,\overline{D}_{-\eta}],\,\Lambda] + \bigg[\overline{D}_{-\eta},\,\overline{D}_{-\eta}^\star + \overline\tau_{-\eta}^\star + \frac{n}{\eta}\,[id_\eta\eta\wedge\cdot,\,\Lambda]\bigg].\end{eqnarray}

Plugging into (\ref{eqn:eta-BKN_proof_1}) the expression given for $-i\,[D_\eta,\,[\Lambda,\,\overline{D}_{-\eta}]]$ in (\ref{eqn:eta-BKN_proof_2}) and using the equality $[\overline{D}_{-\eta},\,\overline{D}_{-\eta}^\star] = \overline{\Delta}_{-\eta}$, we get (\ref{eqn:eta-BKN_rough}). \hfill $\Box$



\vspace{2ex}

We shall now refine the above $\eta$-BKN identity by incorporating some of the first-order terms on the right into a twisted Laplacian. The starting point is the following analogue for a possibly non-constant function $\eta$ of Lemma 4.4. in [BP18].

\begin{Lem}\label{Lem:preliminary_BKN-refined} Let $(X,\,\omega)$ be a complex Hermitian manifold with $\mbox{dim}_\C X=n\geq 2$.  For any $C^\infty$ function $\eta$ on $X$ such that $\eta>0$ or $\eta<0$ and any bidegree $(p,\,q)$, the following identities hold: \\

 \hspace{1ex} (i)\hspace{2ex} $[L,\,\tau_\eta] = 3\,D_\eta\omega\wedge\cdot\,$,  \hspace{6ex}   (ii)\, $[\Lambda,\,\tau_\eta] = 2i\,\overline{\tau}_{-\eta}^\star$,

  \vspace{2ex}
  
\hspace{1ex}  (iii)\hspace{2ex}$[D_\eta,\,\overline{D}_{-\eta}^\star] = -[D_\eta,\,\overline{\tau}_{-\eta}^\star] - \bigg[D_\eta,\,\frac{n}{\eta}\,(\overline{d}_{-\eta}\eta\wedge\cdot)^\star\bigg],$ 

 \begin{eqnarray*}(iv)\hspace{2ex} [D_\eta,\,D_\eta^\star] + [D_\eta,\,\tau_\eta^\star] -[\overline{D}_{-\eta},\,\overline{\tau}_{-\eta}^\star]  = [D_\eta + \tau_\eta,\,D_\eta^\star + \tau_\eta^\star] + S_\omega^{(\eta)} + \frac{n}{\eta}\,[\tau_\eta,\,(d_\eta\eta\wedge\cdot)^\star],\end{eqnarray*} where

\begin{eqnarray*}S_\omega^{(\eta)} : = \frac{i}{2}\,[\Lambda,\,[\Lambda,\,\overline{D}_{-\eta}D_\eta\omega\wedge\cdot\,]] - [D_\eta\omega\wedge\cdot\,,\, (D_\eta\omega\wedge\cdot\,)^\star].\end{eqnarray*}

\end{Lem}

\noindent {\it Proof.} (i)\, The definition of $\tau_\eta$ and the Jacobi identity yield the first and respectively the second identities below: $$[L,\,\tau_\eta] = [L,\,[\Lambda,\,D_\eta\omega\wedge\cdot\,]] = - [\Lambda,\,[D_\eta\omega\wedge\cdot\,,\,L]] - [D_\eta\omega\wedge\cdot\,,\,[L,\,\Lambda]].$$ 

 Now, $[D_\eta\omega\wedge\cdot\,,\,L] = D_\eta\omega\wedge(\omega\wedge\cdot\,) - \omega\wedge(D_\eta\omega\wedge\cdot\,) =0$, so the first term on the r.h.s. above vanishes. Meanwhile, it is standard that $[L,\,\Lambda]=(k-n)\,\mbox{Id}$ on $k$-forms. So for any $k$-form $u$, we get $$[D_\eta\omega\wedge\cdot\,,\,[L,\,\Lambda]]\,u = D_\eta\omega\wedge([L,\,\Lambda]\,u) - [L,\,\Lambda]\,(D_\eta\omega\wedge u) = (k-n)\,D_\eta\omega\wedge u - (k+3-n)\,D_\eta\omega\wedge u = -3\,D_\eta\omega\wedge u.$$ 

 Thus, $[D_\eta\omega\wedge\cdot\,,\,[L,\,\Lambda]] = -3\,D_\eta\omega\wedge\cdot$ and (i) follows.

\vspace{2ex}

(ii)\, We know from the $\eta$-twisted commutation relation (c) of Proposition \ref{Prop:eta-twisted_commutation-relations} that $\tau_\eta = i\,[\overline{D}_{-\eta}^\star,\,L] - D_\eta - \frac{n}{\eta}\,[L,i\,(\overline{d}_{-\eta}\eta\wedge\cdot\,)^\star]$. This implies the former identity below, while (b) of Proposition \ref{Prop:eta-twisted_commutation-relations} yields the latter: \begin{eqnarray}\label{Lem:preliminary_BKN-refined_proof_1}\nonumber [\Lambda,\,\tau_\eta] & = & i\,[\Lambda,\,[\overline{D}_{-\eta}^\star,\,L]] - [\Lambda,\,D_\eta] - \frac{n}{\eta}\,[\Lambda,\,[L,i\,(\overline{d}_{-\eta}\eta\wedge\cdot\,)^\star]] \\
  & = & i\,[\Lambda,\,[\overline{D}_{-\eta}^\star,\,L]] + i\,(\overline{D}_{-\eta}^\star + \overline\tau_{-\eta}^\star) - \frac{n}{\eta}\,[d_\eta\eta\wedge\cdot\,,\,\Lambda] - \frac{n}{\eta}\,[\Lambda,\,[L,i\,(\overline{d}_{-\eta}\eta\wedge\cdot\,)^\star]].\end{eqnarray}

We transform the first term on the right of (\ref{Lem:preliminary_BKN-refined_proof_1}) using the Jacobi identity \begin{eqnarray}\label{Lem:preliminary_BKN-refined_proof_1_bis}[\Lambda,\,[\overline{D}_{-\eta}^\star,\,L]] + [\overline{D}_{-\eta}^\star,\,[L,\,\Lambda]] + [L,\,[\Lambda,\,\overline{D}_{-\eta}^\star]] =0.\end{eqnarray} Since $[L,\,\Lambda] = (k-n)\,\mbox{Id}$ on $k$-forms, we get $[\overline{D}_{-\eta}^\star,\,[L,\,\Lambda]] = \overline{D}_{-\eta}^\star$.

We then transform the last term on the right of (\ref{Lem:preliminary_BKN-refined_proof_1}) using the Jacobi identity \begin{eqnarray*}0 & = & [\Lambda,\,[L,i\,(\overline{d}_{-\eta}\eta\wedge\cdot\,)^\star]] + [L,\,[i\,(\overline{d}_{-\eta}\eta\wedge\cdot\,)^\star,\,\Lambda]] + [i\,(\overline{d}_{-\eta}\eta\wedge\cdot\,)^\star,\,[\Lambda,\,L]] \\
 & = & [\Lambda,\,[L,i\,(\overline{d}_{-\eta}\eta\wedge\cdot\,)^\star]] - [[L,\,i\,\overline{d}_{-\eta}\eta\wedge\cdot\,],\,\Lambda]^\star - i\,(\overline{d}_{-\eta}\eta\wedge\cdot\,)^\star,\end{eqnarray*} where we used again the identity $[L,\,\Lambda] = (k-n)\,\mbox{Id}$ on $k$-forms to get the last term.

Meanwhile, the last term on the left of (\ref{Lem:preliminary_BKN-refined_proof_1_bis}) can be written as $[L,\,[\Lambda,\,\overline{D}_{-\eta}^\star]] = [[\overline{D}_{-\eta},\,L],\,\Lambda]^\star$, so after putting all these pieces of information together, we see that (\ref{Lem:preliminary_BKN-refined_proof_1}) translates to: \begin{eqnarray*}[\Lambda,\,\tau_\eta] = -i\,[[\overline{D}_{-\eta},\,L],\,\Lambda]^\star & - & i\,\overline{D}_{-\eta}^\star + i\,(\overline{D}_{-\eta} + \overline{\tau}_{-\eta})^\star - \frac{n}{\eta}\,[d_\eta\eta\wedge\cdot\,,\,\Lambda] \\
  & - & \frac{n}{\eta}\,[[L,\,i\,\overline{d}_{-\eta}\eta\wedge\cdot\,],\,\Lambda]^\star - \frac{n}{\eta}\,i\,(\overline{d}_{-\eta}\eta\wedge\cdot\,)^\star.\end{eqnarray*} Since $[d_\eta\eta\wedge\cdot\,,\,\Lambda] = -i\,(\overline{d}_{-\eta}\eta\wedge\cdot\,)^\star$ (see (a) of Corollary \ref{Lem:com_d-eta}) and $[L,\,i\,\overline{d}_{-\eta}\eta\wedge\cdot\,] = 0$ (immediate verification), we get: \begin{eqnarray}\label{Lem:preliminary_BKN-refined_proof_2}[\Lambda,\,\tau_\eta] = -i\,[[\overline{D}_{-\eta},\,L],\,\Lambda]^\star + i\,\overline{\tau}_{-\eta}^\star.\end{eqnarray}

Now, for an arbitrary $(p,\,q)$-form $u$, using (\ref{eqn:D_bar_minus-eta}) we get: \begin{eqnarray*}[\overline{D}_{-\eta},\,L]\,u  & = & \overline{D}_{-\eta}(\omega\wedge u) - \omega\wedge\overline{D}_{-\eta}u \\
  & = & (\partial -\eta\,\bar\partial)\omega\wedge u + \omega\wedge(\partial -\eta\,\bar\partial)u - \frac{q+1}{\eta}\,\partial\eta\wedge(\omega\wedge u) + (q+1)\,\bar\partial\eta\wedge(\omega\wedge u) \\
  & - & \omega\wedge(\partial -\eta\,\bar\partial)\,u + \omega\wedge\frac{q}{\eta}\,\partial\eta\wedge u - q\,\omega\wedge\bar\partial\eta\wedge u \\
  & = & \bigg(\partial -\eta\,\bar\partial - \frac{1}{\eta}\,\partial\eta\wedge \cdot + \bar\partial\eta\wedge \cdot\,\bigg)\,\omega\wedge u = \overline{D}_{-\eta}\omega\wedge u.\end{eqnarray*}

This shows that $[\overline{D}_{-\eta},\,L] = \overline{D}_{-\eta}\omega\wedge\cdot$. Combined with the equality $\overline{\tau}_{-\eta} = -[\overline{D}_{-\eta}\omega\wedge\cdot\,,\,\Lambda]$, this transforms (\ref{Lem:preliminary_BKN-refined_proof_2}) to the desired identity $[\Lambda,\,\tau_\eta] = 2i\,\overline{\tau}_{-\eta}^\star$.

This proves (ii).

\vspace{2ex}

(iii)\, The $\eta$-twisted commutation relation (b) of Proposition \ref{Prop:eta-twisted_commutation-relations} implies \begin{eqnarray}\label{Lem:preliminary_BKN-refined_proof_3}[D_\eta,\,\overline{D}_{-\eta}^\star] = -[D_\eta,\,\overline\tau_{-\eta}^\star] + i\,[D_\eta,\,[\Lambda,\,D_\eta]] - [D_\eta,\,\frac{n}{\eta}\,[i\,d_\eta\eta\wedge\cdot\,,\,\Lambda]].\end{eqnarray}

Meanwhile, the Jacobi identity yields $$-[D_\eta,\,[\Lambda,\,D_\eta]] + [\Lambda,\,[D_\eta,\,D_\eta]] + [D_\eta,\,[D_\eta,\,\Lambda]] =0.$$ Since $[D_\eta,\,D_\eta]=0$ (because $D_\eta^2=0$) and $[D_\eta,\,\Lambda] = -[\Lambda,\,D_\eta]$, we get $[D_\eta,\,[\Lambda,\,D_\eta]]=0$.

Together with $[i\,d_\eta\eta\wedge\cdot\,,\,\Lambda] = (\overline{d}_{-\eta}\eta\wedge\cdot\,)^\star$ (see (a) of Corollary \ref{Lem:com_d-eta}), this transforms 
 (\ref{Lem:preliminary_BKN-refined_proof_3}) into the identity claimed under (iii).

\vspace{2ex}

(iv)\, Applying part (ii) and then the Jacobi identity, we get: \begin{eqnarray}\label{Lem:preliminary_BKN-refined_proof_4}[\overline{D}_{-\eta},\,\overline{\tau}_{-\eta}^\star] = -\frac{i}{2}\,[\overline{D}_{-\eta},\,[\Lambda,\,\tau_\eta]] = -\frac{i}{2}\,[\Lambda,\,[\tau_\eta,\,\overline{D}_{-\eta}]] -\frac{i}{2}\,[\tau_\eta,\,[\overline{D}_{-\eta},\,\Lambda]].\end{eqnarray}

On the other hand, \begin{eqnarray}\label{Lem:preliminary_BKN-refined_proof_5}\nonumber[\tau_\eta,\,\overline{D}_{-\eta}] & \stackrel{(a)}{=} & [\overline{D}_{-\eta},\,\tau_\eta] \stackrel{(b)}{=} [\overline{D}_{-\eta},\,[\Lambda,\,D_\eta\omega\wedge\cdot\,]] \stackrel{(c)}{=} [\Lambda,\,[D_\eta\omega\wedge\cdot\,,\,\overline{D}_{-\eta}]] + [D_\eta\omega\wedge\cdot\,,\,[\overline{D}_{-\eta},\, \Lambda]] \\ 
 & \stackrel{(d)}{=} & [\Lambda,\,\overline{D}_{-\eta}D_\eta\omega\wedge\cdot\,] - i\,[D_\eta\omega\wedge\cdot,\, D_\eta^\star + \tau_\eta^\star] - \frac{n}{\eta}\,\bigg[D_\eta\omega\wedge\cdot,\,[\overline{d}_{-\eta}\eta\wedge\cdot\,,\,\Lambda]\bigg],\end{eqnarray}

\noindent where (a) follows from $\tau_\eta$ and $\overline{D}_{-\eta}$ being operators of odd degrees, (b) follows from the definition of $\tau_\eta$, (c) follows from the Jacobi identity, while the sum of the last two terms in (d) follows from \begin{eqnarray*}[\overline{D}_{-\eta},\, \Lambda] = -i\,( D_\eta^\star + \tau_\eta^\star) - \frac{n}{\eta}\,[\overline{d}_{-\eta}\eta\wedge\cdot\,,\,\Lambda]\end{eqnarray*} which is a consequence of the $\eta$-twisted commutation relation (a) of Proposition \ref{Prop:eta-twisted_commutation-relations}.

The computation of the first term in (d) of (\ref{Lem:preliminary_BKN-refined_proof_5}) runs, for any $(p,\,q)$-form $u$, as follows.

\vspace{1ex}

$\bullet$\, First, \begin{eqnarray*}[D_\eta\omega\wedge\cdot\,,\,\overline{D}_{-\eta}]\,u = D_\eta\omega\wedge\overline{D}_{-\eta}u + \overline{D}_{-\eta}(D_\eta\omega\wedge u) = D_\eta\omega\wedge\overline{D}_{-\eta}u + \overline{D}_{-\eta}(D_\eta^{1,\,0}\omega\wedge u) + \overline{D}_{-\eta}(D_\eta^{0,\,1}\omega\wedge u).\end{eqnarray*}

\vspace{1ex}

$\bullet$\, Then, we see that \begin{eqnarray*}\overline{D}_{-\eta}(D_\eta^{1,\,0}\omega\wedge u) = (\partial - \eta\,\bar\partial)(D_\eta^{1,\,0}\omega\wedge u) - \frac{q+1}{\eta}\,\partial\eta\wedge D_\eta^{1,\,0}\omega\wedge u + (q+1)\,\bar\partial\eta\wedge D_\eta^{1,\,0}\omega\wedge u\end{eqnarray*} and that \begin{eqnarray*}\overline{D}_{-\eta}(D_\eta^{0,\,1}\omega\wedge u) = (\partial - \eta\,\bar\partial)(D_\eta^{0,\,1}\omega\wedge u) - \frac{q+2}{\eta}\,\partial\eta\wedge D_\eta^{0,\,1}\omega\wedge u + (q+2)\,\bar\partial\eta\wedge D_\eta^{0,\,1}\omega\wedge u.\end{eqnarray*}

\vspace{1ex}

$\bullet$\, Putting these pieces of information together, we deduce that \begin{eqnarray*}[D_\eta\omega\wedge\cdot\,,\,\overline{D}_{-\eta}]\,u = D_\eta\omega\wedge\overline{D}_{-\eta}u & + & (\partial - \eta\,\bar\partial)\,(D_\eta\omega)\wedge u   \\
   & - & D_\eta^{1,\,0}\omega\wedge\bigg((\partial - \eta\,\bar\partial)\,u - \frac{q+1}{\eta}\,\partial\eta\wedge u + (q+1)\,\bar\partial\eta\wedge u\bigg) \\
  & - & D_\eta^{0,\,1}\omega\wedge\bigg((\partial - \eta\,\bar\partial)\,u - \frac{q+2}{\eta}\,\partial\eta\wedge u + (q+2)\,\bar\partial\eta\wedge u\bigg).\end{eqnarray*}

Hence, \begin{eqnarray*}[D_\eta\omega\wedge\cdot\,,\,\overline{D}_{-\eta}]\,u & = & D_\eta\omega\wedge\overline{D}_{-\eta}u + (\partial - \eta\,\bar\partial)\,(D_\eta^{1,\,0}\omega)\wedge u + (\partial - \eta\,\bar\partial)\,(D_\eta^{0,\,1}\omega)\wedge u \\
  & - & D_\eta\omega\wedge\overline{D}_{-\eta}u +  D_\eta^{1,\,0}\omega\wedge\bigg(\frac{1}{\eta}\,\partial\eta\wedge u -\bar\partial\eta\wedge u\bigg) +  D_\eta^{0,\,1}\omega\wedge\bigg(\frac{2}{\eta}\,\partial\eta\wedge u - 2\,\bar\partial\eta\wedge u\bigg) \\
  & = & \overline{D}_{-\eta}(D_\eta^{1,\,0}\omega)\wedge u + \overline{D}_{-\eta}(D_\eta^{0,\,1}\omega)\wedge u.\end{eqnarray*}

\vspace{1ex}

$\bullet$\, We conclude that the formula \begin{eqnarray}\label{Lem:preliminary_BKN-refined_proof_6}[D_\eta\omega\wedge\cdot\,,\,\overline{D}_{-\eta}] = \overline{D}_{-\eta}D_\eta\omega\wedge\cdot\end{eqnarray} holds in every (bi-)degree.

  This finishes the computation of the first term in (d) of (\ref{Lem:preliminary_BKN-refined_proof_5}).

  \vspace{2ex}

  Taking the bracket with $\Lambda$ in (\ref{Lem:preliminary_BKN-refined_proof_5}), we get: \begin{eqnarray}\label{Lem:preliminary_BKN-refined_proof_7}[\Lambda,\,[\tau_\eta,\,\overline{D}_{-\eta}]] = [\Lambda,\,[\Lambda,\,\overline{D}_{-\eta}D_\eta\omega\wedge\cdot\,]] -i\,[\Lambda,\, [D_\eta\omega\wedge\cdot\,,\, D_\eta^\star + \tau_\eta^\star]] - \frac{n}{\eta}\,\bigg[\Lambda,\,\bigg[D_\eta\omega\wedge\cdot,\,[\overline{d}_{-\eta}\eta\wedge\cdot\,,\,\Lambda]\bigg]\bigg].\end{eqnarray}

 $\bullet$ Now, we compute the second term on the right of (\ref{Lem:preliminary_BKN-refined_proof_7}) starting from the Jacobi formula: \begin{eqnarray}\label{Lem:preliminary_BKN-refined_proof_8}\nonumber[\Lambda,\, [D_\eta\omega\wedge\cdot\,,\, D_\eta^\star + \tau_\eta^\star]] & = & - [D_\eta\omega\wedge\cdot\,,\,[D_\eta^\star + \tau_\eta^\star,\,\Lambda]] + [D_\eta^\star + \tau_\eta^\star,\,[\Lambda,\,D_\eta\omega\wedge\cdot\,]] \\
    & = & - [D_\eta\omega\wedge\cdot\,,\, [L,\,D_\eta+\tau_\eta]^\star] + [D_\eta^\star + \tau_\eta^\star,\,\tau_\eta].\end{eqnarray}

  Moreover, for any $(p,\,q)$-form $u$, we have: \begin{eqnarray*}[L,\,D_\eta]\,u & = & \omega\wedge D_\eta u - D_\eta(\omega\wedge u) \\
    & = & \omega\wedge D_\eta u - (\eta\,\partial + \bar\partial)(\omega\wedge u) + (p+1)\,\partial\eta\wedge\omega\wedge u + \frac{p+1}{\eta}\,\bar\partial\eta\wedge\omega\wedge u \\
    & = & \omega\wedge D_\eta u - \omega\wedge D_\eta u - \bigg((\eta\,\partial + \bar\partial) - \partial\eta\wedge\cdot - \frac{1}{\eta}\,\bar\partial\eta\wedge\cdot\,\bigg)\,\omega\wedge u \\
    & = & -D_\eta\omega\wedge u.\end{eqnarray*}

  We have thus got the formula \begin{eqnarray*}[L,\,D_\eta] = - D_\eta\omega\wedge\cdot,\end{eqnarray*} which, together with the one proved under (i), yields \begin{eqnarray*}[L,\,D_\eta+\tau_\eta] = - D_\eta\omega\wedge\cdot + [L,\,\tau_\eta] = 2\,D_\eta\omega\wedge\cdot\,.\end{eqnarray*}

  Combining this with (\ref{Lem:preliminary_BKN-refined_proof_8}), we transform (\ref{Lem:preliminary_BKN-refined_proof_7}) to \begin{eqnarray}\label{Lem:preliminary_BKN-refined_proof_9}\nonumber[\Lambda,\,[\tau_\eta,\,\overline{D}_{-\eta}]] = [\Lambda,\,[\Lambda,\,\overline{D}_{-\eta}D_\eta\omega\wedge\cdot\,]] & + & 2i\,[D_\eta\omega\wedge\cdot\,,\, (D_\eta\omega\wedge\cdot\,)^\star]] -i\,[D_\eta^\star + \tau_\eta^\star,\,\tau_\eta]   \\
      & - & \frac{n}{\eta}\,\bigg[\Lambda,\,\bigg[D_\eta\omega\wedge\cdot,\,[\overline{d}_{-\eta}\eta\wedge\cdot\,,\,\Lambda]\bigg]\bigg].\end{eqnarray}

  This is the end of the computation of the last but one term in (\ref{Lem:preliminary_BKN-refined_proof_4}).

  \vspace{1ex}

  $\bullet$\, Using (\ref{Lem:preliminary_BKN-refined_proof_9}), (\ref{Lem:preliminary_BKN-refined_proof_4}) becomes \begin{eqnarray}\label{Lem:preliminary_BKN-refined_proof_10}[\overline{D}_{-\eta},\,\overline{\tau}_{-\eta}^\star]  =  -\frac{i}{2}\,[\Lambda,\,[\Lambda,\,\overline{D}_{-\eta}D_\eta\omega\wedge\cdot\,]] & + & [D_\eta\omega\wedge\cdot\,,\, (D_\eta\omega\wedge\cdot\,)^\star] -\frac{1}{2}\,[D_\eta^\star + \tau_\eta^\star,\,\tau_\eta] \\
 \nonumber   & - & \frac{i}{2}\,[\tau_\eta,\,[\overline{D}_{-\eta},\,\Lambda]] + \frac{n}{2\eta}\,\bigg[\Lambda,\,\bigg[D_\eta\omega\wedge\cdot,\,[i\overline{d}_{-\eta}\eta\wedge\cdot\,,\,\Lambda]\bigg]\bigg].\end{eqnarray}

  Now, expressing $i\,[\overline{D}_{-\eta},\,\Lambda]$ by means of the $\eta$-twisted commutation relation (a) of Proposition \ref{Prop:eta-twisted_commutation-relations}, the last but one term in (\ref{Lem:preliminary_BKN-refined_proof_10}) reads: \begin{eqnarray*}-\frac{i}{2}\,[\tau_\eta,\,[\overline{D}_{-\eta},\,\Lambda]] = -\frac{1}{2}\,[\tau_\eta,\,D_\eta^\star + \tau_\eta^\star] + \frac{n}{2\eta}\,[\tau_\eta,\,[i\overline{d}_{-\eta}\eta\wedge\cdot\,,\,\Lambda]].\end{eqnarray*}

  Meanwhile, the Jacobi identity expresses the last term in (\ref{Lem:preliminary_BKN-refined_proof_10}) as \begin{eqnarray*}\frac{n}{2\eta}\,\bigg[\Lambda,\,\bigg[D_\eta\omega\wedge\cdot,\,[i\overline{d}_{-\eta}\eta\wedge\cdot\,,\,\Lambda]\bigg]\bigg] = - \frac{n}{2\eta}\,\bigg[D_\eta\omega\wedge\cdot,\,\bigg[[i\overline{d}_{-\eta}\eta\wedge\cdot\,,\,\Lambda],\,\Lambda\bigg]\bigg] + \frac{n}{2\eta}\,\bigg[[i\overline{d}_{-\eta}\eta\wedge\cdot\,,\,\Lambda],\,\bigg[\Lambda,\,D_\eta\omega\wedge\cdot\bigg]\bigg].\end{eqnarray*}

  Using these equalities, (\ref{Lem:preliminary_BKN-refined_proof_10}) reduces to \begin{eqnarray}\label{Lem:preliminary_BKN-refined_proof_11}-[\overline{D}_{-\eta},\,\overline{\tau}_{-\eta}^\star]  = [D_\eta^\star + \tau_\eta^\star,\,\tau_\eta] & + & \frac{i}{2}\,[\Lambda,\,[\Lambda,\,\overline{D}_{-\eta}D_\eta\omega\wedge\cdot\,]] - [D_\eta\omega\wedge\cdot\,,\, (D_\eta\omega\wedge\cdot\,)^\star] \\
\nonumber    & - & \frac{n}{\eta}\,\bigg[\tau_\eta,\,[i\overline{d}_{-\eta}\eta\wedge\cdot\,,\,\Lambda]\bigg] + \frac{n}{2\eta}\,\bigg[D_\eta\omega\wedge\cdot\,,\,\bigg[[i\overline{d}_{-\eta}\eta\wedge\cdot\,,\,\Lambda],\,\Lambda\bigg]\bigg].\end{eqnarray}

  Now, in the last but one term, (b) of Corollary \ref{Lem:com_d-eta} ensures that $[i\overline{d}_{-\eta}\eta\wedge\cdot\,,\,\Lambda] = -(d_\eta\eta\wedge\cdot\,)^\star$. Meanwhile, the last term in (\ref{Lem:preliminary_BKN-refined_proof_11}) vanishes since \begin{eqnarray*}\bigg[[i\overline{d}_{-\eta}\eta\wedge\cdot\,,\,\Lambda],\,\Lambda\bigg] = -[(d_\eta\eta\wedge\cdot\,)^\star,\,\Lambda] = - [L,\,d_\eta\eta\wedge\cdot\,]^\star = 0,\end{eqnarray*} the last vanishing being a consequence of the equalities \begin{eqnarray*}[L,\,d_\eta\eta\wedge\cdot\,]\,u = \omega\wedge(d_\eta\eta\wedge u) - d_\eta\eta\wedge(\omega\wedge u) = 0\end{eqnarray*} that hold for any form $u$.

    Thus, using these observations and the notation $S_\omega^{(\eta)}$ spelt out in the statement, we see that (\ref{Lem:preliminary_BKN-refined_proof_11}) translates to \begin{eqnarray*}\label{Lem:preliminary_BKN-refined_proof_12}-[\overline{D}_{-\eta},\,\overline{\tau}_{-\eta}^\star]  = [D_\eta^\star + \tau_\eta^\star,\,\tau_\eta] + S_\omega^{(\eta)} + \frac{n}{\eta}\,[\tau_\eta,\,(d_\eta\eta\wedge\cdot\,)^\star].\end{eqnarray*} It remains to add $[D_\eta,\,D_\eta^\star] + [D_\eta,\,\tau_\eta^\star] = [D_\eta,\,D_\eta^\star + \tau_\eta^\star]$ to either side of this equality to get the desired equality (iv) of Lemma \ref{Lem:preliminary_BKN-refined}.  \hfill $\Box$

\vspace{3ex}

We can now give the main result of this section.

\begin{The}\label{The:eta-BKN} Let $(X,\,\omega)$ be a complex Hermitian manifold with $\mbox{dim}_\C X=n\geq 2$.  For any $C^\infty$ function $\eta$ on $X$ such that $\eta>0$ or $\eta<0$ and any $k\in\{0,\dots , 2n\}$, the following {\bf refined $\eta$-Bochner-Kodaira-Nakano ($\eta$-BKN) identity} holds on the $C^\infty$ forms on $X$: \begin{eqnarray}\label{eqn:eta-BKN}\nonumber\Delta_\eta = [\overline{D}_{-\eta} + \overline\tau_{-\eta},\,\overline{D}_{-\eta}^\star + \overline\tau_{-\eta}^\star] & + & T_\omega^{(\eta)} + i\,[[D_\eta,\,\overline{D}_{-\eta}],\,\Lambda] \\
    & + & n\,\bigg[\overline{D}_{-\eta} + \overline\tau_{-\eta},\,\frac{1}{\eta}\,(\overline{d}_{-\eta}\eta\wedge\cdot\,)^\star\bigg] -n\,\bigg[D_\eta,\,\frac{1}{\eta}\,(d_\eta\eta\wedge\cdot\,)^\star\bigg],\end{eqnarray} where $T_\omega^{(\eta)}$ is the zero-th order operator defined by \begin{eqnarray*}T_\omega^{(\eta)}:= \overline{S_\omega^{(-\eta)}} = -\frac{i}{2}\,[\Lambda,\,[\Lambda,\,D_\eta\overline{D}_{-\eta}\omega\wedge\cdot\,]] - [\overline{D}_{-\eta}\omega\wedge\cdot,\, (\overline{D}_{-\eta}\omega\wedge\cdot\,)^\star].\end{eqnarray*}

\end{The}  

\noindent {\it Proof.} Putting together the rough $\eta$-BKN formula (\ref{eqn:eta-BKN_rough}) and equality (iv) of Lemma \ref{Lem:preliminary_BKN-refined}, we get: \begin{eqnarray*} & & \Delta_\eta + [D_\eta + \tau_\eta,\,D_\eta^\star + \tau_\eta^\star] + S_\omega^{(\eta)} + \frac{n}{\eta}\,[\tau_\eta,\,(d_\eta\eta\wedge\cdot\,)^\star] \\
  & = & \overline{\Delta}_{-\eta} + i\,[[D_\eta,\,\overline{D}_{-\eta}],\,\Lambda]  + [\overline{D}_{-\eta},\,\overline\tau_{-\eta}^\star] - [D_\eta,\,\tau_\eta^\star] + [D_\eta,\,D_\eta^\star] + [D_\eta,\,\tau_\eta^\star] -[\overline{D}_{-\eta},\,\overline\tau_{-\eta}^\star] \\
  &  & + n\,\bigg[D_\eta,\,\frac{1}{\eta}\,[i\overline{d}_{-\eta}\eta\wedge\cdot\,,\,\Lambda]\bigg] + n\,\bigg[\overline{D}_{-\eta},\,\frac{1}{\eta}\,[id_\eta\eta\wedge\cdot\,,\,\Lambda]\bigg].\end{eqnarray*}

Since $[D_\eta,\,D_\eta^\star] = \Delta_\eta$ and since the terms $[D_\eta,\,\tau_\eta^\star]$ and $[\overline{D}_{-\eta},\,\overline\tau_{-\eta}^\star]$ reoccur with the opposite signs, this equality reduces to \begin{eqnarray*}\overline{\Delta}_{-\eta} & = & [D_\eta + \tau_\eta,\,D_\eta^\star + \tau_\eta^\star] + S_\omega^{(\eta)} - i\,[[D_\eta,\,\overline{D}_{-\eta}],\,\Lambda] \\
  & + & n\,\bigg[\tau_\eta,\,\frac{1}{\eta}\,(d_\eta\eta\wedge\cdot\,)^\star\bigg] + n\,\bigg[D_\eta,\,\frac{1}{\eta}\,(d_\eta\eta\wedge\cdot\,)^\star\bigg] - n\,\bigg[\overline{D}_{-\eta},\,\frac{1}{\eta}\,(\overline{d}_{-\eta}\eta\wedge\cdot\,)^\star\bigg].\end{eqnarray*} 

Identity (\ref{eqn:eta-BKN}) follows from this by taking conjugates and replacing $\eta$ with $-\eta$. \hfill $\Box$

\vspace{6ex}

We now compute the operator $[D_\eta,\,\overline{D}_{-\eta}]$ that features in identity (\ref{eqn:eta-BKN}).

\begin{Prop}\label{Prop:curvature_computation} Let $(X,\,\omega)$ be a complex Hermitian manifold with $\mbox{dim}_\C X=n\geq 2$.  For any $C^\infty$ function $\eta$ on $X$ such that $\eta>0$ or $\eta<0$ and any bidegree $(p,\,q)$, the following identities hold: \begin{eqnarray}\label{eqn:curvature_computation_first}\nonumber[D_\eta,\,\overline{D}_{-\eta}]\,(u^{p,\,q}) & = & \bigg(\eta - \frac{1}{\eta}\bigg)\,\bar\partial\eta\wedge\partial u^{p,\,q} - \bigg(\eta + \frac{1}{\eta}\bigg)\,\partial\eta\wedge\bar\partial u^{p,\,q} \\
    \nonumber  & + & 2\,\frac{p-q}{\eta^2}\,\partial\eta\wedge\bar\partial\eta\wedge u^{p,\,q} + (q-p)\,\bigg(\eta + \frac{1}{\eta}\bigg)\,\partial\bar\partial\eta\wedge u^{p,\,q},\end{eqnarray} and \begin{eqnarray}\label{eqn:curvature_computation_second}\nonumber i\,\bigg[[D_\eta,\,\overline{D}_{-\eta}],\,\Lambda\bigg]\,(u^{p,\,q}) & = & \bigg(\eta - \frac{1}{\eta}\bigg)\,\bigg([i\bar\partial\eta\wedge\cdot\,,\,\Lambda]\,\partial u^{p,\,q} + i\bar\partial\eta\wedge[\partial,\,\Lambda]\,u^{p,\,q}\bigg) \\
 & - & \bigg(\eta + \frac{1}{\eta}\bigg)\,\bigg([i\partial\eta\wedge\cdot\,,\,\Lambda]\,\bar\partial u^{p,\,q} + i\partial\eta\wedge[\bar\partial,\,\Lambda]\,u^{p,\,q}\bigg) \\
\nonumber & + & 2\frac{p-q}{\eta^2}\,\bigg[i\partial\eta\wedge\bar\partial\eta\wedge\cdot\,,\,\Lambda\bigg]\,u^{p,\,q} + (q-p)\,\bigg(\eta + \frac{1}{\eta}\bigg)\,\bigg[i\partial\bar\partial\eta\wedge\cdot\,,\,\Lambda\bigg]\,u^{p,\,q}\end{eqnarray} for any $C^\infty$ $(p,\,q)$-form $u^{p,\,q}$.

\end{Prop}  

\noindent {\it Proof.} We compute separately the two terms in the sum $[D_\eta,\,\overline{D}_{-\eta}]\,(u^{p,\,q}) = (D_\eta\overline{D}_{-\eta})\,(u^{p,\,q}) + (\overline{D}_{-\eta}D_\eta)\,(u^{p,\,q})$.

\vspace{1ex}

$\bullet$ We use (\ref{eqn:D_eta_explicit}) to get the first line and then (\ref{eqn:D_bar_minus-eta}) to get the next four lines below: \begin{eqnarray*}(\overline{D}_{-\eta}D_\eta)\,(u^{p,\,q}) & = & \overline{D}_{-\eta}(\eta\,\partial u^{p,\,q}) + \overline{D}_{-\eta}(\bar\partial u^{p,\,q}) -p\,\overline{D}_{-\eta}(\partial\eta\wedge u^{p,\,q}) - p\,\overline{D}_{-\eta}\bigg(\frac{1}{\eta}\,\,\bar\partial\eta\wedge u^{p,\,q}\bigg) \\
   & = & \partial(\eta\,\partial u^{p,\,q}) -\eta\bar\partial(\eta\,\partial u^{p,\,q}) - \frac{q}{\eta}\,\eta\,\partial\eta\wedge\partial u^{p,\,q} + q\eta\,\bar\partial\eta\wedge\partial u^{p,\,q}  \\
    & + & \partial\bar\partial u^{p,\,q} - \frac{q+1}{\eta}\,\partial\eta\wedge\bar\partial u^{p,\,q} + (q+1)\,\bar\partial\eta\wedge\bar\partial u^{p,\,q}   \\
    & - & p\partial(\partial\eta\wedge u^{p,\,q}) + p\eta\,\bar\partial(\partial\eta\wedge u^{p,\,q}) - pq\,\bar\partial\eta\wedge\partial\eta\wedge u^{p,\,q} \\
   & - & p\partial\bigg(\frac{1}{\eta}\,\bar\partial\eta\wedge u^{p,\,q}\bigg) + p\eta\,\bar\partial\bigg(\frac{1}{\eta}\,\bar\partial\eta\wedge u^{p,\,q}\bigg) + \frac{p(q+1)}{\eta^2}\,\partial\eta\wedge\bar\partial\eta\wedge u^{p,\,q}.\end{eqnarray*}

This leads, through straightforward computations, to \begin{eqnarray}\label{eqn:curvature_computation_proof_1}\nonumber(\overline{D}_{-\eta}D_\eta)\,(u^{p,\,q}) & = & p\,\bigg(q + \frac{q+2}{\eta^2}\bigg)\,\partial\eta\wedge\bar\partial\eta\wedge u^{p,\,q} - p\,\bigg(\eta + \frac{1}{\eta}\bigg)\,\partial\bar\partial\eta\wedge u^{p,\,q} \\
  \nonumber & + &  (p-q+1)\,\partial\eta\wedge\partial u^{p,\,q} + \bigg(\frac{p}{\eta} + (q+1)\,\eta\bigg)\,\bar\partial\eta\wedge\partial u^{p,\,q} \\
  \nonumber & - & \bigg(p\eta + \frac{q+1}{\eta}\bigg)\,\partial\eta\wedge\bar\partial u^{p,\,q} + (q-p+1)\,\bar\partial\eta\wedge\bar\partial u^{p,\,q} \\
 & + & (1+\eta^2)\,\partial\bar\partial u^{p,\,q}.\end{eqnarray}

$\bullet$ Meanwhile, we use (\ref{eqn:D_bar_minus-eta}) to get the first line and then (\ref{eqn:D_eta_explicit}) to get the next four lines below: \begin{eqnarray*}(D_\eta\overline{D}_{-\eta})\,(u^{p,\,q}) & = & D_\eta(\partial u^{p,\,q}) - D_\eta(\eta\,\bar\partial u^{p,\,q}) - q\,D_\eta\bigg(\frac{1}{\eta}\,\partial\eta\wedge u^{p,\,q}\bigg) + q\,D_\eta(\bar\partial\eta\wedge u^{p,\,q}) \\
  & = & \bar\partial\partial u^{p,\,q} - (p+1)\,\partial\eta\wedge\partial u^{p,\,q} - \frac{p+1}{\eta}\,\bar\partial\eta\wedge\partial u^{p,\,q} \\
  & - & \eta\,\partial(\eta\,\bar\partial u^{p,\,q}) - \bar\partial(\eta\,\bar\partial u^{p,\,q}) + p\eta\,\partial\eta\wedge\bar\partial u^{p,\,q} + p\,\bar\partial\eta\wedge\bar\partial u^{p,\,q} \\
  & - & q\eta\,\partial\bigg(\frac{1}{\eta}\,\partial\eta\wedge u^{p,\,q}\bigg) - q\,\bar\partial\bigg(\frac{1}{\eta}\,\partial\eta\wedge u^{p,\,q}\bigg) + \frac{(p+1)q}{\eta^2}\,\bar\partial\eta\wedge\partial\eta\wedge u^{p,\,q} \\
  & + & q\eta\,\partial(\bar\partial\eta\wedge u^{p,\,q}) + q\,\bar\partial(\bar\partial\eta\wedge u^{p,\,q}) - pq\,\partial\eta\wedge\bar\partial\eta\wedge u^{p,\,q}.\end{eqnarray*}

This leads, through straightforward computations, to \begin{eqnarray}\label{eqn:curvature_computation_proof_2}\nonumber (D_\eta\overline{D}_{-\eta})\,(u^{p,\,q}) & = & - \bigg(\frac{(p+2)q}{\eta^2} + pq\bigg)\,\partial\eta\wedge\bar\partial\eta\wedge u^{p,\,q} + q\,\bigg(\eta + \frac{1}{\eta}\bigg)\,\partial\bar\partial\eta\wedge u^{p,\,q} \\
  \nonumber & + &  (q-p-1)\,\partial\eta\wedge\partial u^{p,\,q} - \bigg(\frac{p+1}{\eta} + q\,\eta\bigg)\,\bar\partial\eta\wedge\partial u^{p,\,q} \\
  \nonumber & + & \bigg((p-1)\,\eta + \frac{q}{\eta}\bigg)\,\partial\eta\wedge\bar\partial u^{p,\,q} + (p-q-1)\,\bar\partial\eta\wedge\bar\partial u^{p,\,q} \\      & - & (1+\eta^2)\,\partial\bar\partial u^{p,\,q}.  \end{eqnarray} 
   
\vspace{2ex}

Adding up (\ref{eqn:curvature_computation_proof_1}) and (\ref{eqn:curvature_computation_proof_2}) yields the former desired equality.

\vspace{1ex}

Now, starting from $i\,\bigg[[D_\eta,\,\overline{D}_{-\eta}],\,\Lambda\bigg]\,(u^{p,\,q}) = i\,[D_\eta,\,\overline{D}_{-\eta}]\,(\Lambda u^{p,\,q}) -\Lambda\bigg(i\,[D_\eta,\,\overline{D}_{-\eta}]\,(u^{p,\,q})\bigg)$ and using the equality we just proved, we get \begin{eqnarray}\label{eqn:curvature_computation_second_proof_1}\nonumber i\,\bigg[[D_\eta,\,\overline{D}_{-\eta}],\,\Lambda\bigg]\,(u^{p,\,q}) & = & \bigg(\eta - \frac{1}{\eta}\bigg)\,\bigg(i\bar\partial\eta\wedge\partial\Lambda\, u^{p,\,q} - \Lambda(i\bar\partial\eta\wedge\partial u^{p,\,q})\bigg) \\
  \nonumber & - & \bigg(\eta + \frac{1}{\eta}\bigg)\,\bigg(i\partial\eta\wedge\bar\partial\Lambda\, u^{p,\,q} - \Lambda(i\partial\eta\wedge\bar\partial u^{p,\,q})\bigg) \\
 \nonumber & + & 2\frac{p-q}{\eta^2}\,\bigg[i\partial\eta\wedge\bar\partial\eta\wedge\cdot\,,\,\Lambda\bigg]\,u^{p,\,q} + (q-p)\,\bigg(\eta + \frac{1}{\eta}\bigg)\,\bigg[i\partial\bar\partial\eta\wedge\cdot\,,\,\Lambda\bigg]\,u^{p,\,q}.\end{eqnarray}

It remains to notice that on the first line on the r.h.s. above we can write: \begin{eqnarray*}i\bar\partial\eta\wedge\partial\Lambda\, u^{p,\,q} - \Lambda(i\bar\partial\eta\wedge\partial u^{p,\,q}) = [i\bar\partial\eta\wedge\cdot\,,\,\Lambda]\,\partial u^{p,\,q} + i\bar\partial\eta\wedge[\partial,\,\Lambda]\,u^{p,\,q},\end{eqnarray*} while on the second line on the r.h.s. above we can write: \begin{eqnarray*}i\partial\eta\wedge\bar\partial\Lambda\, u^{p,\,q} - \Lambda(i\partial\eta\wedge\bar\partial u^{p,\,q}) = [i\partial\eta\wedge\cdot\,,\,\Lambda]\,\bar\partial u^{p,\,q} + i\partial\eta\wedge[\bar\partial,\,\Lambda]\,u^{p,\,q}\end{eqnarray*} to get (\ref{eqn:curvature_computation_second}). \hfill $\Box$

\section{Twisted commutation relations for $D_\eta^{1,\,0}$ and $D_\eta^{0,\,1}$}\label{section:twisted-commutation_1-0_0-1} The setting is the same as in the previous section $\S$\ref{section:twisted-commutation_D-eta}. Besides the differential operators defined by: \begin{eqnarray*}D_\eta^{1,\,0} = \eta\partial -p\partial\eta\wedge\cdot  \hspace{3ex} \mbox{and} \hspace{3ex} D_\eta^{0,\,1} = \bar\partial -\frac{p}{\eta}\,\bar\partial\eta\wedge\cdot\, ,\end{eqnarray*} on $(p,\,q)$-forms, we now also consider the operators $d'_\eta = \eta\partial$ and $d''_\eta = \bar\partial$, as well as \begin{eqnarray*}\tau_\eta^{1,\,0}:=[\Lambda,\,D_\eta^{1,\,0}\omega\wedge\cdot\,] = \eta\tau - \bigg[\Lambda,\,\partial\eta\wedge\omega\wedge\cdot\,\bigg]  \hspace{3ex} \mbox{and} \hspace{3ex} \tau_\eta^{0,\,1}:=[\Lambda,\,D_\eta^{0,\,1}\omega\wedge\cdot\,] = \overline\tau - \bigg[\Lambda,\,\frac{1}{\eta}\,\bar\partial\eta\wedge\omega\wedge\cdot\,\bigg].\end{eqnarray*}

Thus, we have: \begin{eqnarray*}D_\eta = D_\eta^{1,\,0} + D_\eta^{0,\,1}; \hspace{5ex} d_\eta = d'_\eta + d''_\eta; \hspace{5ex} \tau_\eta = \tau_\eta^{1,\,0} + \tau_\eta^{0,\,1}.\end{eqnarray*}

\begin{Lem}\label{Lem:formulae_D_eta_1-0_0-1_tau_eta_1-0_0-1} (i)\, The following formulae hold in any bidegree: \begin{eqnarray}\label{eqn:formulae_D_eta_1-0_0-1_tau_eta_1-0_0-1_no1}\tau = \overline{\tau_\eta^{0,\,1}} + \bigg[\Lambda,\,\frac{1}{\eta}\,\partial\eta\wedge\omega\wedge\cdot\,\bigg] \hspace{5ex}\mbox{and}\hspace{5ex} \overline\tau = \tau_\eta^{0,\,1} + \bigg[\Lambda,\,\frac{1}{\eta}\,\bar\partial\eta\wedge\omega\wedge\cdot\,\bigg].\end{eqnarray}

 \vspace{1ex}

 (ii)\, For any $p,q$, the following formulae hold in bidegree $(p,\,q)$: \begin{eqnarray}\label{eqn:formulae_D_eta_1-0_0-1_tau_eta_1-0_0-1_no2}\nonumber D_\eta^{0,\,1} + \tau_\eta^{0,\,1} & = & (\bar\partial + \bar\tau) - \frac{p}{\eta}\,\bar\partial\eta\wedge\cdot - \bigg[\Lambda,\,\frac{1}{\eta}\,\bar\partial\eta\wedge\omega\wedge\cdot\,\bigg] \\
 \overline{D_\eta^{0,\,1}} + \overline{\tau_\eta^{0,\,1}} & = & (\partial + \tau) - \frac{q}{\eta}\,\partial\eta\wedge\cdot - \bigg[\Lambda,\,\frac{1}{\eta}\,\partial\eta\wedge\omega\wedge\cdot\,\bigg].\end{eqnarray}

\end{Lem}

\noindent {\it Proof.} Immediate verification.  \hfill $\Box$

\vspace{2ex}

Equating the bidegrees in $D_\eta^2=0$, we get: \begin{eqnarray}\label{eqn:D_eta_1-0_0-1_prop}(i)\,(D_\eta^{1,\,0})^2=0; \hspace{3ex} (ii)\,D_\eta^{1,\,0}D_\eta^{0,\,1} = - D_\eta^{0,\,1}D_\eta^{1,\,0};  \hspace{3ex} (iii)\,(D_\eta^{0,\,1})^2=0.\end{eqnarray} The anti-commutation of $D_\eta^{1,\,0}$ and $D_\eta^{0,\,1}$ is an advantage these operators have over $D_\eta$ and $\overline{D}_{-\eta}$ (that do not anti-commute -- see Proposition \ref{Prop:curvature_computation}) considered in $\S.$\ref{section:twisted-commutation_1-0_0-1}.

\vspace{2ex}

We also consider the Laplacians: \begin{eqnarray*}\Delta''_\eta:=[D_\eta^{0,\,1},\,(D_\eta^{0,\,1})^\star]  \hspace{3ex}\mbox{and}\hspace{3ex}  \Delta'_\eta:=[\overline{D_\eta^{0,\,1}},\,\overline{D_\eta^{0,\,1}}^{\,\star}],\end{eqnarray*} where the conjugate operator $\overline{D_\eta^{0,\,1}}$ is defined by requiring $\overline{D_\eta^{0,\,1}}u=\overline{D_\eta^{0,\,1}\bar{u}}$ for every form $u$.

An immediate computation shows that, in bidegree $(p,\,q)$, we have: \begin{eqnarray}\label{eqn:conjugate_D-01_eta}\overline{D_\eta^{0,\,1}} = \partial - \frac{q}{\eta}\,\partial\eta\wedge\cdot\,, \hspace{7ex}\mbox{hence also}\hspace{7ex} \overline{D^{0,\,1}_\eta}^{\,\star} = \partial^\star - \frac{q}{\eta}\,(\partial\eta\wedge\cdot\,)^\star.\end{eqnarray}


\begin{Lem}\label{Lem:conjugate-star_D-01_eta_conjugate-star_tau-01_eta} The following formula holds in any bidegree $(p,\,q)$: \begin{eqnarray}\label{eqn:conjugate-star_D-01_eta_conjugate-star_tau-01_eta}\overline{D^{0,\,1}_\eta}^{\,\star} + \overline{\tau^{0,\,1}_\eta}^{\,\star} = (\partial^\star + \tau^\star) - \frac{q}{\eta}\,(\partial\eta\wedge\cdot\,)^\star - \bigg[\bigg(\frac{1}{\eta}\,\partial\eta\wedge\omega\wedge\cdot\,\bigg)^\star,\,\omega\wedge\cdot\,\bigg].\end{eqnarray}

\end{Lem}

\noindent {\it Proof.} It follows at once by conjugating the latter equality in (\ref{eqn:formulae_D_eta_1-0_0-1_tau_eta_1-0_0-1_no2}).  \hfill $\Box$

\vspace{3ex}

On the other hand, the identity $D_\eta = \theta_\eta d\theta_\eta^{-1}$ (see (\ref{eqn:D_eta})) yields: \begin{eqnarray}\label{eqn:D_eta_1-0_0-1}D^{1,\,0}_\eta = \theta_\eta\partial\theta_\eta^{-1} \hspace{3ex}\mbox{and}\hspace{3ex} D^{0,\,1}_\eta = \theta_\eta\bar\partial\theta_\eta^{-1}.\end{eqnarray}

In particular, for every bidegree $(p,\,q)$, the $D^{1,\,0}_\eta$-cohomology and $D^{0,\,1}_\eta$-cohomology spaces: \begin{eqnarray}\label{eqn:D_h-cohomology_1-0_0-1_def}\nonumber H^{p,\,q}_{D^{1,\,0}_\eta}(X,\,\C) & := & \frac{\ker\bigg(D^{1,\,0}_\eta:C^\infty_{p,\,q}(X,\,\C)\longrightarrow C^\infty_{p+1,\,q}(X,\,\C)\bigg)}{\mbox{Im}\,\bigg(D^{1,\,0}_\eta:C^\infty_{p-1,\,q}(X,\,\C)\longrightarrow C^\infty_{p,\,q}(X,\,\C)\bigg)}, \\
  H^{p,\,q}_{D^{0,\,1}_\eta}(X,\,\C) & := & \frac{\ker\bigg(D^{0,\,1}_\eta:C^\infty_{p,\,q}(X,\,\C)\longrightarrow C^\infty_{p,\,q+1}(X,\,\C)\bigg)}{\mbox{Im}\,\bigg(D^{0,\,1}_\eta:C^\infty_{p,\,q-1}(X,\,\C)\longrightarrow C^\infty_{p,\,q}(X,\,\C)\bigg)}\end{eqnarray} are isomorphic to the corresponding $\partial$-cohomology, resp. $\bar\partial$-cohomology, spaces via the following {\bf isomorphisms} induced in cohomology by $\theta_\eta$ as in Proposition \ref{Prop:cohomologies_isom}: \begin{eqnarray}\label{eqn:cohomologies_1-0_0-1_isom}\nonumber\theta_\eta:H^{p,\,q}_\partial(X,\,\C)\longrightarrow H^{p,\,q}_{D^{1,\,0}_\eta}(X,\,\C), & & \hspace{5ex} \{u\}_\partial\longmapsto\{\theta_\eta u\}_{D^{1,\,0}_\eta}, \\
 \theta_\eta:H^{p,\,q}_{\bar\partial}(X,\,\C)\longrightarrow H^{p,\,q}_{D^{0,\,1}_\eta}(X,\,\C), & & \hspace{5ex} \{u\}_{\bar\partial}\longmapsto\{\theta_\eta u\}_{D^{0,\,1}_\eta}.\end{eqnarray} 

Another immediate observation is that \begin{eqnarray}\label{eqn:eta_minus-eta_indices}D_\eta^{0,\,1} = D_{-\eta}^{0,\,1}, \hspace{3ex}\mbox{hence also}\hspace{3ex} \tau_\eta^{0,\,1} = \tau_{-\eta}^{0,\,1} \hspace{3ex}\mbox{and}\hspace{3ex} \Delta'_\eta = \Delta'_{-\eta}.\end{eqnarray}

\vspace{2ex}

On the other hand, the principal part of $\Delta''_\eta$ is the classical $\bar\partial$-Laplacian $\Delta''$, while the principal part of $\Delta'_\eta$ is $\Delta'$. This shows that $\Delta''_\eta$ and $\Delta'_\eta$ are {\bf elliptic}. Consequently, if $X$ is {\bf compact}, we get, in every bidegree $(p,\,q)$, the {\bf Hodge isomorphism} \begin{eqnarray}\label{eqn:Hodge_isom_Delta''-Delta'_eta}H^{p,\,q}_{D^{0,\,1}_\eta}(X,\,\C)\simeq{\cal H}_{\Delta''_\eta}^{p,\,q}(X,\,\C):=\ker\bigg(\Delta''_\eta:C^\infty_{p,\,q}(X,\,\C)\longrightarrow C^\infty_{p,\,q}(X,\,\C)\bigg)\end{eqnarray} mapping every $D^{0,\,1}_\eta$-cohomology class to its unique $\Delta''_\eta$-harmonic representative, as well as the analogous statement for $\Delta'_\eta$ and $\overline{D^{0,\,1}_\eta}$.

\vspace{2ex}

\vspace{2ex}

Splitting each of the $\eta$-twisted commutation relations (\ref{eqn:eta-twisted_commutation-relations}) into two identities according to the bidegrees, we see that Proposition \ref{Prop:eta-twisted_commutation-relations} translates to

\begin{Prop}\label{Prop:eta-twisted_commutation-relations_1-0_0-1} Let $(X,\,\omega)$ be a complex Hermitian manifold with $\mbox{dim}_\C X=n\geq 2$. For any $C^\infty$ function $\eta$ on $X$ such that $\eta>0$ or $\eta<0$, the following {\bf $\eta$-twisted commutation relations} hold on differential forms of any degree on $X$: \begin{eqnarray}\label{eqn:eta-twisted_commutation-relations_1-0_0-1}\nonumber & (a')& (D_\eta^{1,\,0})^\star + (\tau_\eta^{1,\,0})^\star = -i\,[\Lambda,\, \overline{D^{1,\,0}_{-\eta}}] + \frac{n}{\eta}\,[i\,\overline{d'_{-\eta}}\eta\wedge\cdot\,,\Lambda]; \\
    \nonumber  & (a'') & (D_\eta^{0,\,1})^\star + (\tau_\eta^{0,\,1})^\star = -i\,[\Lambda,\, \overline{D^{0,\,1}_{-\eta}}] + \frac{n}{\eta}\,[i\,\overline{d''_{-\eta}}\eta\wedge\cdot\,,\Lambda]; \\
    \nonumber & (b') & \overline{D^{1,\,0}_{-\eta}}^{\,\star} + \overline{\tau^{1,\,0}_{-\eta}}^{\,\star} = i\,[\Lambda,\, D^{1,\,0}_\eta] - \frac{n}{\eta}\,[i\,d'_\eta\eta\wedge\cdot\,,\Lambda];\\
  \nonumber & (b'') & \overline{D^{0,\,1}_{-\eta}}^{\,\star} + \overline{\tau^{0,\,1}_{-\eta}}^{\,\star} = i\,[\Lambda,\, D^{0,\,1}_\eta] - \frac{n}{\eta}\,[i\,d''_\eta\eta\wedge\cdot\,,\Lambda];\\
  \nonumber & (c') & D^{1,\,0}_\eta + \tau^{1,\,0}_\eta = i\,[\overline{D^{1,\,0}_{-\eta}}^{\,\star},\,\omega\wedge\cdot\,] - \frac{n}{\eta}\,[\omega\wedge\cdot\,,i\,(\overline{d'_{-\eta}}\eta\wedge\cdot\,)^\star];  \\
  \nonumber & (c'') & D^{0,\,1}_\eta + \tau^{0,\,1}_\eta = i\,[\overline{D^{0,\,1}_{-\eta}}^{\,\star},\,\omega\wedge\cdot\,] - \frac{n}{\eta}\,[\omega\wedge\cdot\,,i\,(\overline{d''_{-\eta}}\eta\wedge\cdot\,)^\star];  \\
 \nonumber & (d') & \overline{D^{1,\,0}_{-\eta}} + \overline{\tau^{1,\,0}_{-\eta}} = -i\,[(D_\eta^{1,\,0})^\star,\,\omega\wedge\cdot\,] + \frac{n}{\eta}\,[\omega\wedge\cdot\,,i\,(d'_\eta\eta\wedge\cdot\,)^\star]; \\
  & (d'') & \overline{D^{0,\,1}_{-\eta}} + \overline{\tau^{0,\,1}_{-\eta}} = -i\,[(D_\eta^{0,\,1})^\star,\,\omega\wedge\cdot\,] + \frac{n}{\eta}\,[\omega\wedge\cdot\,,i\,(d''_\eta\eta\wedge\cdot\,)^\star].\end{eqnarray}

\end{Prop}

As a consequence of these $\eta$-twisted commutation relations, we get the following analogue of Proposition \ref{Prop:eta-BKN_rough}.

\begin{Prop}\label{Prop:eta-BKN_rough_1-0_0-1} Let $(X,\,\omega)$ be a complex Hermitian manifold with $\mbox{dim}_\C X=n\geq 2$.  For any $C^\infty$ function $\eta$ on $X$ such that $\eta>0$ or $\eta<0$ and any $k\in\{0,\dots , 2n\}$, the following {\bf rough $\eta$-Bochner-Kodaira-Nakano ($\eta$-BKN) identity} holds on the $C^\infty$ forms of any degree on $X$: \begin{eqnarray}\label{eqn:eta-BKN_rough_1-0_0-1}\nonumber\Delta''_\eta = \Delta'_\eta + i\,\bigg[[D^{0,\,1}_\eta,\,\overline{D^{0,\,1}_\eta}],\,\Lambda\bigg] & + & [\overline{D^{0,\,1}_\eta},\,\overline{\tau^{0,\,1}_\eta}^{\,\star}] - \bigg[D^{0,\,1}_\eta,\,(\tau^{0,\,1}_\eta)^\star\bigg] \\
    & + & n\,\bigg[D^{0,\,1}_\eta,\,\frac{i}{\eta}\,[\partial\eta\wedge\cdot\,,\,\Lambda]\bigg] + n\,\bigg[\overline{D^{0,\,1}_\eta},\,\frac{i}{\eta}\,[\bar\partial\eta\wedge\cdot\,,\,\Lambda]\bigg].\end{eqnarray}

\end{Prop}  

\noindent {\it Proof.} Using the expression for $(D^{0,\,1}_\eta)^\star$ given in (a'') of (\ref{eqn:eta-twisted_commutation-relations_1-0_0-1}), we get the second equality below: \begin{eqnarray}\label{eqn:eta-BKN_proof_1_bis}\Delta''_\eta & = & [D_\eta^{0,\,1},\,(D^{0,\,1}_\eta)^\star] = -i\,\bigg[D^{0,\,1}_\eta,\,[\Lambda,\,\overline{D^{0,\,1}_{-\eta}}]\bigg] - [D^{0,\,1}_\eta,\,(\tau^{0,\,1}_\eta)^\star] + n\,\bigg[D_\eta^{0,\,1},\,\frac{i}{\eta}\,[\partial\eta\wedge\cdot\,,\,\Lambda]\bigg].\end{eqnarray}

Now, the Jacobi identity yields the former equality below (expressing the first term on the r.h.s. of (\ref{eqn:eta-BKN_proof_1_bis})), while the $\eta$-twisted commutation relation (b'') of (\ref{eqn:eta-twisted_commutation-relations_1-0_0-1}) yields the latter equality: \begin{eqnarray}\label{eqn:eta-BKN_proof_2_bis}\nonumber -i\,\bigg[D^{0,\,1}_\eta,\,[\Lambda,\,\overline{D^{0,\,1}_{-\eta}}]\bigg] & = & \bigg[\overline{D^{0,\,1}_{-\eta}},\,i\,[\Lambda,\,D^{0,\,1}_\eta]\bigg] + i\,\bigg[[\overline{D^{0,\,1}_{-\eta}},\,D^{0,\,1}_\eta],\,\Lambda\bigg] \\
  & = & i\,\bigg[[D^{0,\,1}_\eta,\,\overline{D^{0,\,1}_{-\eta}}],\,\Lambda\bigg] + \bigg[\overline{D^{0,\,1}_{-\eta}},\,\overline{D^{0,\,1}_{-\eta}}^\star + \overline{\tau^{0,\,1}_{-\eta}}^\star + \frac{n}{\eta}\,[i\bar\partial\eta\wedge\cdot,\,\Lambda]\bigg].\end{eqnarray}

Plugging into (\ref{eqn:eta-BKN_proof_1_bis}) the expression given for $-i\,\bigg[D^{0,\,1}_\eta,\,[\Lambda,\,\overline{D^{0,\,1}_{-\eta}}]\bigg]$ in (\ref{eqn:eta-BKN_proof_2_bis}) and using the equality $[\overline{D^{0,\,1}_{-\eta}},\,\overline{D^{0,\,1}_{-\eta}}^{\,\star}] = \Delta'_{-\eta}$ as well as (\ref{eqn:eta_minus-eta_indices}), we get (\ref{eqn:eta-BKN_rough_1-0_0-1}). \hfill $\Box$

\vspace{3ex}

The first step towards refining the above $\eta$-BKN identity by incorporating some of the first-order terms on the right into a twisted Laplacian is the following analogue of Lemma \ref{Lem:preliminary_BKN-refined}.

\begin{Lem}\label{Lem:preliminary_BKN-refined_1-0_0-1} Let $(X,\,\omega)$ be a complex Hermitian manifold with $\mbox{dim}_\C X=n\geq 2$.  For any $C^\infty$ function $\eta$ on $X$ such that $\eta>0$ or $\eta<0$ and any bidegree $(p,\,q)$, the following identities hold: \\

 \hspace{1ex} (i)\hspace{2ex} $[L,\,\tau^{0,\,1}_\eta] = 3\,D^{0,\,1}_\eta\omega\wedge\cdot\,$,  \hspace{6ex}   (ii)\, $[\Lambda,\,\tau^{0,\,1}_\eta] = 2i\,\overline{\tau^{0,\,1}_{-\eta}}^{\,\star}$,

  \vspace{2ex}
  
\hspace{1ex}  (iii)\hspace{2ex}$[D^{0,\,1}_\eta,\,\overline{D^{0,\,1}_{-\eta}}^{\,\star}] = -[D^{0,\,1}_\eta,\,\overline{\tau^{0,\,1}_{-\eta}}^{\,\star}] - \bigg[D^{0,\,1}_\eta,\,\frac{n}{\eta}\,(\partial\eta\wedge\cdot\,)^\star\bigg],$ 

\begin{eqnarray*}(iv)\hspace{2ex} [D^{0,\,1}_\eta,\,(D^{0,\,1}_\eta)^\star] + [D^{0,\,1}_\eta,\,(\tau^{0,\,1}_\eta)^\star] - \bigg[\overline{D^{0,\,1}_{-\eta}},\,\overline{\tau^{0,\,1}_{-\eta}}^{\,\star}\bigg] &  = & [D^{0,\,1}_\eta + \tau^{0,\,1}_\eta,\,(D^{0,\,1}_\eta)^\star + (\tau^{0,\,1}_\eta)^\star] + S_\omega^{(\eta)''} \\
  & + & \frac{n}{\eta}\,[\tau^{0,\,1}_\eta,\,(\bar\partial\eta\wedge\cdot)^\star],\end{eqnarray*} where

\begin{eqnarray*}S_\omega^{(\eta)''} : = \frac{i}{2}\,[\Lambda,\,[\Lambda,\,\overline{D^{0,\,1}_\eta}D^{0,\,1}_\eta\omega\wedge\cdot\,]] - [D^{0,\,1}_\eta\omega\wedge\cdot\,,\, (D^{0,\,1}_\eta\omega\wedge\cdot\,)^\star].\end{eqnarray*}

\end{Lem}

\noindent {\it Proof.} The computations are the exact analogues of those forming the proof of Lemma \ref{Lem:preliminary_BKN-refined} when $D_\eta$ is replaced by $D^{0,\,1}_\eta$, $\tau_\eta$ is replaced by $\tau^{0,\,1}_\eta$, etc. The details will not be repeated.  \hfill $\Box$

\vspace{3ex}

The main result of this section is the following analogue of Theorem \ref{The:eta-BKN}.

\begin{The}\label{The:eta-BKN_1-0_0-1} Let $(X,\,\omega)$ be a complex Hermitian manifold with $\mbox{dim}_\C X=n\geq 2$.  For any $C^\infty$ function $\eta$ on $X$ such that $\eta>0$ or $\eta<0$ and any $k\in\{0,\dots , 2n\}$, the following {\bf refined $\eta$-Bochner-Kodaira-Nakano ($\eta$-BKN) identity} holds on the $C^\infty$ forms on $X$: \begin{eqnarray}\label{eqn:eta-BKN_1-0_0-1}\Delta''_\eta = \bigg[\overline{D^{0,\,1}_\eta} + \overline{\tau^{0,\,1}_\eta},\,\overline{D^{0,\,1}_\eta}^{\,\star} + \overline{\tau^{0,\,1}_\eta}^{\,\star}\bigg] & + & T_\omega^{(\eta)''} + i\,\bigg[[D^{0,\,1}_\eta,\,\overline{D^{0,\,1}_\eta}],\,\Lambda\bigg] \\
 \nonumber   & + & n\,\bigg[\overline{D^{0,\,1}_\eta} + \overline{\tau^{0,\,1}_\eta},\,\frac{1}{\eta}\,(\partial\eta\wedge\cdot\,)^\star\bigg] -n\,\bigg[D^{0,\,1}_\eta,\,\frac{1}{\eta}\,(\bar\partial\eta\wedge\cdot\,)^\star\bigg],\end{eqnarray} where $T_\omega^{(\eta)''}$ is the zero-th order operator defined by \begin{eqnarray*}T_\omega^{(\eta)''}:= \overline{S_\omega^{(-\eta)''}} = -\frac{i}{2}\,[\Lambda,\,[\Lambda,\,D^{0,\,1}_\eta\overline{D^{0,\,1}_\eta}\omega\wedge\cdot\,]] - [\overline{D^{0,\,1}_\eta}\omega\wedge\cdot,\, (\overline{D^{0,\,1}_\eta}\omega\wedge\cdot\,)^\star].\end{eqnarray*}

\end{The}  

\noindent {\it Proof.} It is the exact analogue of the proof of Theorem \ref{The:eta-BKN}, based on the use of Proposition \ref{Prop:eta-BKN_rough_1-0_0-1} instead of Proposition \ref{Prop:eta-BKN_rough} and the use of Lemma \ref{Lem:preliminary_BKN-refined_1-0_0-1} instead of Lemma \ref{Lem:preliminary_BKN-refined}. We have also replaced $-\eta$ by $\eta$ in the subscripts of the statement thanks to the already noticed equalities $D_\eta^{0,\,1} = D_{-\eta}^{0,\,1}$ and $\tau_\eta^{0,\,1} = \tau_{-\eta}^{0,\,1}$.  \hfill $\Box$

\vspace{2ex}

We now compute the operator $i\,\bigg[[D^{0,\,1}_\eta,\,\overline{D^{0,\,1}_\eta}],\,\Lambda\bigg]$ that plays in (\ref{eqn:eta-BKN_1-0_0-1}) the role of the curvature operator of the classical Bochner-Kodaira-Nakano identity and the $(1,\,1)$-form $D^{0,\,1}_\eta\overline{D^{0,\,1}_\eta}\omega$ that features in the definition of $T_\omega^{(\eta)''}$. The following statement is the analogue of Proposition \ref{Prop:curvature_computation}.

\begin{Prop}\label{Prop:curvature_computation_1-0_0-1} Let $(X,\,\omega)$ be a complex Hermitian manifold with $\mbox{dim}_\C X=n\geq 2$. Fix any $C^\infty$ function $\eta$ on $X$ such that $\eta>0$ or $\eta<0$.

\vspace{1ex}

  (i)\, For any bidegree $(p,\,q)$, the following identities hold: \begin{eqnarray}\label{eqn:curvature_computation_1-0_0-1_first}\nonumber i\,[D^{0,\,1}_\eta,\,\overline{D^{0,\,1}_\eta}]\,(u^{p,\,q}) & = & - \frac{1}{\eta}\,\bigg(i\bar\partial\eta\wedge\partial u^{p,\,q} +i\partial\eta\wedge\bar\partial u^{p,\,q}\bigg) + \frac{p-q}{\eta}\,\bigg(\frac{2}{\eta}\,i\partial\eta\wedge\bar\partial\eta - i\partial\bar\partial\eta\bigg)\wedge u^{p,\,q},\end{eqnarray} and \begin{eqnarray}\label{eqn:curvature_computation_1-0_0-1_second}\nonumber i\,\bigg[[D^{0,\,1}_\eta,\,\overline{D^{0,\,1}_\eta}],\,\Lambda\bigg]\,(u^{p,\,q}) & = & \frac{1}{\eta}\,[\Lambda,\,i\bar\partial\eta\wedge\cdot\,]\,\partial u^{p,\,q} + \frac{1}{\eta}\,[\Lambda,\,i\partial\eta\wedge\cdot\,]\,\bar\partial u^{p,\,q} \\
    \nonumber & + & \frac{1}{\eta}\,i\bar\partial\eta\wedge[\Lambda,\,\partial]\,u^{p,\,q} + \frac{1}{\eta}\,i\partial\eta\wedge[\Lambda,\,\bar\partial]\,u^{p,\,q}  \\
    & + & (p-q)\,\bigg[\bigg(\frac{2}{\eta^2}\,i\partial\eta\wedge\bar\partial\eta - \frac{1}{\eta}\,i\partial\bar\partial\eta\bigg)\wedge\cdot\,,\,\Lambda\bigg]\,u^{p,\,q}\end{eqnarray} for any $C^\infty$ $(p,\,q)$-form $u^{p,\,q}$ on $X$.

\vspace{1ex}

(ii)\, The following identities hold: \begin{eqnarray}\label{eqn:eta-hessian_omega}\nonumber D^{0,\,1}_\eta\overline{D^{0,\,1}_\eta}\omega & = & \eta^2\,\bar\partial\bigg(\frac{1}{\eta}\,\partial\bigg(\frac{1}{\eta}\,\omega\bigg)\bigg) \\
  & = & -\partial\bar\partial\omega - \frac{1}{\eta}\,\bigg(2\,\bar\partial\eta\wedge\partial\omega - \partial\eta\wedge\bar\partial\omega\bigg) + \frac{1}{\eta}\,\bigg(\partial\bar\partial\eta - \frac{3}{\eta}\,\partial\eta\wedge\bar\partial\eta\bigg)\wedge\omega.\end{eqnarray}

\end{Prop}  

\noindent {\it Proof.} (i)\, We compute separately the two terms in the sum $[D^{0,\,1}_\eta,\,\overline{D^{0,\,1}_\eta}]\,(u^{p,\,q}) = (D^{0,\,1}_\eta\overline{D^{0,\,1}_\eta})\,(u^{p,\,q}) + (\overline{D^{0,\,1}_\eta}D^{0,\,1}_\eta)\,(u^{p,\,q})$.

\vspace{1ex}

$\bullet$ We use (\ref{eqn:conjugate_D-01_eta}) to get the first line and then (\ref{eqn:D_eta_1-0-1_def}) to get the next line below: \begin{eqnarray*}(D^{0,\,1}_\eta\overline{D^{0,\,1}_\eta})\,(u^{p,\,q}) & = & D^{0,\,1}_\eta\bigg(\partial u^{p,\,q} - \frac{q}{\eta}\,\partial\eta\wedge u^{p,\,q}\bigg) \\
  & = & \bar\partial\bigg(\partial u^{p,\,q} - \frac{q}{\eta}\,\partial\eta\wedge u^{p,\,q}\bigg) - \frac{p+1}{\eta}\,\bar\partial\eta\wedge\bigg(\partial u^{p,\,q} - \frac{q}{\eta}\,\partial\eta\wedge u^{p,\,q}\bigg).\end{eqnarray*}

This leads to \begin{eqnarray}\label{eqn:curvature_computation_1-0_0-1_first_proof_1}\nonumber(D^{0,\,1}_\eta\overline{D^{0,\,1}_\eta})\,(u^{p,\,q}) = -\partial\bar\partial u^{p,\,q} & + & \frac{1}{\eta}\,\bigg(q\,\partial\eta\wedge\bar\partial u^{p,\,q} - (p+1)\,\bar\partial\eta\wedge\partial u^{p,\,q}\bigg) \\
  & + & \frac{q}{\eta}\,\bigg(\partial\bar\partial\eta - \frac{p+2}{\eta}\,\partial\eta\wedge\bar\partial\eta\bigg)\wedge u^{p,\,q}.\end{eqnarray}

$\bullet$ Meanwhile, we use (\ref{eqn:D_eta_1-0-1_def}) to get the first line and then (\ref{eqn:conjugate_D-01_eta}) to get the next line below: \begin{eqnarray*}(\overline{D^{0,\,1}_\eta}D^{0,\,1}_\eta)\,(u^{p,\,q}) & = & \overline{D^{0,\,1}_\eta}\bigg(\bar\partial u^{p,\,q} - \frac{p}{\eta}\,\bar\partial\eta\wedge u^{p,\,q}\bigg) \\
  & = & \partial\bigg(\bar\partial u^{p,\,q} - \frac{p}{\eta}\,\bar\partial\eta\wedge u^{p,\,q}\bigg) - \frac{q+1}{\eta}\,\partial\eta\wedge\bigg(\bar\partial u^{p,\,q} - \frac{p}{\eta}\,\bar\partial\eta\wedge u^{p,\,q}\bigg).\end{eqnarray*}

This leads to \begin{eqnarray}\label{eqn:curvature_computation_1-0_0-1_first_proof_2}\nonumber(\overline{D^{0,\,1}_\eta}D^{0,\,1}_\eta)\,(u^{p,\,q}) = \partial\bar\partial u^{p,\,q} & + & \frac{1}{\eta}\,\bigg(p\,\bar\partial\eta\wedge\partial u^{p,\,q} - (q+1)\,\partial\eta\wedge\bar\partial u^{p,\,q}\bigg) \\
  & + & \frac{p}{\eta}\,\bigg(\frac{q+2}{\eta}\,\partial\eta\wedge\bar\partial\eta - \partial\bar\partial\eta\bigg)\wedge u^{p,\,q}.\end{eqnarray}

\vspace{2ex}

Adding up (\ref{eqn:curvature_computation_1-0_0-1_first_proof_1}) and (\ref{eqn:curvature_computation_1-0_0-1_first_proof_2}) yields the first equality claimed under (i).

\vspace{1ex}

Now, starting from $i\,\bigg[[D^{0,\,1}_\eta,\,\overline{D^{0,\,1}_\eta}],\,\Lambda\bigg]\,(u^{p,\,q}) = i\,[D^{0,\,1}_\eta,\,\overline{D^{0,\,1}_\eta}]\,(\Lambda u^{p,\,q}) -\Lambda\bigg(i\,[D^{0,\,1}_\eta,\,\overline{D^{0,\,1}_\eta}]\,(u^{p,\,q})\bigg)$ and using the first equality under (i), we get \begin{eqnarray}\label{eqn:curvature_computation_1-0_0-1_second_proof_1}\nonumber i\,\bigg[[D^{0,\,1}_\eta,\,\overline{D^{0,\,1}_\eta}],\,\Lambda\bigg]\,(u^{p,\,q}) & = & \bigg(-\frac{1}{\eta}\,i\bar\partial\eta\wedge\partial\Lambda u^{p,\,q} + \frac{1}{\eta}\,\Lambda(i\bar\partial\eta\wedge\partial u^{p,\,q})\bigg) \\
  \nonumber  & + & \bigg(-\frac{1}{\eta}\,i\partial\eta\wedge\bar\partial\Lambda u^{p,\,q} + \frac{1}{\eta}\,\Lambda(i\partial\eta\wedge\bar\partial u^{p,\,q})\bigg) \\
  & + & \frac{p-q}{\eta}\,\bigg[\bigg(\frac{2}{\eta}\,i\partial\eta\wedge\bar\partial\eta - i\partial\bar\partial\eta\bigg)\wedge\cdot\,,\Lambda\bigg]\,u^{p,\,q}.\end{eqnarray}

It remains to notice that the first parenthesis on the right of (\ref{eqn:curvature_computation_1-0_0-1_second_proof_1}) equals \begin{eqnarray*}\frac{1}{\eta}\,[\Lambda,\,i\bar\partial\eta\wedge\cdot\,]\,\partial u^{p,\,q} + \frac{1}{\eta}\,i\bar\partial\eta\wedge[\Lambda,\,\partial]\,u^{p,\,q},\end{eqnarray*} while the second parenthesis on the right on the right of (\ref{eqn:curvature_computation_1-0_0-1_second_proof_1}) equals \begin{eqnarray*}\frac{1}{\eta}\,[\Lambda,\,i\partial\eta\wedge\cdot\,]\,\bar\partial u^{p,\,q} + \frac{1}{\eta}\,i\partial\eta\wedge[\Lambda,\,\bar\partial]\,u^{p,\,q},\end{eqnarray*} to get (\ref{eqn:curvature_computation_1-0_0-1_second}).

\vspace{1ex}

(ii)\, Using (\ref{eqn:D_eta_1-0_0-1}), we get: \begin{eqnarray*}D^{0,\,1}_\eta\overline{D^{0,\,1}_\eta}\omega =\theta_\eta\bar\partial\theta_\eta^{-1}\overline{\theta_\eta}\partial\bigg(\frac{1}{\eta}\omega\bigg) = (\theta_\eta\bar\partial\theta_\eta^{-1})\bigg(\eta\,\partial\bigg(\frac{1}{\eta}\omega\bigg)\bigg) = (\theta_\eta\bar\partial)\bigg(\frac{1}{\eta}\,\partial\bigg(\frac{1}{\eta}\omega\bigg)\bigg) = \eta^2\,\bar\partial\bigg(\frac{1}{\eta}\,\partial\bigg(\frac{1}{\eta}\omega\bigg)\bigg),\end{eqnarray*} which is the former identity in (\ref{eqn:eta-hessian_omega}).

The latter identity in (\ref{eqn:eta-hessian_omega}) can be proved by a direct computation. We get $\overline{D^{0,\,1}_\eta}\omega = \partial\omega - \frac{1}{\eta}\,\partial\eta\wedge\omega$, hence \begin{eqnarray*}D^{0,\,1}_\eta\overline{D^{0,\,1}_\eta}\omega & = & \bar\partial\bigg(\partial\omega - \frac{1}{\eta}\,\partial\eta\wedge\omega\bigg) - \frac{2}{\eta}\,\bar\partial\eta\wedge\bigg(\partial\omega - \frac{1}{\eta}\,\partial\eta\wedge\omega\bigg) \\
  & = & -\partial\bar\partial\omega + \frac{1}{\eta^2}\,\bar\partial\eta\wedge\partial\eta\wedge\omega - \frac{1}{\eta}\,\bar\partial\partial\eta\wedge\omega + \frac{1}{\eta}\,\partial\eta\wedge\bar\partial\omega - \frac{2}{\eta}\,\bar\partial\eta\wedge\partial\omega + \frac{2}{\eta^2}\,\bar\partial\eta\wedge\partial\eta\wedge\omega.\end{eqnarray*} Collecting terms, we get the latter identity in (\ref{eqn:eta-hessian_omega}).    \hfill $\Box$

\vspace{2ex}

If $\varphi$ is the real-valued $C^\infty$ function on $X$ defined by $\eta=e^{-\varphi}$, the quantities depending on $\eta$ featuring in (\ref{eqn:curvature_computation_1-0_0-1_second}) read: \begin{eqnarray}\label{eqn:eta-derivatives_phi}\nonumber \frac{1}{\eta}\,i\partial\eta = -i\partial\varphi; \hspace{3ex} \frac{1}{\eta}\,i\bar\partial\eta & = & -i\bar\partial\varphi; \hspace{3ex} \frac{1}{\eta}\,i\partial\bar\partial\eta = i\partial\varphi\wedge\bar\partial\varphi - i\partial\bar\partial\varphi; \\
  \frac{2}{\eta^2}\,i\partial\eta\wedge\bar\partial\eta - \frac{1}{\eta}\,i\partial\bar\partial\eta & = & i\partial\varphi\wedge\bar\partial\varphi + i\partial\bar\partial\varphi.\end{eqnarray}

\section{Vanishing of certain $L^2$ $\Delta''_\eta$-harmonic spaces on certain non-compact complete complex manifolds}\label{section:vanishing_harmonic} Let $X$ be a complex manifold with $\mbox{dim}_\C X=n$. With every $C^\infty$ function $\eta:X\longrightarrow(0,\,\infty)$, we associate the $C^\infty$ real $(1,\,1)$-form \begin{eqnarray*}\gamma_\eta:=\frac{2}{\eta^2}\,i\partial\eta\wedge\bar\partial\eta - \frac{1}{\eta}\,i\partial\bar\partial\eta\end{eqnarray*} featuring in the ``curvature'' term (\ref{eqn:curvature_computation_1-0_0-1_second}) of the refined $\eta$-BKN identity (\ref{eqn:eta-BKN_1-0_0-1}). Note that the former term in $\gamma_\eta$ is $\geq 0$ on $X$, while the latter term may be signless in general and, by the maximum principle, is signless or vanishes identically if $X$ is compact.

\begin{Lem}\label{Lem:gamma_eta_eta-closed}$\overline{D_\eta^{0,\,1}}\gamma_\eta = 0$. Hence, $d(\frac{1}{\eta}\,\gamma_\eta) = 0$. 

\end{Lem}    

\noindent {\it Proof.} Straightforward calculations yield: \begin{eqnarray*}\partial\gamma_\eta = -\frac{1}{\eta^2}\,\partial\eta\wedge i\partial\bar\partial\eta \hspace{5ex}\mbox{and}\hspace{5ex} \frac{1}{\eta}\,\partial\eta\wedge\gamma_\eta = - \frac{1}{\eta^2}\,\partial\eta\wedge i\partial\bar\partial\eta.\end{eqnarray*} Hence, using (\ref{eqn:conjugate_D-01_eta}) for the former equality below, we get: \begin{eqnarray*}\overline{D_\eta^{0,\,1}}\gamma_\eta = \partial\gamma_\eta - \frac{1}{\eta}\,\partial\eta\wedge\gamma_\eta = 0.\end{eqnarray*}

This proves the former claimed equality, which in turn implies the latter one thanks to (\ref{eqn:D_eta_1-0_0-1}) .  \hfill $\Box$

\vspace{2ex}

Henceforth, we will make the assumption that $\gamma_\eta$ is {\it positive definite} at every point of $X$: \begin{eqnarray}\label{eqn:assumption_H}\gamma_\eta:=\frac{2}{\eta^2}\,i\partial\eta\wedge\bar\partial\eta - \frac{1}{\eta}\,i\partial\bar\partial\eta > 0.\end{eqnarray} This assumption is made necessary by the need to create positivity in the refined $\eta$-BKN identity (\ref{eqn:eta-BKN_1-0_0-1}), specifically in its ``curvature'' term (\ref{eqn:curvature_computation_1-0_0-1_second}).

Assumption (\ref{eqn:assumption_H}) means that $\gamma_\eta$ defines a Hermitian metric on $X$. Together with Lemma \ref{Lem:gamma_eta_eta-closed}, this further implies that $\frac{1}{\eta}\,\gamma_\eta$ defines a {\it K\"ahler metric} on $X$. Meanwhile, assumption (\ref{eqn:assumption_H}) is never satisfied on a {\it compact} manifold $X$ since, otherwise, the maximum principle would imply that $\eta$ is constant, hence $\gamma_\eta$ would vanish identically, contradicting the strict positivity assumption (\ref{eqn:assumption_H}).

Summing up, assumption (\ref{eqn:assumption_H}) implies that $X$ is a {\it non-compact K\"ahler} complex manifold. We now choose the Hermitian metric $\omega$ on $X$ to be $\gamma_\eta$, i.e. $\omega := \gamma_\eta$. It is with respect to this metric that all the norms, inner products, formal adjoints and other objects will be considered in this $\S$\ref{section:vanishing_harmonic}. Note that  $\omega = \gamma_\eta$ is not K\"ahler, only $\frac{1}{\eta}\,\gamma_\eta$ is.

\vspace{2ex}

We can now state the main result of this section, an application of Theorem \ref{The:eta-BKN_1-0_0-1}.

\begin{The}\label{The:vanishing_harmonic} Let $X$ be a (non-compact) complex manifold with $\mbox{dim}_\C X=n$. Suppose there exists a $C^\infty$ function $\eta:X\longrightarrow(0,\,\infty)$ satisfying the following three conditions:

\vspace{1ex}

(i)\, the $C^\infty$ $(1,\,1)$-form $\gamma_\eta:=\frac{2}{\eta^2}\,i\partial\eta\wedge\bar\partial\eta - \frac{1}{\eta}\,i\partial\bar\partial\eta$ is {\bf positive definite} at every point of $X$;

\vspace{1ex}

(ii)\, the Hermitian metric $\gamma_\eta$ defined on $X$ under (i) is {\bf complete};

\vspace{1ex}

(iii)\, the pointwise $\gamma_\eta$-norm $|\partial\eta| = |\partial\eta|_{\gamma_\eta}$ of the $(1,\,0)$-form $\partial\eta$ is {\bf small} relative to $\eta$ in that \begin{eqnarray}\label{eqn:vanishing_harmonic_H1}C_1(\eta):=\sup\limits_X\frac{|\partial\eta|}{\eta} < \frac{1}{10n + 4n\,\sqrt{n} + 8n\,C(\varphi)},\end{eqnarray} where $\varphi:=-\log\eta$ and $C(\varphi):=\sup\limits_X|i\partial\bar\partial\varphi|_{\gamma_\eta}$.



\vspace{2ex}

Then, for any bidegree $(p,\,q)$ such that either $\bigg(p>q \hspace{1ex}\mbox{and}\hspace{1ex} p+q\geq n+1\bigg)$ or $\bigg(p<q \hspace{1ex}\mbox{and}\hspace{1ex} p+q\leq n-1\bigg)$, the space of $\Delta''_\eta$-harmonic $L^2_{\gamma_\eta}$-forms on $X$ of bidegree $(p,\,q)$ {\bf vanishes}: \begin{eqnarray}\label{eqn:vanishing_harmonic_conclusion}{\cal H}^{p,\,q}_{\Delta''_\eta}(X,\,\C):=\ker\bigg(\Delta''_\eta:\mbox{Dom}_{p,\,q}(\Delta''_\eta)\longrightarrow L^2_{p,\,q}(X,\,\C)\bigg) = \{0\},\end{eqnarray} where $\Delta''_\eta$ is the closed and densely defined unbounded extension to the space of $L^2_{\gamma_\eta}$-forms of bidegree $(p,\,q)$ of the operator $\Delta''_\eta$ previously defined on $C^\infty$ forms w.r.t. the metric $\omega=\gamma_\eta$. 

\end{The}

\vspace{2ex}

Note that the metric $\gamma_\eta$ and the positive real $C_1(\eta)$ are invariant under rescalings of $\eta$ by positive constants $\lambda$: \begin{eqnarray*}\gamma_{\lambda\eta} = \gamma_\eta; \hspace{5ex} C_1(\lambda\eta) = C_1(\eta).\end{eqnarray*}

\vspace{2ex}

\noindent {\it Proof of Theorem \ref{The:vanishing_harmonic}.} With the choice of metric $\omega=\gamma_\eta$, Lemma \ref{Lem:gamma_eta_eta-closed} implies that $D_\eta^{0,\,1}\omega = 0$, hence also $D_\eta^{0,\,1}\overline{D_\eta^{0,\,1}}\omega = 0$. Consequently, we get: \begin{eqnarray*}\tau_\eta^{0,\,1} = [\Lambda,\,D_\eta^{0,\,1}\omega\wedge\cdot\,] = 0 \hspace{3ex}\mbox{and}\hspace{3ex} T_\omega^{(\eta)''}= -\frac{i}{2}\,[\Lambda,\,[\Lambda,\,D^{0,\,1}_\eta\overline{D^{0,\,1}_\eta}\omega\wedge\cdot\,]] - [\overline{D^{0,\,1}_\eta}\omega\wedge\cdot,\, (\overline{D^{0,\,1}_\eta}\omega\wedge\cdot\,)^\star] = 0,\end{eqnarray*} as well as $0 = [\Lambda,\,\overline{D_\eta^{0,\,1}}\omega\wedge\cdot\,] = [\Lambda,\,\partial\omega\wedge\cdot\,] - [\Lambda,\, \frac{1}{\eta}\,\partial\eta\wedge\omega\wedge\cdot\,]$. This last equality translates to \begin{eqnarray}\label{eqn:tau_omega-eta}\tau = [\Lambda,\, \frac{1}{\eta}\,\partial\eta\wedge\omega\wedge\cdot\,].\end{eqnarray}

These equalities reduce (\ref{eqn:formulae_D_eta_1-0_0-1_tau_eta_1-0_0-1_no2}) and (\ref{eqn:conjugate-star_D-01_eta_conjugate-star_tau-01_eta}) to the following formulae holding in any bidegree $(p,\,q)$: \begin{eqnarray}\label{eqn:formulae_D_eta_1-0_0-1_tau_eta_1-0_0-1_no2_bis}\nonumber D_\eta^{0,\,1} = \bar\partial - \frac{p}{\eta}\,\bar\partial\eta\wedge\cdot\,  \hspace{5ex} & \mbox{and} & \hspace{5ex} \overline{D_\eta^{0,\,1}} = \partial - \frac{q}{\eta}\,\partial\eta\wedge\cdot\, \\
 \overline{D^{0,\,1}_\eta}^{\,\star} = \partial^\star - \frac{q}{\eta}\,(\partial\eta\wedge\cdot\,)^\star \hspace{5ex} & \mbox{and} & \hspace{5ex} ({D^{0,\,1}_\eta})^{\,\star} =  \bar\partial^\star - \frac{p}{\eta}\,(\bar\partial\eta\wedge\cdot\,)^\star.\end{eqnarray}

Similarly, our refined $\eta$-BKN identity (\ref{eqn:eta-BKN_1-0_0-1}) reduces to \begin{eqnarray}\label{eqn:eta-BKN_1-0_0-1_simplified}\Delta''_\eta =  \Delta'_\eta + i\,\bigg[[D^{0,\,1}_\eta,\,\overline{D^{0,\,1}_\eta}],\,\Lambda\bigg] + n\,\bigg[\overline{D^{0,\,1}_\eta} + \overline{\tau^{0,\,1}_\eta},\,\frac{1}{\eta}\,(\partial\eta\wedge\cdot\,)^\star\bigg] -n\,\bigg[D^{0,\,1}_\eta,\,\frac{1}{\eta}\,(\bar\partial\eta\wedge\cdot\,)^\star\bigg],\end{eqnarray} where, in bidegree $(p,\,q)$, the curvature term is given by \begin{eqnarray}\label{eqn:eta-BKN_1-0_0-1_curvature_1} i\,\bigg[[D^{0,\,1}_\eta,\,\overline{D^{0,\,1}_\eta}],\,\Lambda\bigg] = A + (p-q)\,[\gamma_\eta\wedge\cdot\,,\,\Lambda] = A + (p-q)(p+q-n)\,\mbox{Id} ,\end{eqnarray} while the first-order operator $A$ is given by \begin{eqnarray}\label{eqn:eta-BKN_1-0_0-1_curvature_2}\nonumber A & = & \frac{1}{\eta}\,[\Lambda,\,i\bar\partial\eta\wedge\cdot\,]\,\partial + \frac{1}{\eta}\,[\Lambda,\,i\partial\eta\wedge\cdot\,]\,\bar\partial  + \frac{1}{\eta}\,i\bar\partial\eta\wedge[\Lambda,\,\partial] + \frac{1}{\eta}\,i\partial\eta\wedge[\Lambda,\,\bar\partial] \\
  \nonumber & = & \frac{1}{\eta}\,[\Lambda,\,i\bar\partial\eta\wedge\cdot\,]\,\partial + \frac{1}{\eta}\,[\Lambda,\,i\partial\eta\wedge\cdot\,]\,\bar\partial  + \frac{1}{\eta}\,i\bar\partial\eta\wedge i(\bar\partial^\star + \bar\tau^\star) - \frac{1}{\eta}\,i\partial\eta\wedge i(\partial^\star + \tau^\star) \\
  \nonumber & = & \frac{1}{\eta}\,[\Lambda,\,i\bar\partial\eta\wedge\cdot\,]\,\bigg(\overline{D_\eta^{0,\,1}} + \frac{q}{\eta}\,\partial\eta\wedge\cdot\,\bigg) + \frac{1}{\eta}\,[\Lambda,\,i\partial\eta\wedge\cdot\,]\,\bigg(D_\eta^{0,\,1} + \frac{p}{\eta}\,\bar\partial\eta\wedge\cdot\,\bigg) \\
  & - & \frac{1}{\eta}\,\bar\partial\eta\wedge\bigg(({D^{0,\,1}_\eta})^{\,\star} + \frac{p}{\eta}\,(\bar\partial\eta\wedge\cdot\,)^\star + \bar\tau^\star\bigg) + \frac{1}{\eta}\,\partial\eta\wedge\bigg(\overline{D^{0,\,1}_\eta}^{\,\star} + \frac{q}{\eta}\,(\partial\eta\wedge\cdot\,)^\star + \tau^\star\bigg),\end{eqnarray} where we have used the standard Hermitian commutation relations (\ref{eqn:standard-comm-rel}) to get the second equality and (\ref{eqn:formulae_D_eta_1-0_0-1_tau_eta_1-0_0-1_no2_bis}) to get the third one.

Now, the pointwise operator norm of $\Lambda = \Lambda_\omega$ induced by the pointwise norm $|\,\,| = |\,\,|_\omega$ defined by any metric $\omega$ on $X$ satisfies \begin{eqnarray*}|\Lambda| = |L| = \sup_{|u|=1}|\omega\wedge u|\leq\sqrt{n}\end{eqnarray*} since $|\omega|_\omega = \sqrt{n}$ at every point of $X$, where $L = L_\omega = \omega\wedge\cdot$.

Consequently, the pointwise operator norm of $\tau = \tau_\omega$ induced by the pointwise norm $|\,\,| = |\,\,|_\omega$ satisfies \begin{eqnarray*}|\tau| =  \sup_{|u|=1}|\tau u| \leq \sup_{|u|=1}|\Lambda(\partial\omega\wedge u)| + \sup_{|u|=1}|\partial\omega\wedge\Lambda u| \leq 2\sqrt{n}\,|\partial\omega|.\end{eqnarray*} Since in our case $\omega = \gamma_\eta$, we have $\partial\omega = \partial\gamma_\eta = -(2/\eta^2)\,\partial\eta\wedge i\partial\bar\partial\eta$. Hence, if we put $C_2(\eta):=\sup\limits_X\frac{|\partial\bar\partial\eta|}{\eta}$, we get: \begin{eqnarray*}|\tau| \leq 4\sqrt{n}\,\frac{|\partial\eta|}{\eta}\,\frac{|\partial\bar\partial\eta|}{\eta} \leq 4\sqrt{n}\,C_1(\eta)\,C_2(\eta)\end{eqnarray*} at every point of $X$.

We will now estimate each of the last three terms on the r.h.s. of (\ref{eqn:eta-BKN_1-0_0-1_simplified}). 

\vspace{1ex}

$\bullet$ {\it Estimating $A$.}

\vspace{1ex}

From (\ref{eqn:eta-BKN_1-0_0-1_curvature_2}) combined with these remarks, we deduce that, for every $(p,\,q)$-form $u$, we have: \begin{eqnarray*}|\langle\langle Au,\,u\rangle\rangle| \leq \sqrt{n}\,C_1(\eta)\,||u||\, & \bigg(& ||\overline{D_\eta^{0,\,1}}u|| + q\,C_1(\eta)\,||u|| + ||D_\eta^{0,\,1}u|| + p\,C_1(\eta)\,||u|| \\
  & & + ||({D^{0,\,1}_\eta})^{\,\star} u|| + p\,C_1(\eta)\,||u|| + 4\sqrt{n}\,C_1(\eta)\,C_2(\eta)\,||u|| \\
  & & + ||\overline{D^{0,\,1}_\eta}^{\,\star}u|| + q\,C_1(\eta)\,||u|| + 4\sqrt{n}\,C_1(\eta)\,C_2(\eta)\,||u||\bigg).\end{eqnarray*} This means that \begin{eqnarray*}|\langle\langle Au,\,u\rangle\rangle| \leq \sqrt{n}\,C_1(\eta)\,\bigg(||\overline{D_\eta^{0,\,1}}u||\,||u|| & + & ||D_\eta^{0,\,1}u||\,||u|| + ||({D^{0,\,1}_\eta})^{\,\star} u||\,||u|| + ||\overline{D^{0,\,1}_\eta}^{\,\star}u||\,||u|| \\
  & + & 2(p+q)\,C_1(\eta)\,||u||^2 +8\sqrt{n}\,C_1(\eta)\,C_2(\eta)\,||u||^2\bigg).\end{eqnarray*}

Applying the elementary inequality $ab\leq (a^2 + b^2)/2$ with $b=|u|$ and $a$ each of the norms of $\overline{D^{0,\,1}_\eta}^{\,\star}u$, $D_\eta^{0,\,1}u$, $({D^{0,\,1}_\eta})^{\,\star} u$ and $\overline{D^{0,\,1}_\eta}^{\,\star}u$, we further get the inequality: \begin{eqnarray*}|\langle\langle Au,\,u\rangle\rangle| \leq \frac{\sqrt{n}}{2}\,C_1(\eta)\,\bigg(||\overline{D_\eta^{0,\,1}}u||^2 + ||\overline{D^{0,\,1}_\eta}^{\,\star}u||^2 + ||D_\eta^{0,\,1}u||^2 + ||({D^{0,\,1}_\eta})^{\,\star} u||^2 + C_{1,\,2}(p,\,q,\,n)\,||u||^2\bigg)\end{eqnarray*} for every $(p,\,q)$-form $u$ on $X$, where we put \begin{eqnarray*}C_{1,\,2}(p,\,q,\,n):= 4\,\bigg(1 + (p+q)\,C_1(\eta) + 4\sqrt{n}\,C_1(\eta)\,C_2(\eta)\bigg).\end{eqnarray*}

Thus, we get: \begin{eqnarray}\label{eqn:vanishing_harmonic_proof_1}|\langle\langle Au,\,u\rangle\rangle| \leq \frac{\sqrt{n}}{2}\,C_1(\eta)\,\bigg(\langle\langle\Delta'_\eta u,\,u\rangle\rangle + \langle\langle\Delta''_\eta u,\,u\rangle\rangle + C_{1,\,2}(p,\,q,\,n)\,||u||^2\bigg)\end{eqnarray} for every $(p,\,q)$-form $u$ on $X$ in the relevant domains.

Note that we have used the {\it completeness} of the metric $\omega=\gamma_\eta$ (hypothesis (ii)) in order to have \begin{eqnarray*}||\overline{D_\eta^{0,\,1}}u||^2 + ||\overline{D^{0,\,1}_\eta}^{\,\star}u||^2 = \langle\langle\Delta'_\eta u,\,u\rangle\rangle \hspace{5ex}\mbox{and}\hspace{5ex} ||D_\eta^{0,\,1}u||^2 + ||({D^{0,\,1}_\eta})^{\,\star} u||^2 = \langle\langle\Delta''_\eta u,\,u\rangle\rangle\end{eqnarray*} for every form $u\in\mbox{Dom}\,\Delta'_\eta\cap\mbox{Dom}\,\Delta''_\eta$.

\vspace{1ex}

$\bullet$ {\it Estimating the last but one term on the r.h.s. of (\ref{eqn:eta-BKN_1-0_0-1_simplified})}.

\vspace{1ex}

Recall that $\tau^{0,\,1}_\eta = 0$. Thus, w.r.t. the $L^2_\omega$-inner product, for every $(p,\,q)$-form $u$ we get: \begin{eqnarray*}\bigg|\bigg\langle\bigg\langle n\,\bigg[\overline{D^{0,\,1}_\eta} + \overline{\tau^{0,\,1}_\eta},\,\frac{1}{\eta}\,(\partial\eta\wedge\cdot\,)^\star\bigg] u,\,u\bigg\rangle\bigg\rangle\bigg| & \leq & n\,\bigg|\bigg\langle\bigg\langle\overline{D^{0,\,1}_\eta}\,u,\,\frac{1}{\eta}\,\partial\eta\wedge u\bigg\rangle\bigg\rangle\bigg| + n\,\bigg|\bigg\langle\bigg\langle\frac{1}{\eta}\,(\partial\eta\wedge\cdot\,)^\star u,\, \overline{D^{0,\,1}_\eta}^\star\,u\bigg\rangle\bigg\rangle\bigg| \\
  & \leq & n\,C_1(\eta)\,\bigg(\bigg|\bigg|\overline{D^{0,\,1}_\eta}\,u\bigg|\bigg|\,||u|| + \bigg|\bigg|\overline{D^{0,\,1}_\eta}^\star\,u\bigg|\bigg|\,||u||\bigg) \\
& \leq & \frac{n}{2}\,C_1(\eta)\,\bigg(\bigg|\bigg|\overline{D^{0,\,1}_\eta}\,u\bigg|\bigg|^2 + \bigg|\bigg|\overline{D^{0,\,1}_\eta}^\star\,u\bigg|\bigg|^2 + 2\,||u||^2\bigg).\end{eqnarray*}

After using the completeness of $\omega=\gamma_\eta$, we get: \begin{eqnarray}\label{eqn:vanishing_harmonic_proof_2}\bigg|\bigg\langle\bigg\langle n\,\bigg[\overline{D^{0,\,1}_\eta} + \overline{\tau^{0,\,1}_\eta},\,\frac{1}{\eta}\,(\partial\eta\wedge\cdot\,)^\star\bigg] u,\,u\bigg\rangle\bigg\rangle\bigg| \leq \frac{n}{2}\,C_1(\eta)\,\bigg(\bigg\langle\bigg\langle\Delta'_\eta u,\,u\bigg\rangle\bigg\rangle  + 2\,||u||^2\bigg)\end{eqnarray} for every $(p,\,q)$-form $u$ on $X$ in the relevant domains.

\vspace{1ex}

$\bullet$ {\it Estimating the last term on the r.h.s. of (\ref{eqn:eta-BKN_1-0_0-1_simplified})}.

\vspace{1ex}

W.r.t. the $L^2_\omega$-inner product, for every $(p,\,q)$-form $u$ we get: \begin{eqnarray*}\bigg|\bigg\langle\bigg\langle n\,\bigg[D^{0,\,1}_\eta,\,\frac{1}{\eta}\,(\bar\partial\eta\wedge\cdot\,)^\star\bigg]\,u,\,u\bigg\rangle\bigg\rangle\bigg| & \leq &  n\,\bigg|\bigg\langle\bigg\langle D^{0,\,1}_\eta\,u,\,\frac{1}{\eta}\,\bar\partial\eta\wedge u\bigg\rangle\bigg\rangle\bigg| + n\,\bigg|\bigg\langle\bigg\langle\frac{1}{\eta}\,(\bar\partial\eta\wedge\cdot\,)^\star u,\, (D^{0,\,1}_\eta)^\star\,u\bigg\rangle\bigg\rangle\bigg| \\
  & \leq & n\,C_1(\eta)\,\bigg(\bigg|\bigg|D^{0,\,1}_\eta\,u\bigg|\bigg|\,||u|| + \bigg|\bigg|(D^{0,\,1}_\eta)^\star\,u\bigg|\bigg|\,||u||\bigg) \\
  & \leq & \frac{n}{2}\,C_1(\eta)\,\bigg(\bigg|\bigg|D^{0,\,1}_\eta\,u\bigg|\bigg|^2 + \bigg|\bigg|(D^{0,\,1}_\eta)^\star u\bigg|\bigg|^2 + 2\,||u||^2\bigg).\end{eqnarray*}

After using the completeness of $\omega=\gamma_\eta$, we get: \begin{eqnarray}\label{eqn:vanishing_harmonic_proof_3}\bigg|\bigg\langle\bigg\langle n\,\bigg[D^{0,\,1}_\eta,\,\frac{1}{\eta}\,(\bar\partial\eta\wedge\cdot\,)^\star\bigg] u,\,u\bigg\rangle\bigg\rangle\bigg| \leq \frac{n}{2}\,C_1(\eta)\,\bigg(\langle\langle\Delta''_\eta u,\,u\rangle\rangle  + 2\,||u||^2\bigg)\end{eqnarray} for every $(p,\,q)$-form $u$ on $X$ in the relevant domains.

\vspace{1ex}

$\bullet$ Putting together (\ref{eqn:eta-BKN_1-0_0-1_simplified})--(\ref{eqn:vanishing_harmonic_proof_3}), we get: \begin{eqnarray*}\langle\langle\Delta''_\eta u,\,u\rangle\rangle & \geq & \langle\langle\Delta'_\eta u,\,u\rangle\rangle + (p-q)(p+q-n)\,||u||^2 \\
   & - & \frac{\sqrt{n}}{2}\,C_1(\eta)\,\bigg(\langle\langle\Delta'_\eta u,\,u\rangle\rangle + \langle\langle\Delta''_\eta u,\,u\rangle\rangle + C_{1,\,2}(p,\,q,\,n)\,||u||^2\bigg) \\
  & - & \frac{n}{2}\,C_1(\eta)\,\bigg(\langle\langle\Delta'_\eta u,\,u\rangle\rangle  + 2\,||u||^2\bigg) - \frac{n}{2}\,C_1(\eta)\,\bigg(\langle\langle\Delta''_\eta u,\,u\rangle\rangle  + 2\,||u||^2\bigg),\end{eqnarray*} which amounts to \begin{eqnarray}\label{eqn:eta-BKN_1-0_0-1_simplified_ineq}\bigg(1 + \frac{n+\sqrt{n}}{2}\,C_1(\eta)\bigg)\,\langle\langle\Delta''_\eta u,\,u\rangle\rangle & \geq & \bigg(1 - \frac{n+\sqrt{n}}{2}\,C_1(\eta)\bigg)\,\langle\langle\Delta'_\eta u,\,u\rangle\rangle \\
  \nonumber  & + & \bigg[(p-q)(p+q-n) - C_1(\eta)\,\bigg(2n + \frac{\sqrt{n}}{2}\,C_{1,\,2}(p,\,q,\,n)\bigg)\bigg]\,||u||^2\end{eqnarray} for every $(p,\,q)$-form $u$ on $X$ in the relevant domains.

  Since $\langle\langle\Delta''_\eta u,\,u\rangle\rangle, \langle\langle\Delta'_\eta u,\,u\rangle\rangle\geq 0$, (\ref{eqn:eta-BKN_1-0_0-1_simplified_ineq}) shows that the vanishing of $\langle\langle\Delta''_\eta u,\,u\rangle\rangle$ (which is equivalent to $u\in\ker\Delta''_\eta$) implies the vanishing of $\langle\langle\Delta'_\eta u,\,u\rangle\rangle$ (which is irrelevant to us here) and the vanishing of $u$ whenever the coefficients of $\langle\langle\Delta'_\eta u,\,u\rangle\rangle$ and $||u||^2$ are positive. For this to happen, we need \begin{eqnarray*}\label{eqn:vanishing_cond_1}(a)\hspace{2ex} C_1(\eta)<\frac{2}{n+\sqrt{n}}   \hspace{5ex}\mbox{and}\hspace{5ex} (b)\hspace{2ex} C_1(\eta)\,\bigg(2n + \frac{\sqrt{n}}{2}\,C_{1,\,2}(p,\,q,\,n)\bigg)< (p-q)(p+q-n).\end{eqnarray*}

  Since $p,q$ are integers and we require that either ($p>q$ and $p+q\geq n+1$) or ($p<q$ and $p+q\leq n-1$), the r.h.s. of (b) is a positive integer. Hence, (b) is satisfied whenever the inequality \begin{eqnarray*}(b')\hspace{2ex} C_1(\eta)\,\bigg(2n + \frac{\sqrt{n}}{2}\,C_{1,\,2}(p,\,q,\,n)\bigg)< 1\end{eqnarray*} is satisfied. Meanwhile, $p+q\leq 2n$, so \begin{eqnarray*}C_{1,\,2}(p,\,q,\,n) \leq 4\bigg(1 + 2n\,C_1(\eta) + 4\sqrt{n}\,C_1(\eta)\,C_2(\eta)\bigg).\end{eqnarray*} Thus, for (b') to hold, it suffices to have \begin{eqnarray*}C_1(\eta)\,\bigg[2n + \frac{\sqrt{n}}{2}\,4\,\bigg(1 + 2n\,C_1(\eta) + 4\sqrt{n}\,C_1(\eta)\,C_2(\eta)\bigg)\bigg] <1,\end{eqnarray*} which is equivalent to \begin{eqnarray*}C_1(\eta)^2\,4n\,\bigg(\sqrt{n} + 2\,C_2(\eta)\bigg) + 2n\,C_1(\eta) <1.\end{eqnarray*} Since $C_1(\eta)^2<C_1(\eta)$ (indeed, (a) implies $C_1(\eta)<1$), for this to happen it suffices that \begin{eqnarray*}(b'')\hspace{2ex} C_1(\eta)\,\bigg[4n\,\bigg(\sqrt{n} + 2\,C_2(\eta)\bigg) + 2n\bigg] <1.\end{eqnarray*}

    Now, recall that $\varphi$ is the function such that $\eta = e^{-\varphi}$. Thus, $C_1(\eta) = \sup\limits_X(|\partial\eta|/\eta) = \sup\limits_X|\partial\varphi|$ and \begin{eqnarray*}C_2(\eta) = \sup\limits_X\frac{|\partial\bar\partial\eta|}{\eta} = \sup\limits_X|\partial\varphi\wedge\bar\partial\varphi - \partial\bar\partial\varphi| \leq \sup\limits_X|\partial\varphi|^2 + \sup\limits_X|\partial\bar\partial\varphi| \leq C_1(\eta)^2 + C(\varphi).\end{eqnarray*} This upper estimate for $C_2(\eta)$ shows that (b'') holds whenever \begin{eqnarray*}2nC_1(\eta)\,\bigg[1 + 2\,\bigg(\sqrt{n} + 2\,C_1(\eta)^2 + 2\,C(\varphi)\bigg)\bigg]<1.\end{eqnarray*} This is equivalent to \begin{eqnarray*}2n\,C_1(\eta) + 4n\sqrt{n}\,C_1(\eta) + 8n\,C_1(\eta)^3 + 8n\,C_1(\eta)\,C(\varphi)<1.\end{eqnarray*} Since $C_1(\eta)^3<C_1(\eta)$, for this to happen it suffices that \begin{eqnarray*}C_1(\eta)\,\bigg(10n + 4n\sqrt{n} + 8n\,C(\varphi)\bigg) <1.\end{eqnarray*}

    This holds thanks to hypothesis (\ref{eqn:vanishing_harmonic_H1}), which also implies (a). The proof of Theorem \ref{The:vanishing_harmonic} is complete.  \hfill $\Box$

\vspace{2ex}

 We now discuss two variants of Theorem \ref{The:vanishing_harmonic}. As noticed in (\ref{eqn:eta-derivatives_phi}), if we write $\eta=e^{-\varphi}$ we have 
 \begin{eqnarray*}\gamma_\eta = \frac{1}{\eta^2}\,i\partial\eta\wedge\bar\partial\eta + \bigg(\frac{1}{\eta^2}\,i\partial\eta\wedge\bar\partial\eta - \frac{1}{\eta}\,i\partial\bar\partial\eta\bigg)  =  i\partial\varphi\wedge\bar\partial\varphi + i\partial\bar\partial\varphi \geq i\partial\bar\partial\varphi,\end{eqnarray*} the last inequality being a consequence of the standard inequality $i\alpha\wedge\overline\alpha\geq 0$ for any $(1,0)$-form $\alpha$. In particular, if we assume $i\partial\bar\partial\varphi > 0$ on $X$, we also have $\gamma_\eta >0$ on $X$. Moreover, if we further assume the metric $i\partial\bar\partial\varphi$ to be complete on $X$, the metric $\gamma_\eta$ is complete as well.

 Thus, if $i\partial\bar\partial\varphi$ is supposed to be a complete metric on $X$, hypotheses (i) and (ii) of Theorem \ref{The:vanishing_harmonic} are satisfied, while $|\partial\eta|_{\gamma_\eta}\leq|\partial\eta|_{i\partial\bar\partial\varphi}$ and \begin{eqnarray*}|i\partial\bar\partial\varphi|_{\gamma_\eta}\leq|i\partial\bar\partial\varphi|_{i\partial\bar\partial\varphi} = \sqrt{n}\end{eqnarray*} at every point of $X$. In particular, $C(\varphi) = \sup\limits_X|i\partial\bar\partial\varphi|_{\gamma_\eta}\leq\sqrt{n}$, so we get the following consequence of Theorem \ref{The:vanishing_harmonic}.

 \begin{Cor}\label{Cor:vanishing_harmonic_light} Let $X$ be a (non-compact) complex manifold with $\mbox{dim}_\C X=n$. Suppose there exists a $C^\infty$ function $\eta=e^{-\varphi}:X\longrightarrow(0,\,\infty)$ such that the $(1,\,1)$-form $i\partial\bar\partial\varphi$ is {\bf positive definite} at every point and the Hermitian metric it defines on $X$ is {\bf complete}. Suppose, moreover, that \begin{eqnarray}\label{eqn:vanishing_harmonic_light_H1}\sup\limits_X|\partial\varphi|_{i\partial\bar\partial\varphi} < \frac{1}{10n + 12n\,\sqrt{n}}.\end{eqnarray}

   Then, for any bidegree $(p,\,q)$ such that either $\bigg(p>q \hspace{1ex}\mbox{and}\hspace{1ex} p+q\geq n+1\bigg)$ or $\bigg(p<q \hspace{1ex}\mbox{and}\hspace{1ex} p+q\leq n-1\bigg)$, the space of $\Delta''_\eta$-harmonic $L^2_{\gamma_\eta}$-forms on $X$ of bidegree $(p,\,q)$ {\bf vanishes}: \begin{eqnarray}\label{eqn:vanishing_harmonic_light_conclusion}{\cal H}^{p,\,q}_{\Delta''_\eta}(X,\,\C):=\ker\bigg(\Delta''_\eta:\mbox{Dom}_{p,\,q}(\Delta''_\eta)\longrightarrow L^2_{p,\,q}(X,\,\C)\bigg) = \{0\},\end{eqnarray} where $\Delta''_\eta$ is the closed and densely defined unbounded extension to the space of $L^2_{\gamma_\eta}$-forms of bidegree $(p,\,q)$ of the operator $\Delta''_\eta$ defined on $C^\infty$ forms w.r.t. the metric $\omega=\gamma_\eta = i\partial\varphi\wedge\bar\partial\varphi + i\partial\bar\partial\varphi$.

 \end{Cor}

 The next observation is that the upper bound in hypothesis (\ref{eqn:vanishing_harmonic_light_H1}) can be made as small as we wish after possibly replacing $\varphi$ with $\varepsilon\varphi$ and choosing the constant $\varepsilon>0$ small enough. Indeed, \begin{eqnarray*}|\partial(\varepsilon\varphi)|_{i\partial\bar\partial(\varepsilon\varphi)} = \varepsilon\,|\partial\varphi|_{\varepsilon\,i\partial\bar\partial\varphi} = \frac{\varepsilon}{\sqrt{\varepsilon}}\,|\partial\varphi|_{i\partial\bar\partial\varphi} = \sqrt{\varepsilon}\,|\partial\varphi|_{i\partial\bar\partial\varphi}.\end{eqnarray*} Thus, the following consequence of Corollary \ref{Cor:vanishing_harmonic_light} is a manifestation of the {\it twisted adiabatic limit} (as $\varepsilon\downarrow 0$ in this case) introduced in this paper.

 \begin{Cor}\label{Cor:vanishing_harmonic_family} Let $X$ be a (non-compact) complex manifold with $\mbox{dim}_\C X=n$. Suppose there exists a $C^\infty$ function $\eta=e^{-\varphi}:X\longrightarrow(0,\,\infty)$ such that the $(1,\,1)$-form $i\partial\bar\partial\varphi$ is {\bf positive definite} at every point and the Hermitian metric it defines on $X$ is {\bf complete}. Suppose, moreover, that \begin{eqnarray}\label{eqn:vanishing_harmonic_family_H1}\sup\limits_X|\partial\varphi|_{i\partial\bar\partial\varphi} < \infty.\end{eqnarray}

   For every constant $\varepsilon>0$, let $\eta_\varepsilon = e^{-\varepsilon\varphi}$.

   Then, for any bidegree $(p,\,q)$ such that either $\bigg(p>q \hspace{1ex}\mbox{and}\hspace{1ex} p+q\geq n+1\bigg)$ or $\bigg(p<q \hspace{1ex}\mbox{and}\hspace{1ex} p+q\leq n-1\bigg)$ and for every $\varepsilon>0$ {\bf small enough}, the space of $\Delta''_{\eta_\varepsilon}$-harmonic $L^2_{\gamma_{\eta_\varepsilon}}$-forms on $X$ of bidegree $(p,\,q)$ {\bf vanishes}: \begin{eqnarray}\label{eqn:vanishing_harmonic_family_conclusion}{\cal H}^{p,\,q}_{\Delta''_{\eta_\varepsilon}}(X,\,\C):=\ker\bigg(\Delta''_{\eta_\varepsilon}:\mbox{Dom}_{p,\,q}(\Delta''_{\eta_\varepsilon})\longrightarrow L^2_{p,\,q}(X,\,\C)\bigg) = \{0\},\end{eqnarray} where $\Delta''_{\eta_\varepsilon}$ is the closed and densely defined unbounded extension to the space of $L^2_{\gamma_{\eta_\varepsilon}}$-forms of bidegree $(p,\,q)$ of the operator $\Delta''_{\eta_\varepsilon}$ defined on $C^\infty$ forms w.r.t. the metric $\omega_\varepsilon=\gamma_{\eta_\varepsilon} = \varepsilon^2\,i\partial\varphi\wedge\bar\partial\varphi + \varepsilon\,i\partial\bar\partial\varphi$.

 \end{Cor}

 \noindent {\it Proof.} The positive definiteness and the completeness of $i\partial\varphi\wedge\bar\partial\varphi$ are preserved when $\varphi$ is replaced by $\varepsilon\varphi$. Thus, the claim follows at once from Corollary \ref{Cor:vanishing_harmonic_light} and from the short discussion that preceded the statement.  \hfill $\Box$

\section{Vanishing of certain Dolbeault cohomology groups on certain compact complex manifolds}\label{section:vanishing_cohomology_compact} The discussion in this section is the analogue in the compact setting of the one we had in $\S$\ref{section:vanishing_harmonic}.

      Let $(X,\,\omega)$ be a {\bf compact} Hermitian manifold with $\mbox{dim}_\C X=n$. Fix an arbitrary $C^\infty$ function $\eta:X\longrightarrow(0,\,\infty)$. For any $(p,\,q)$, we consider the curvature-like ($1$-st order differential) operator \begin{eqnarray*}\label{eqn:F_eta-operator}F_\eta:=i\,\bigg[[D^{0,\,1}_\eta,\,\overline{D^{0,\,1}_\eta}],\,\Lambda\bigg]:C^\infty_{p,\,q}(X,\,\C)\longrightarrow C^\infty_{p,\,q}(X,\,\C).\end{eqnarray*}

      We will use the following standard terminology: if $A$ and $B$ are linear operators acting on the differential forms on $X$, we will say that $A\geq B$ (resp. $A>B$) if $\langle\langle Au,\,u\rangle\rangle\geq\langle\langle Bu,\,u\rangle\rangle$ for every $u$ (resp. if $\langle\langle Au,\,u\rangle\rangle>\langle\langle Bu,\,u\rangle\rangle$ for every $u\neq 0$).

\vspace{2ex}

 As in $\S$\ref{section:vanishing_harmonic}, we consider the constant $C_1(\eta):=\sup\limits_X\frac{|\partial\eta|}{\eta}<\infty$. Unlike in $\S$\ref{section:vanishing_harmonic}, the pointwise norm $|\partial\eta|$, like all the norms and all the inner products in this section, is the one induced by $\omega$ and, due to the compactness of $X$,  the constant $C_1(\eta)$ is finite. We start by noticing that the pointwise operator norm of $\tau_\eta^{0,\,1} = [\Lambda,\,D_\eta^{0,\,1}\omega\wedge\cdot\,]$ w.r.t. $\omega$ satisfies \begin{eqnarray*}\bigg|\tau_\eta^{0,\,1}\bigg| = \sup_{|u|=1}\bigg|\tau_\eta^{0,\,1} u\bigg| \leq 2\sqrt{n}\,\bigg|D_\eta^{0,\,1}\omega\bigg| =  2\sqrt{n}\,\bigg|\bar\partial\omega - \frac{1}{\eta}\,\bar\partial\eta\wedge\omega\bigg| \leq 2\sqrt{n}\,\bigg(|\bar\partial\omega| + C_1(\eta)\,\sqrt{n}\bigg):=C_3(\eta).\end{eqnarray*}

\begin{The}\label{The:vanishing_cohomology_compact} Let $(X,\,\omega)$ be a compact complex Hermitian manifold with $\mbox{dim}_\C X=n$. If there exists a $C^\infty$ function $\eta:X\longrightarrow(0,\,\infty)$ such that:

\vspace{1ex}

(i)\, $C_3(\eta) + n\sqrt{n}\,C_1(\eta)<2$;

\vspace{1ex}

(ii)\, for some bidegree $(p,\,q)$, $F_\eta>2\,\bigg(C_3(\eta) + n\sqrt{n}\,C_1(\eta)\bigg)\,\mbox{Id}$ in bidegree $(p,\,q)$,

\vspace{1ex}

\noindent then $H^{p,\,q}_{\bar\partial}(X,\,\C) = \{0\}$.

\end{The}

\noindent {\it Proof.} The rough $\eta$-BKN identity (\ref{eqn:eta-BKN_rough_1-0_0-1}) will suffice for our purposes. We will estimate the last four terms on its r.h.s. by using the Cauchy-Schwarz inequality and the elementary inequality $ab\leq (a^2 + b^2)/2$ for non-negative reals $a, b$. 

Let $u\in C^\infty_{p,\,q}(X,\,\C)$. For the $L^2$-inner products w.r.t. $\omega$, we get: \begin{eqnarray*}\bigg|\bigg\langle\bigg\langle\bigg[\overline{D^{0,\,1}_\eta},\,\overline{\tau^{0,\,1}_\eta}^{\,\star}\bigg]\,u,\,u\bigg\rangle\bigg\rangle\bigg| & \leq & \bigg|\bigg|\overline{D^{0,\,1}_\eta}\,u\bigg|\bigg|\,\bigg|\bigg|\overline{\tau^{0,\,1}_\eta}\,u\bigg|\bigg| + \bigg|\bigg|\overline{D^{0,\,1}_\eta}^\star\,u\bigg|\bigg|\,\bigg|\bigg|\overline{\tau^{0,\,1}_\eta}^\star\,u\bigg|\bigg| \\
  & \leq & C_3(\eta)\,\bigg(\bigg|\bigg|\overline{D^{0,\,1}_\eta}\,u\bigg|\bigg|\,||u|| + \bigg|\bigg|\overline{D^{0,\,1}_\eta}^\star\,u\bigg|\bigg|\,||u||\bigg) \\
  & \leq & \frac{C_3(\eta)}{2}\,\bigg(\bigg|\bigg|\overline{D^{0,\,1}_\eta}\,u\bigg|\bigg|^2 + \bigg|\bigg|\overline{D^{0,\,1}_\eta}^\star\,u\bigg|\bigg|^2 + 2||u||^2\bigg).\end{eqnarray*}

Similarly, we get: \begin{eqnarray*}\bigg|\bigg\langle\bigg\langle\bigg[D^{0,\,1}_\eta,\,(\tau^{0,\,1}_\eta)^{\,\star}\bigg]\,u,\,u\bigg\rangle\bigg\rangle\bigg| & \leq & \frac{C_3(\eta)}{2}\,\bigg(\bigg|\bigg|D^{0,\,1}_\eta\,u\bigg|\bigg|^2 + \bigg|\bigg|(D^{0,\,1}_\eta)^\star\,u\bigg|\bigg|^2 + 2\,||u||^2\bigg), \\
  \bigg|\bigg\langle\bigg\langle n\,\bigg[D^{0,\,1}_\eta,\,\frac{i}{\eta}\,[\partial\eta\wedge\cdot\,,\,\Lambda]\bigg]\,u,\,u\bigg\rangle\bigg\rangle\bigg| & \leq & \frac{n\sqrt{n}\,C_1(\eta)}{2}\,\bigg(\bigg|\bigg|D^{0,\,1}_\eta\,u\bigg|\bigg|^2 + \bigg|\bigg|(D^{0,\,1}_\eta)^\star\,u\bigg|\bigg|^2 + 2\,||u||^2\bigg), \\
  \bigg|\bigg\langle\bigg\langle n\,\bigg[\overline{D^{0,\,1}_\eta},\,\frac{i}{\eta}\,[\bar\partial\eta\wedge\cdot\,,\,\Lambda]\bigg]\,u,\,u\bigg\rangle\bigg\rangle\bigg| & \leq & \frac{n\sqrt{n}\,C_1(\eta)}{2}\, \bigg(\bigg|\bigg|\overline{D^{0,\,1}_\eta}\,u\bigg|\bigg|^2 + \bigg|\bigg|\overline{D^{0,\,1}_\eta}^\star\,u\bigg|\bigg|^2 + 2\,||u||^2\bigg).\end{eqnarray*}

From these four inequalities and the rough $\eta$-BKN identity (\ref{eqn:eta-BKN_rough_1-0_0-1}), we get: \begin{eqnarray*}\langle\langle\Delta''_\eta u,\,u\rangle\rangle \geq \langle\langle\Delta'_\eta u,\,u\rangle\rangle & + & \langle\langle F_\eta u,\,u\rangle\rangle - \frac{C_3(\eta)}{2}\,\bigg(\langle\langle\Delta'_\eta u,\,u\rangle\rangle + 2\,||u||^2\bigg) - \frac{C_3(\eta)}{2}\,\bigg(\langle\langle\Delta''_\eta u,\,u\rangle\rangle + 2\,||u||^2\bigg) \\
  & - & \frac{n\sqrt{n}\,C_1(\eta)}{2}\,\bigg(\langle\langle\Delta''_\eta u,\,u\rangle\rangle + 2\,||u||^2\bigg) - \frac{n\sqrt{n}\,C_1(\eta)}{2}\,\bigg(\langle\langle\Delta'_\eta u,\,u\rangle\rangle + 2\,||u||^2\bigg)\end{eqnarray*} for every $u\in C^\infty_{p,\,q}(X,\,\C)$.

This amounts to \begin{eqnarray*}\bigg(1 + \frac{C_3(\eta) + n\sqrt{n}\,C_1(\eta)}{2}\bigg)\,\langle\langle\Delta''_\eta u,\,u\rangle\rangle & \geq & \bigg(1 - \frac{C_3(\eta) + n\sqrt{n}\,C_1(\eta)}{2}\bigg)\,\langle\langle\Delta'_\eta u,\,u\rangle\rangle + \langle\langle F_\eta u,\,u\rangle\rangle \\
  & - & 2\,\bigg(C_3(\eta) + n\sqrt{n}\,C_1(\eta)\bigg)\,||u||^2\end{eqnarray*} for every $u\in C^\infty_{p,\,q}(X,\,\C)$.

  This inequality, together with the hypotheses (i) and (ii), implies that whenever $\Delta''_\eta u = 0$ we must have $u=0$. The result follows from this, from the Hodge isomorphism (\ref{eqn:Hodge_isom_Delta''-Delta'_eta}) and from the cohomology isomorphism (\ref{eqn:cohomologies_1-0_0-1_isom}).  \hfill $\Box$

\section{Appendix: review of standard commutation relations}\label{Appendix}

 We briefly recall here some standard formulae that were used throughout the paper.

 \begin{Lem}\label{Lem:com_1} Let $(X,\,\omega)$ be a compact complex Hermitian manifold. The following standard Hermitian commutation relations ([Dem84], see also [Dem97, VII, $\S.1$]) hold: \begin{eqnarray}\label{eqn:standard-comm-rel}\nonumber &  & (i)\,\,(\partial + \tau)^{\star} = i\,[\Lambda,\,\bar\partial];  \hspace{3ex} (ii)\,\,(\bar\partial + \bar\tau)^{\star} = - i\,[\Lambda,\,\partial]; \\
&  & (iii)\,\, \partial + \tau = -i\,[\bar\partial^{\star},\,L]; \hspace{3ex} (iv)\,\,
\bar\partial + \bar\tau = i\,[\partial^{\star},\,L],\end{eqnarray}

   \noindent where the upper symbol $\star$ stands for the formal adjoint w.r.t. the $L^2$ inner product induced by $\omega$, $L=L_{\omega}:=\omega\wedge\cdot$ is the Lefschetz operator of multiplication by $\omega$, $\Lambda=\Lambda_{\omega}:=L^{\star}$ and $\tau:=[\Lambda,\,\partial\omega\wedge\cdot]$ is the torsion operator (of order zero and type $(1,\,0)$) associated with the metric $\omega$.

 \end{Lem}  

 Again following [Dem97, VII, $\S.1$], recall that the commutation relations $(1)$ immediately induce via the Jacobi identity the Bochner-Kodaira-Nakano-type identity \begin{equation}\label{eqn:BKN_demailly1}\Delta'' = \Delta' + [\partial,\,\tau^{\star}] - [\bar\partial,\,\bar{\tau}^{\star}]\end{equation} 

\noindent relating the $\bar\partial$-Laplacian $\Delta''=[\bar\partial,\,\bar\partial^{\star}]=\bar\partial\bar\partial^{\star} + \bar\partial^{\star} \bar\partial$ and the $\partial$-Laplacian $\Delta'=[\partial,\,\partial^{\star}]=\partial\partial^{\star} + \partial^{\star}\partial$. This, in turn, induces the following Bochner-Kodaira-Nakano-type identity (cf. [Dem84]) in which the first-order terms have been absorbed in the twisted Laplace-type operator $\Delta'_{\tau}:=[\partial+\tau,\, (\partial+\tau)^{\star}]$: \begin{equation}\label{eqn:BKN_demailly2}\Delta'' = \Delta'_{\tau} + T_{\omega},\end{equation} 

\noindent where $T_{\omega}:=\bigg[\Lambda,\,[\Lambda,\,\frac{i}{2}\,\partial\bar\partial\omega]\bigg] - [\partial\omega\wedge\cdot,\,(\partial\omega\wedge\cdot)^{\star}]$ is a zeroth order operator of type $(0,\,0)$ associated with the torsion of $\omega$. Formula (\ref{eqn:BKN_demailly2}) is obtained from (\ref{eqn:BKN_demailly1}) via the following identities (cf. [Dem84] or [Dem97, VII, $\S.1$]) which have an interest of their own: \begin{eqnarray}\label{eqn:BKN_demailly_auxiliary}\nonumber &  & (i)\,\,[L,\,\tau] = 3\,\partial\omega\wedge\cdot,  \hspace{3ex} (ii)\,\, [\Lambda,\,\tau] = -2i\,\bar{\tau}^{\star},\\
 &  & (iii)\,\, [\partial,\,\bar{\tau}^{\star}] = - [\partial,\,\bar\partial^{\star}] = [\tau,\,\bar\partial^{\star}],  \hspace{3ex}  (iv)\,\, -[\bar\partial,\,\bar\tau^{\star}] = [\tau,\, (\partial+\tau)^{\star}] + T_{\omega}.\end{eqnarray}

\noindent Note that $(iii)$ yields, in particular, that $\partial$ and $\bar\partial^{\star} + \bar\tau^{\star}$ anti-commute, hence by conjugation, $\bar\partial$ and $\partial^{\star} + \tau^{\star}$ anti-commute, i.e. \begin{equation}\label{eqn:basic-anticommutation}[\partial,\,\bar\partial^{\star} + \bar\tau^{\star}] = 0 \hspace{2ex} \mbox{and} \hspace{2ex} [\bar\partial,\,\partial^{\star} + \tau^{\star}] = 0.\end{equation}

\vspace{2ex}

The following formulae can be viewed as commutation relations for zeroth-order operators (see [Pop03, $\S.1.0.2$] or [Pop23, Appendix]). 

\begin{Lem}\label{Lem:com_2} Let $(X,\,\omega)$ be a complex Hermitian manifold and $\eta$ a {\bf real-valued} ${\cal C}^\infty$ function on $X$. 

The following identities hold pointwise for arbitrary differential forms of any degree on $X$. \\

\hspace{2ex} (a)\,\,$[\partial\eta\wedge\cdot, \, \Lambda]=i\, (\bar\partial\eta\wedge\cdot)^{\star}$, \hspace{4ex}  $[\bar\partial\eta\wedge\cdot, \, \Lambda] = -i\, (\partial\eta\wedge\cdot)^{\star}$ \\

\hspace{6ex} $[L, \, (\partial\eta\wedge\cdot)^{\star}] = -i\, \bar\partial\eta\wedge\cdot$,  \hspace{4ex} $[(\bar\partial\eta\wedge\cdot)^{\star}, \, L] = -i\, \partial\eta\wedge\cdot$. 

\noindent \begin{eqnarray}\nonumber (b)\,\, [i\, \partial\eta\wedge\bar\partial\eta\wedge\cdot, \, \Lambda] & = & (\partial\eta\wedge\cdot) \, (\partial\eta\wedge\cdot)^{\star} - (\bar\partial\eta\wedge\cdot)^{\star}\, (\bar\partial\eta\wedge\cdot) \\
  \nonumber & = & (\bar\partial\eta\wedge\cdot)\, (\bar\partial\eta\wedge\cdot)^{\star} - (\partial\eta\wedge\cdot)^{\star}  \, (\partial\eta\wedge\cdot).\end{eqnarray}

\end{Lem}

\begin{Cor}\label{Lem:com_d-eta} In the setting of Lemma \ref{Lem:com_2}, the next pointwise identities hold: \\

  \hspace{2ex} (a)\,\,$[\Lambda,\,id_\eta\eta\wedge\cdot\,]=-(\bar{d}_{-\eta}\eta\wedge\cdot)^{\star}$, \hspace{4ex}  (b)\,\, $[\Lambda,\,i\bar{d}_{-\eta}\eta\wedge\cdot\,] = (d_\eta\eta\wedge\cdot)^{\star}$.

\end{Cor}  

\noindent {\it Proof.} This follows right away from Lemma \ref{Lem:com_2} after we use the identities: \begin{eqnarray*}d_\eta\eta = \eta\,\partial\eta + \bar\partial\eta  \hspace{6ex}\mbox{and}\hspace{6ex} \bar{d}_{-\eta}\eta = \partial\eta - \eta\,\bar\partial\eta.\end{eqnarray*} \hfill $\Box$ 

\vspace{6ex}

\noindent {\bf References} \\

\noindent [BP18]\, H. Bellitir, D. Popovici --- {\it Positivity Cones under Deformations of Complex Structures} --- Riv. Mat. Univ. Parma, Vol. {\bf 9} (2018), 133-176.

\vspace{1ex}

\noindent [COUV16]\, M. Ceballos, A. Otal, L. Ugarte, R. Villacampa --- {\it Invariant complex structures on $6$-nilmanifolds: classification, Fr\"olicher spectral sequence and special Hermitian metrics} --- J. Geom. Anal. {\bf 26} (2016), no. 1, 252–286. 

\vspace{1ex}

\noindent [Dem84]\, J.-P. Demailly --- {\it Sur l'identit\'e de Bochner-Kodaira-Nakano en g\'eom\'etrie hermitienne} --- S\'eminaire d'analyse P. Lelong, P. Dolbeault, H. Skoda (editors) 1983/1984, Lecture Notes in Math., no. {\bf 1198}, Springer Verlag (1986), 88-97.

\vspace{1ex}

\noindent [Dem 96]\, J.-P. Demailly --- {\it Th\'eorie de Hodge $L^2$ et th\'eor\`emes d'annulation} --- in ``Introduction \`a la th\'eorie de Hodge'', J. Bertin, J.-P. Demailly, L. Illusie, C. Peters, Panoramas et Synth\`eses {\bf 3}, SMF (1996).

\vspace{1ex}

\noindent [Dem 97]\, J.-P. Demailly --- {\it Complex Analytic and Algebraic Geometry}---http://www-fourier.ujf-grenoble.fr/~demailly/books.html

\vspace{1ex}

\noindent [Gri69]\, P. Griffiths --- {\it Hermitian Differential Geometry, Chern Classes and Positive Vector Bundles} --- Global Analysis, papers in honor of K. Kodaira, Univ. of Tokyo Press, Tokyo (1969) 185-251.

\vspace{1ex}

\noindent [Gro91]\, M. Gromov --- {\it K\"ahler Hyperbolicity and $L^2$ Hodge Theory} --- J. Diff. Geom. {\bf 33} (1991), 263-292.

\vspace{1ex}

\noindent [KP23]\, H. Kasuya, D. Popovici --- {\it Partially Hyperbolic Compact Complex Manifolds} --- arXiv e-print DG 2304.01697v1.

\vspace{1ex}

\noindent [MM90]\, R. R. Mazzeo, R. B. Melrose --- {\it The Adiabatic Limit, Hodge Cohomology and Leray's Spectral Sequence} --- J. Diff. Geom. {\bf 31} (1990) 185-213.

\vspace{1ex}

\noindent [MP21]\, S. Marouani, D. Popovici --- {\it Balanced Hyperbolic and Divisorially Hyperbolic Compact Complex Manifolds} --- Math. Res. Lett., Vol. {\bf 30}, No. 6, 1813-1855 (2023),

\noindent DOI: https://dx.doi.org/10.4310/MRL.2023.v30.n6.a7

\vspace{1ex}

\noindent [MP21]\, S. Marouani, D. Popovici --- {\it Some Properties of Balanced Hyperbolic Compact Complex Manifolds} --- Internat. J. Math. {\bf 33} (2022), no. 3, 2250019 (39 pages), DOI: 10.1142/S0129167X22500197.

\vspace{1ex}

\noindent [Ohs82]\, T. Ohsawa --- {\it Isomorphism Theorems for Cohomology Groups of Weakly $1$-Complete Manifolds} --- Publ. Res. Inst. Math. Sci. {\bf 18} (1982) 191-232.

\vspace{1ex}

\noindent [Pop03]\, D. Popovici --- {\it Quelques applications des m\'ethodes effectives en g\'eom\'etrie analytique} --- PhD thesis, University Joseph Fourier (Grenoble 1), http://tel.ccsd.cnrs.fr/documents/

\vspace{1ex}

\noindent [Pop18]\, D. Popovici --- {\it Non-K\"ahler Mirror Symmetry of the Iwasawa Manifold} --- Int. Math. Res. Not. IMRN 2020, no. {\bf 23}, 9471–9538.

\vspace{1ex}

\noindent [Pop19]\, D. Popovici --- {\it Adiabatic Limit and the Fr\"olicher Spectral Sequence} --- Pacific J. Math. {\bf 300} (2019), no. 1, 121-158.

\vspace{1ex}

\noindent [PSU21a]\, D. Popovici, J. Stelzig, L. Ugarte --- {\it Higher-Page Bott-Chern and Aeppli Cohomologies and Applications} --- J. reine angew. Math. (Crelle) {\bf 777} (2021), 157-194.

\vspace{1ex}

\noindent [PSU21b]\, D. Popovici, J. Stelzig, L. Ugarte --- {\it Deformations of Higher-Page Analogues of $\partial\bar\partial$-Manifolds} --- Math. Z. {\bf 300} (2022), no. 3, 2611–2635.

\vspace{1ex}

\noindent [PSU22]\, D. Popovici, J. Stelzig, L. Ugarte --- {\it Higher-Page Hodge Theory of Compact Complex Manifolds} --- Ann. Sc. Norm. Super. Pisa Cl. Sci. (5) Vol. {\bf XXV} (2024), 1431-1464.



\vspace{1ex}

\noindent [Wit85]\, E. Witten --- {\it Global Gravitational Anomalies} --- Commun. Math. Phys, {\bf 100}, 197-229 (1985).

\vspace{6ex}

\noindent Universit\'e Paul Sabatier, Institut de Math\'ematiques de Toulouse

\noindent 118, route de Narbonne, 31062, Toulouse Cedex 9, France

\noindent Email: popovici@math.univ-toulouse.fr

\end{document}